\documentclass[12pt,reqno]{amsart}
\usepackage{amssymb}
\usepackage{mathtools}
\usepackage[dvipsnames]{xcolor}
\usepackage{hyperref}
\usepackage{dsfont}
\usepackage{enumerate}
\usepackage{dsfont}
\usepackage[pdftex]{graphicx}
\usepackage{caption}
\usepackage[labelfont=]{subcaption}
\usepackage{marginnote}

\usepackage{tikz}
\usepackage{verbatim}
\usetikzlibrary{arrows}
\usetikzlibrary{calc,decorations.pathreplacing,positioning,shapes.multipart}

\newtheorem{theorem}{Theorem}[section]
\newtheorem{definition}[theorem]{Definition}

\newtheorem{lemma}[theorem]{Lemma}
\newtheorem{corollary}[theorem]{Corollary}

\newtheorem{example}[theorem]{Example}
\newtheorem{remark}[theorem]{Remark}
\newtheorem{assumption}[theorem]{Assumption}

\newtheorem*{proposition*}{Proposition}

\newcommand*\E{\mathbb{E}}
\renewcommand*\P{\mathbb{P}}
\newcommand{\R}{\mathbb{R}}  
\newcommand{\Z}{\mathbb{Z}}

\newcommand{\N}{\mathbb{N}}
\newcommand{\NN}{\mathcal{N}}

\newcommand*\as{\textnormal{(a.s.)}}

\newcommand*\todo[2][]{{\color{red} \ifx\\#1\\ TODO\else TODO(#1)\fi : #2}}
\definecolor{darkgreen}{rgb}{0,0.5,0}
\newcommand*\tocheck[1]{\marginpar{ \raggedright \color{darkgreen} #1}}
\newcommand*\oldtodo[2][]{\tocheck{%
\ifx\\#1\\Old TODO\else Old TODO(#1)\fi : #2}}
\newcommand*\defn[1]{\textit{#1}} 
\newcommand*\an{\textnormal{\texttt{an}}}
\newcommand*\aff{\textnormal{\texttt{aff}}}
\newcommand*\tmax{\textnormal{\texttt{max}}}
\newcommand*\up{\textnormal{\texttt{up}}}
\newcommand*\down{\textnormal{\texttt{down}}}

\newcounter{locallabel}[theorem]
\newcommand*\llabel[1]{\label{locallabel:sec-\thesection:thm-\thetheorem:loc-\thelocallabel:#1}}

\newcommand*\leqref[1]{\eqref{locallabel:sec-\thesection:thm-\thetheorem:loc-\thelocallabel:#1}}

\makeatletter
\newcommand*\vbig[1]{\bBigg@{#1}}
\makeatother

\DeclareMathOperator{\ent}{ent}
\DeclareMathOperator{\Ent}{Ent}

\DeclareMathOperator{\argmin}{argmin}
\DeclareMathOperator{\argmax}{argmax}

\DeclareMathOperator{\Adm}{Adm}
\DeclareMathOperator{\Lip}{Lip}
\DeclareMathOperator{\Sym}{Sym}
\DeclareMathOperator{\diam}{diam}

\DeclarePairedDelimiter\ceil{\lceil}{\rceil}
\DeclarePairedDelimiter\floor{\lfloor}{\rfloor}
\DeclarePairedDelimiter\abs{\lvert}{\rvert}
\DeclarePairedDelimiter\norm{\lVert}{\rVert}

\mathtoolsset{showonlyrefs}

\begin{document}

\title[Homogenization of the variational principle]{Homogenization of the variational principle for discrete random maps.}

\author{Andrew Krieger}
\address{Department of Mathematics, University of California, Los Angeles}
\email{akrieger@math.ucla.edu}

\author{Georg Menz}
\address{Department of Mathematics, University of California, Los Angeles}
\email{gmenz@math.ucla.edu}

\author{Martin Tassy}
\address{Dartmouth College, Hanover}
\email{mtassy@math.dartmouth.edu}


\date{\today}

\begin{abstract}
We consider homogenization of random surfaces
and study the variational principle for graph homomorphisms
from subsets of~$\Z^m$ into~$\Z$,
where the underlying uniform measure is perturbed by a random field.
Motivated by the theories of random walks in random potentials,
we assume that random field is stationary, ergodic, and bounded in~$L^1$.
We show that the variational principle holds in probability
and that the entropy functional homogenizes,
i.e.\ is independent of the values taken by the random field.
The main ingredients in the argument are
the existence of the quenched surface tension,
the equivalence of the quenched and the annealed surface tension,
and robustness of the surface tension under change in boundary data.
These ingredients are deduced by a combination of
a superadditive ergodic theorem and combinatorial results,
especially the Kirszbraun theorem.
\end{abstract}

\maketitle


\tableofcontents

\section{Introduction}


The broader scope of this article is the study of limit shapes
as a limiting behavior of discrete systems.
Limit shapes are a well-known and studied phenomenon
in statistical physics and combinatorics (e.g.~\cite{Georgii}).
Among others, models that exhibit limit shapes include
domino tilings and dimer models (e.g.~\cite{Kas63,CEP96,CKP01}),
polymer models (e.g.~\cite{BiPr18,BeYa19}),
lozenge tilings (e.g.~\cite{Des98,LRS01,Wil04}),
Ginzburg-Landau models (e.g.~\cite{DeGiIo00,FuOs04}),
Gibbs models (e.g.~\cite{She05}),
the Ising model (e.g.~\cite{DKS92,Cer06}),
asymmetric exclusion processes (e.g.~\cite{FS06}),
sandpile models (e.g.\cite{LP08}),
the six vertex model (e.g.~\cite{BCG2016,CoSp16,ReSr16}),
and the Young tableaux (e.g.~\cite{LS77,VK77,PR07}).
\\

\begin{figure}[b]
    \centering
 \includegraphics[width=0.65\textwidth]{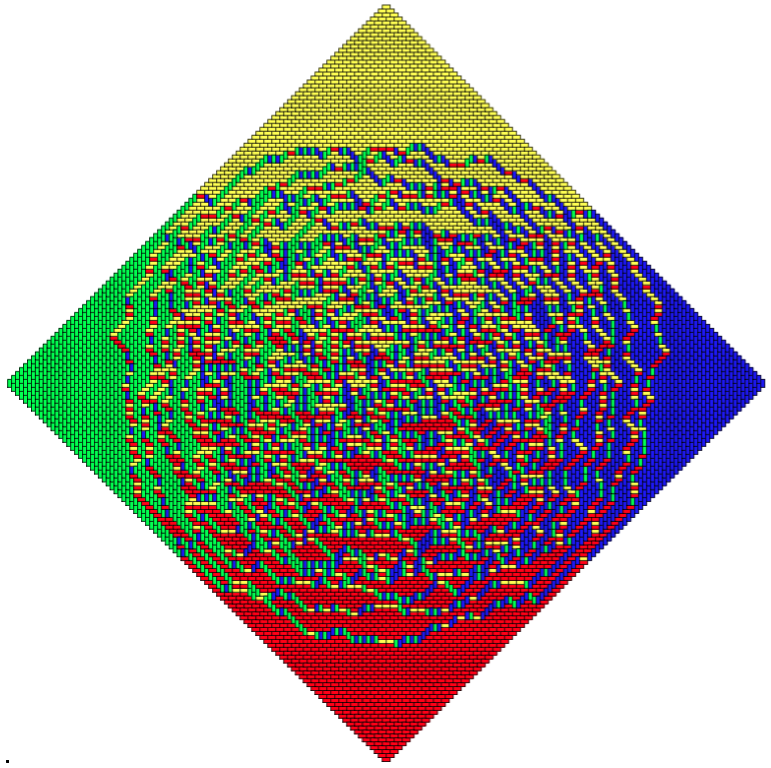}
\caption{ An Aztec diamond for domino tilings. The combinatorics of the model is similar to Lipschitz functions from $\Z^2$ to $\Z$ (see~\cite{CKP01}).}\label{f_aztec_domino_tilings}
\end{figure}

Limit shapes appear in stiff models whenever fixed boundary conditions
force a certain response of the system.
The numerous examples in the literature and many simulations
show that the existence of limit shapes is a universal phenomenon.
Among many possible references, let us just mention
\cite{Ken00,BoGoGu15,MeTa16,MoPaTa18,KeSt18}.
Several new approaches were developed recently to make methods more robust;
see for example~\cite{ChJoYo15,BuKn16,CoSp16,Agg19}.
A part of the effort to develop robust and universal methods
to deduce variational principles are the articles~\cite{MeTa16}.
In~\cite{MeTa16} variational principles were studied
in target spaces where the usual cluster swapping methods do not work.
\\

In this article we explore a new direction
and show the robustness of the variational principle in a random potential.
The basic objects for our model are graph homomorphisms
from finite subsets of the $m$-dimensional lattice $\Z^m$ into $\Z$,
also called height functions.
In two dimensions and without random potential,
$\Z$-homomorphisms are equivalent to a special case of the six-vertex model,
where all vertex weights are identical, i.e.\@ the square-ice model.
The limiting behavior of the $\Z$-homomorphism model without random potential
is well-studied; see for example \cite{BeHaMo00,PeSaYe12,Pel17},
as well as the companion article~\cite{KrMeTa19} in which the current authors
describe a robust method of proof for the variational principle and large
deviations principle for $\Z$-homomorphisms.
\\

The limit results alluded to above are about a $\Z$-homomorphism
that is chosen uniformly at random from some (finite) set
of admissible homomorphisms
(e.g.\ the set of homomorphisms with prescribed boundary values).
The random potential that we introduce in this article
perturbs the uniform measure $\mu$ in the above results,
replacing it by a weighted measure $\mu_\omega$,
where $\omega$ denotes the random potential.
Our reason for introducing a random potential is to test the robustness
of the methods used to prove variational principles and similar results.
Several nice properties do not carry over from the unperturbed model:
exact computations like those in~\cite{CKP01} are prohibitively difficult,
\textit{a priori} proofs of concentration
(e.g.\ the martingale method of~\cite{CEP96}) do not seem to apply,
and there are no obvious global symmetries.
To overcome those obstacles we make use of ergodicity and homogenization.
\\

The random potential is inspired by
homogenization of random walks in random environment
(see e.g.\ the survey~\cite{Bis11}).
Indeed, the bridge model of~\cite{GP11},
i.e.\@ transient random walks in random environment
conditioned to start and end at prescribed boundary values,
is a special case of the $\Z$-homomorphism in random potential model
with dimension $m=1$.
The bridge model exhibits asymptotically different maximal order statistics
than does the bridge model originating from the simple random walk.
The model considered in this article is a natural extension of bridges
to ``random sheets.''
It would be interesting to extend the results of~\cite{GP11}
to higher dimensions and compare against the Gaussian free field.
\\

\begin{figure}
	\subcaptionbox{
		A height function sampled without random potential.
		\label{f_sim_zero}
	}[0.3\textwidth]{%
		\includegraphics[width=0.3\textwidth]{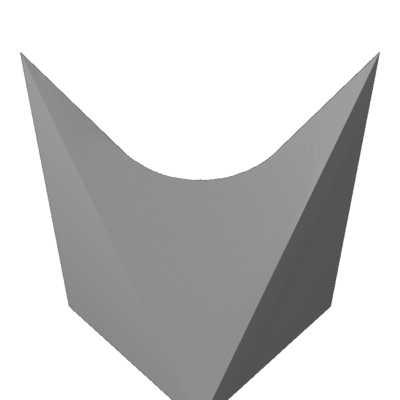}%
		\hfil%
	} \hfil
	\subcaptionbox{
		Two height functions sampled from $\mu_\omega$,
		under two different samples of $\omega$
		with a.s.\ bounded distribution.
		Specifically, the random variables
		$\{\omega_e \,|\, e \in E(\Z)\}$ are i.i.d.\
		with $\P(\omega_e = 1) = \P(\omega_e = -1) = \tfrac{1}{2}$.
		\label{f_sim_bernoulli}
	}[0.6\textwidth]{%
		\includegraphics[width=0.3\textwidth]{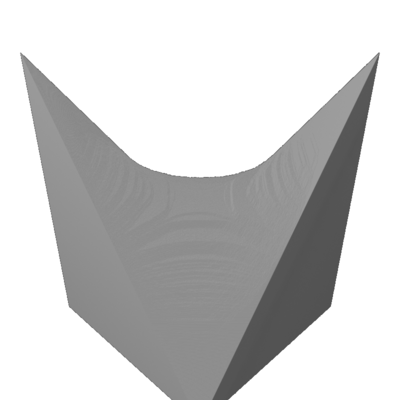}%
		\includegraphics[width=0.3\textwidth]{bernoulli-29-40B-3d-lowres}%
	} \vskip1em
	\subcaptionbox{
		Three height functions sampled from $\mu_\omega$,
		under three different samples of $\omega$
		with a.s.\ unbounded distribution.
		Specifically, the random variables
		$\{\omega_e \,|\, e \in E(\Z)\}$ are i.i.d.\
		standard normal variables.
		\label{f_sim_gauss}
	}[0.9\textwidth]{%
		\includegraphics[width=0.3\textwidth]{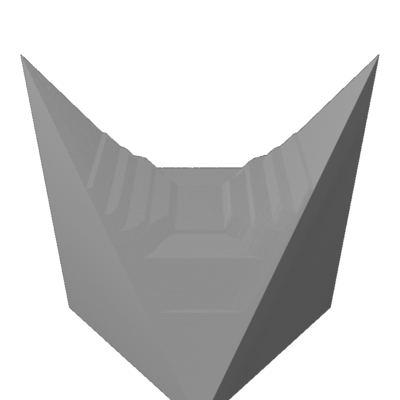}%
		\includegraphics[width=0.3\textwidth]{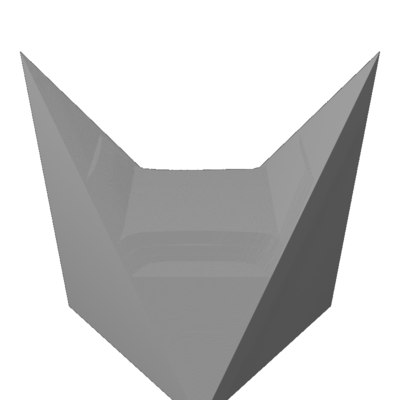}%
		\includegraphics[width=0.3\textwidth]{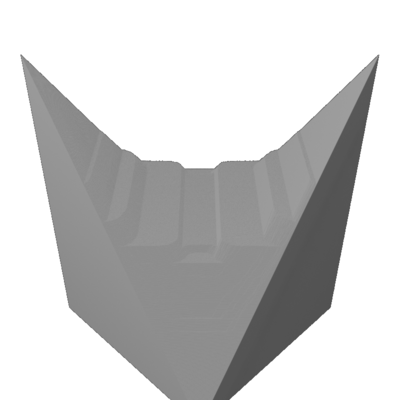}%
	}
	\caption{
		Height functions $h$, sampled from random measures $\mu_\omega$,
		which in turn are derived from randomly sampled fields $\omega$.
		The height functions are rendered as 3D solids
		with the random surface $\{(x,y, h(x,y)) \,|\, (x,y) \in S_n\}$
		as their ``top'' face.
		Boundary values of $h$ are fixed,
		and the behavior of $h$ on the interior of the domain $S_n$
		depends on $\omega$.
		If the law of $\omega$ is bounded
		(as is the case in Figure~\ref{f_sim_bernoulli}),
		then the results in this article
		imply that with high probability,
		the macroscopic behavior of $h$ does not depend on the specific
		realization of $\omega$, but only on its distribution.
	}
	\label{f_sim}
\end{figure}

Before summarizing the mathematical results of the article,
let us motivate the model further by discussing empirical results
from computer simulation; cf.\ Figure~\ref{f_sim}.
We generated random environments~$(\omega_e)_{e \in E(S_n)}$
according to various distributions (e.g.\ i.i.d.~Gaussian)
and for various box sizes (up to a $1000 \! \times \!1000$-vertex box).
We chose boundary data~$h_{\partial S_n}$,
then we sampled a height function~$h \in M(S_n, h_{\partial S_n})$
according to the random measure~$\mu_\omega$
using the Markov chain Monte Carlo method.
\\

We call attention to a few details from the simulations.
The two height functions in Figure~\ref{f_sim_bernoulli}
are drawn from two different measures $\mu_\omega$,
where the random potentials $\omega$ are sampled such that
$\{\omega_e \,|\, e \in E(\Z)\}$ are
i.i.d.\ with $\P(\omega_e = 1) = \P(\omega_e = -1) = \tfrac{1}{2}$.
Although the exact value of $\omega$ varies in the two samples,
the randomly chosen height functions in the pictures
appear to be macroscopically identical.
This is a good indicator that this model homogenizes.
By this we mean that the macroscopic features measure $\mu_\omega$ do not depend
(in the limit, except with negligible probability)
on the exact choice of $\omega$.
Rather those macroscopic features of $\mu_\omega$
only depend on the distribution of $\omega$ and the boundary data.
Indeed, the main results of this article apply to the random potential
from Figure~\ref{f_sim_bernoulli},
so we know that this model homogenizes.
\\

The three height functions in Figure~\ref{f_sim_gauss}
are sampled from three different measures $\mu_\omega$,
where $\omega_e \sim \NN(0,1)$ are i.i.d.\ standard normal variables.
Notice that these three height functions differ macroscopically,
depending on the realization of $\omega$.
This does not contradict the results of this article
because the random potential is unbounded.
We expect (but have not proven) that this model fails to homogenize
when the random potential is unbounded,
with energetic effects from $\omega$ overwhelming
the entropic effects from the underlying combinatorial $\Z$-homomorphism model.
If the random measure $\mu_\omega$ does not homogenize,
then the limit shape under $\mu_\omega$ may depend
on the actual values of $\omega$,
and thus be may be a non-trivial random variable
(with respect to the randomness that determines $\omega$).
\\

In order to understand the behavior
underlying the simulations in Figure~\ref{f_sim},
we prove two main results:
a profile theorem (see Theorem~\ref{th_profile}) and a variational principle (see Theorem~\ref{th_varprin}).
A third related result, namely a large deviations principle (see Theorem~\ref{th_ldp}),
is not difficult to prove by the same methods,
but we omit it for brevity.
These results hold with high probability for fixed environments~$\omega$.
They establish that, for typical samples of $\omega$,
there holds a conclusion similar to the profile theorem or variational principle
for the non-random model studied in our companion article \cite{KrMeTa19}.
Indeed the purpose of the companion article was
to distill, simplify, and explain the steps involved in proving these results.
Understanding the methods in the companion article
will help to understand the general outline of the proofs in this article.
From a high-level perspective,
the main results are similar to the simpler case studied in~\cite{KrMeTa19}.
\\

Let us now briefly discuss the main results of this article. We start with the profile theorem. It asymptotically characterizes the cardinality of the set of height functions~$h_{R_n}$ that are uniformly close to a particular macroscopic state~$h_R$ (also called asymptotic height profile later on).
Without random potential, the profile theorem states
(cf.\ \cite[Theorem 15]{KrMeTa19}) that
\begin{equation} \label{e_intro_profthm_nopot} \begin{aligned}
	\hskip3em&\hskip-3em
	\Ent_{R_n} \Bigl( \Bigl\{
		\text{height functions $h_{R_n}: R_n \to \Z$} \\
	& \hskip6em \text{with $\norm{ h_{R_n} - \tilde h_R }_\infty
		< \varepsilon$}
	\Bigr\} \Bigr) \\
	&\approx \Ent_R(h_R) \,,
\end{aligned} \end{equation}
where (for details see Section~\ref{s_main_results}):
\begin{align}
	R_n &\subset \Z^m \text{ is such that $\tfrac{1}{n} R_n$
		converges to $R$} \,, \\
	\tilde h_R(z) &= n h_R(\tfrac{z}{n})
	\text{ is a rescaled version of $h_R$} \,, \\
	\label{e_intro_profthm_nopot_micro}
	\Ent_{R_n}(M) &= - \frac{1}{\abs{R_n}} \log \abs{M} \,, \\
	\label{e_intro_profthm_nopot_macro}
	\Ent_R(h_R) &= \int_R \ent(\nabla h_R(x)) \,,\, \text{and} \\
	\ent: [-1,1]^m \to \R &\text{ is determined by the combinatorics}\\
	&\qquad \text{of the $\Z$-homomorphism model.}
\end{align}
In the setting of homogenization we substitute the uniform measure on the set of microscopic height functions with a random measure~$\mu_{\omega}$ that is characterized by the random potential~$\omega$. The quantity $\Ent_{R_n}(M)$ from~\eqref{e_intro_profthm_nopot_micro}
is dependent on $\omega$, and is therefore a random variable.
Specifically, the cardinality $\abs{M}$ is replaced by a partition function of~$\mu_\omega$ on the set of height functions.
The quantity $\Ent_R(h_R)$
from the right-hand side of~\eqref{e_intro_profthm_nopot}
is replaced by $\Ent_{R,\an}(h_R)$,
the annealed macroscopic entropy (see Definition~\ref{d_ann_macro_ent}).
Likewise the local surface tension $\ent(\cdot)$
in the definition~\eqref{e_intro_profthm_nopot_macro}
is replaced by the annealed local surface tension $\ent_\an(\cdot)$.
In both cases, ``annealed'' means that the influence of the random field
$\omega$ is averaged out i.e.~$\ent_{\an}(s) := \mathbb{E}[\ent(s,\omega)]$.
Therefore $\ent_\an(\cdot)$ and $\Ent_{R,\an}(\cdot)$ are non-random.
Turning back to the conclusion~\eqref{e_intro_profthm_nopot}
of the profile theorem,
the left-hand side is a non-trivial random variable,
but the limiting quantity on the right-hand side is not random.
\medskip

Let's turn to the second main result, namely the variational principle. Recall that the profile theorem measures the set of height functions
that stay close to a target asymptotic height profile over the entire domain. The variational principle instead measures the whole set of height functions with certain boundary values.
Without random potential (cf.\ \cite[Theorem 16]{KrMeTa19}), the result is
\begin{equation} \label{e_intro_varprin_nopot} \begin{aligned}
	\hskip3em&\hskip-3em
	\Ent_{R_n} \bigl( \bigl\{ h_{R_n}: R_n \to \Z \,\big|\,
		\text{$h_{R_n}|_{\partial R_n}$ is close to $h_{\partial R}$}
	\bigr\} \bigr) \\
	&\approx \inf_{h_R} \Ent_R(h_R) \,,
\end{aligned} \end{equation}
where ``close'' means close in the supremum norm after rescaling,
and where the infimum runs over all asymptotic height profiles
consistent with the given boundary data $h_{\partial R}: \partial R \to \R$.
In the setting of homogenization, i.e.~adding a random potential to the uniform measure,
the necessary modifications to this approximate identity
are analogous to those for the profile theorem above:
$\Ent_{R_n}(\cdot)$ becomes a random variable dependent on $\omega$
just as above,
and $\Ent_R(h_R)$ is again replaced by the non-random quantity
$\Ent_{R,\an}(h_R)$.
Hence, it follows from our main result that the variational principle homogenizes.\medskip

Let's now discuss the large deviations principle.
Let $h_{\partial R}$ be an asymptotic boundary height function
and let $A$ be a Borel set in the space of asymptotic height functions
with boundary values given by $h_{\partial R}$, equipped with the supremum norm.
Without random potential the large deviations principle states
(cf.\ \cite[Theorem 17]{KrMeTa19}):
\begin{equation} \label{e_intro_ldp_nopot} \begin{aligned}
	\hskip3em&\hskip-3em
	-\frac{1}{\abs{R_n}} \log \mu_n \bigl( \bigl\{
		h_{R_n}: R_n \to \Z \,\big|\,
		\text{after rescaling, $h_{R_n} \in A$} \bigl\} \bigl) \\
	&\approx \inf_{h_R \in A} \Ent_R(h_R) - E_0\,,
\end{aligned} \end{equation}
where $\mu_n$ is the uniform measure on the set of (microscopic) height functions
with appropriate boundary values,
\and where $E_0 := \inf_{h_R} \Ent_R(h_R)$ is the infimum of the entropy
over all asymptotic height functions with boundary values $h_{\partial R}$.

In the setting of homogenization, the large deviation principle needs
to be adapted in an analogous way as
for the profile theorem and the variational principle:
The uniform measure~$\mu_n$ is replaced by the random measure~$\mu_{n,\omega}$,
and $\Ent_R(\cdot)$ is replaced by the annealed macroscopic entropy.
Then the large deviation principle holds again
with respect to sample $\omega$ of the random field with high probability.
Because the rate functional homogenizes,
i.e.~it is independent of the realization~$\omega$ of the random field,
the large deviation principle would also homogenize.
\\

It is natural to ask whether the infima in~\eqref{e_intro_varprin_nopot}
and~\eqref{e_intro_ldp_nopot} admit a minimizer,
and if so whether the minimizer is unique.
Existence follows from convexity of the integrand function $\ent_\an(\cdot)$
(Lemma~\ref{p_convexity} establishes convexity of $\ent_\an(\cdot)$).
Uniqueness of $h_{R,\min}$ follows from strict convexity of $\ent_\an(\cdot)$;
see for example \cite{CKP01}.
Strict convexity remains an open question for this model.
\\

The proofs of the main results are based on two main ingredients:
existence and characterization of the quenched local surface tension
and robustness of the entropy.
\\

The first main ingredient is
the existence of the quenched local surface tension $\ent_\an(s,\omega)$.
Without random potential existence follows from superadditivity
by application of Fekete's lemma.
With random potential we turn to a superadditive ergodic theorem instead.
superadditivity and translation invariance are enough to establish existence
of the quenched local surface tension.
Ergodicity is used to characterize the quenched local surface tension.
At slopes $s \ne 0$, translating a domain $R_n$ by $z \in \Z^m$
implies shifting the boundary heights by $s \cdot z$,
and the random potential is ergodic with respect to this kind of height shift.
Since the quenched local surface tension is translation invariant,
it follows that it is almost surely equal to its expectation,
the annealed local surface tension.
The same conclusion holds in the case $s=0$, which we show using an argument
with credit to Marek Biskup.
\\

The second main ingredient in proving the main results of this article
is robustness.
In Section~\ref{s_local_surface_tension} we prove several results,
which serve to control the change in the microscopic entropy
$\Ent_{R_n}(A,\omega)$ as the set of height functions $A$ changes.
For example, when
\[
	A = \bigl\{
		h_{R_n} \,\big|\, h_{R_n}|_{\partial R_n} = h_{\partial R_n}
	\bigr\}
\]
is defined by boundary data $h_{\partial R_n}$,
we consider the effect of changing the boundary data.

The main idea used to control the change of microscopic entropy is to use the Kirszbraun theorem (see Theorem~\ref{p_kirszbraun}). It allows to extend height functions on a domain
to height functions on the larger domain.
This provides an injection between the two sets of height functions,
and it remains to control the energetic effect
contributed by the newly added edges in the larger domain.
When the Kirszbraun theorem is not useful,
we fall back to combinatorial results.
\\

The proof of the robustness results illustrate a primary source of difficulty:
passing from combinatorial estimates on the number of height functions
to control over energetic effects arising from the random potential.
In the example discussed above, every height function in the smaller
(in the sense of cardinality) set admits an extension in the larger set.
It is not difficult to compare the total energy of an extension
to that of the original height function,
using the assumption that the random potential is bounded.
\\

After applying the two ingredients listed above,
it remains to apply approximations of Lipschitz functions
and compactness of the space of asymptotic height functions
(with fixed boundary values).
For these last steps of the argument
we follow Sections~5 through~8 of~\cite{KrMeTa19},
with some modifications needed to account for the random potential.
Because the proof is largely the same as in our companion article we do not go into great detail for these steps.
\\

The rest of this article is organized as follows.
\begin{itemize}
	\item In Section~\ref{s_main_results} we define the precise setting and state the main results.
	\item In Section~\ref{s_local_surface_tension} we state and prove
	key results about the local surface tension.
	\item In Section~\ref{s_profile} we prove the first main result,
	namely the profile theorem.
	\item In Section~\ref{s_varprin} we prove the second main result,
	namely the variational principle.
	\item In Section~\ref{s_open} we state a few open problems
	and directions for further research.
\end{itemize}

\section*{Notation and conventions}
For the convenience of the reader,
we summarize the basic notation that we use throughout this article.

\begin{itemize}

\item $\abs{A}$ denotes either the cardinality or the Lebesgue measure of the
      set~$A$, depending on context.
\item $S_n : = \left\{ -n, -(n-1), \ldots, n-1, n \right\}^m \subset \Z^m$
      denotes a hypercube in the lattice, centered at the origin.
\item For $z, z' \in \Z^m$, $z \sim z'$ means that $z$ and $z'$ are
      nearest-neighbors (i.e.\ the $\ell^1$ distance
      $\sum_{i=1}^m \abs{z_i - z'_i}_1$ is exactly 1).
\item For $S \subset \Z^m$, $\partial S
          := \{ z \in S
          \,|\, \exists \tilde z \in \Z^m \setminus S, \, \tilde z \sim z \}
      $ is the (interior) boundary of $S$.
\item $e_{zz'}$ is the unoriented edge between neighbors $z \sim z'$ in $\Z^m$.
\item For $h: \Z^m \to \Z$ and $e = e_{zz'} \in E(\Z^m)$,
      we abuse notation and write $h(e)$ for the edge $e_{h(z), h(z')} \in E(\Z)$.
\item $\tau_w$ denotes the shift by~$w \in \Z^m$ on edges of the graph~$\Z^m$.
      That is, $\tau_w e_{zz'} = e_{z+w, z'+w}$.
\item $s \in \R^m$ denotes a vector satisfying~$\abs{s}_\infty \le 1$.
\item $\theta(\varepsilon) \ge 0$ denotes a smooth function
      with~$\lim_{\varepsilon \to 0} \theta(\varepsilon) = 0$.
      More precisely, we use the same notation/convention for~$\theta$
      as explained in~\cite[Section 2.4]{KrMeTa19}
\end{itemize}

\section{Setting and main results} \label{s_main_results}

In this section we describe the model under study,
introduce related notation,
and state the main results of this article.
The setting, notation, and main results are similar to
those of the companion article~\cite{KrMeTa19}.

\subsection{Basic definitions}

Throughout the sequel, we fix a dimension $m \in \N$,
a macroscopic domain $R \subset \R^m$,
and a sequence of microscopic domains $R_n \subset \Z^m$
satisfying these assumptions:

\begin{assumption}[Assumptions on domain $R$ and $R_n$] \label{a_domain}
We assume that $R \subset \R^m$ is compact and connected,
that $R$ is the closure of its interior,
and that the boundary of $R$ has zero Lebesgue measure.
We assume that $R_n \subset \Z^m$
is contained in $R$ after rescaling,
i.e.\ that $\frac{1}{n} R_n \subset R$,
although this is just a simplifying assumption.
Moreover, we assume that $\frac{1}{n} R_n \to R$ in the Hausdorff metric,
i.e.\ the metric on $\{A \subset \R^m\}$ defined by
\begin{equation} \label{e_hausdorff}
	d_H(A, B)
	:= \biggl( \: \sup_{x \in A} \inf_{y \in B} \abs{x-y}_1 \biggr)
	\vee \biggl( \: \sup_{y \in B} \inf_{x \in A} \abs{x-y}_1 \biggr) \,.
\end{equation}
\end{assumption}

Now, we define precisely the height functions in our model.

\begin{definition}[Height function] \label{d_ht_func}
A \defn{height function} on $R_n$
is a graph homomorphism $h_{R_n}: R_n \to \Z$.
In other words, if $z, w \in R_n$ and $z \sim w$,
then $\abs{ h_{R_n}(z) - h_{R_n}(w) } = 1$,
and for any $z = (z_1, \dotsc, z_m) \in R_n$,
\begin{equation} \label{e_ht_func_parity}
	h_{R_n}(z) \equiv z \pmod 2 \,,
	\quad \text{i.e. } h_{R_n}(z) \equiv \sum_{i=1}^m z_i \pmod 2 \,.
\end{equation}
\end{definition}
The condition~\eqref{e_ht_func_parity} states that a height function
preserves the parity of the lattice~$\mathbb{Z}^m$.
Indeed, every graph homomorphism either preserves parity at all points
or inverts parity at all points,
since the source space $\Z^m$ and the target space $\Z$ are both bipartite.
Our main results are also valid without the parity-preserving condition,
but for the same reasons as outlined in~\cite[Section 2.1]{KrMeTa19}
we include it for simplicity.

We introduce the following symbols to refer to sets of height functions:

\begin{definition}[Sets of height functions] \label{d_ht_func_sets}
Let $R_n$ be a microscopic domain as above,
let $h_{R_n}: R_n \to \Z$ be a boundary height function,
and let $\delta > 0$.
We define:
\begin{align}
	M(R_n)
	&:= \bigl\{ h_{R_n}: R_n \to \Z
	\, \big| \, \text{$h_{R_n}$ is a height function} \bigr\} \,, \\
	M(R_n, h_{\partial R_n})
	&:= \bigl\{ h_{R_n} \in M(R_n)
	\, \big| \, h_{R_n}|_{\partial R_n} = h_{\partial R_n} \bigr\} \,, \\
	M(R_n, h_{\partial R_n}, \delta)
	&:= \bigl\{ h_{R_n} \in M(R_n)
	\, \big| \, \sup_{z \in \partial R_n}
		\abs{ h_{R_n}(z) - h_{\partial R_n}(z) } < \delta n \bigr\}
	\,,\, \text{and} \\
	B(R_n, h_R, \delta)
	&:= \bigl\{ h_{R_n} \in M(R_n)
	\, \big| \, \sup_{z \in R_n}
		\abs{ h_R(\tfrac{1}{n} z) - \tfrac{1}{n} h_{R_n}(z) } < \delta
	\bigr\} \,.
\end{align}

In the last definition, the expression ``$h_R(\tfrac{1}{n}z)$'' makes sense
because of the assumption that $\tfrac{1}{n}R_n \subset R$
in Assumption~\ref{a_domain}.
\end{definition}

The limiting object for convergent sequences of height functions is:

\begin{definition}[Asymptotic height function] \label{d_asymp_ht_func}
We call a function $h_R: R \to \R$ an \defn{asymptotic height function}
if $h_R$ is Lipschitz with Lipschitz constant at most $1$,
with respect to the $\ell^1$-norm on $\R^m$;
that is, if
\begin{equation} \label{e_asymp_ht_func}
	\operatorname{Lip}(h_R)
	:= \sup_{x \ne y \in R} \frac{\abs{h_R(x) - h_R(y)}}{\abs{x-y}_1}
	\le 1 \,.
\end{equation}

Likewise, if $h_{\partial R}: \partial R \to \R$ is $1$-Lipschitz
(with respect to the $\ell^1$-norm),
we call $h_{\partial R}$ an \defn{asymptotic boundary height function}.
\end{definition}

The limit of height functions is defined as follows.

\begin{definition}[Convergence of height functions] \label{d_limit_ht_func}
Given a sequence of height functions $h_{R_n}: R_n \to \Z$
and an asymptotic height function $h_R: R \to \R$,
we say that $h_{R_n}$ converges in the scaling limit to $h_R$ if
\begin{equation} \label{e_limit_ht_func}
	\lim_{n \to \infty} \: \sup_{z \in R_n} \:
		\sup_{\substack{
			x \in R \\
			\abs{ x - \frac{1}{n} z } _1 \le d_n
		}} \:
			\abs[\Big]{ \, \frac{1}{n} h_{\partial R_n}(z)
			- h_{\partial R}(x) \, } = 0 \,,
\end{equation}
where $d_n := d_H(\tfrac{1}{n} R_n, R)$.
\end{definition}

Finally, we define the following sets of asymptotic height functions:

\begin{definition}[Sets of asymptotic height functions]
\label{d_asymp_ht_func_sets}
Let $R \subset \R^m$ be a domain satisfying Assumption~\ref{a_domain},
let $h_{\partial R}: \partial R \to \R$
be an asymptotic boundary height function,
and let $\delta > 0$.
We define:
\begin{align}
	M(R)
	&:= \bigl\{ h_R: R \to \R \, \big| \,
		\text{$h_R$ is an asymptotic height function} \bigr\} \,, \\
	M(R, h_{\partial R})
	&:= \bigl\{ h_R: R \to \R \, \big| \,
		h_R|_{\partial R} = h_{\partial R} \bigr\} \,, \\
	M(R, h_{\partial R}, \delta)
	&:= \bigl\{ h_R: R \to \R \, \big| \,
		\forall x \in \partial R \,,\,
			\abs{ h_R(x) - h_{\partial R}(x) } \le \delta
	\bigr\} \,,\, \text{and} \\
	B(R, \tilde h_R, \delta)
	&:= \bigl\{ h_R: R \to \R \, \big| \,
		\forall x \in R \,,\, \abs{ h_R(x) - \tilde h_R(x) } < \delta
	\bigr\} \,.
\end{align}
\end{definition}

\subsection{Defining the entropy}

In order to define the local surface tension, both quenched and annealed,
we fix a family of canonical height functions with fixed slope.
These are the linear and affine height functions,
so called because they approximate
linear and affine functions of real variables.

\begin{definition}[Affine and linear height functions] \label{d_aff_ht_funcs}
For $s \in [-1,1]^m$, $b \in \R$, and $n \in \N$,
we define the affine height function $h^{s \cdot x + b} \in M(\Z^m)$ as
\begin{equation}
	h^{s \cdot x + b}(z) := [ s \cdot z + b ]_{z \bmod 2}
	\quad \text{for all $z \in \Z^m$} \,,
\end{equation}
where for $t \in \R$ and $z \in \Z^m$, $[t]_{z \bmod 2}$
is the integer with the same parity as $z$
that is closest to $z$.
(In the ambiguous case,
namely when $t$ is an integer having opposite parity as $z$,
we choose arbitrarily but consistently to ``round up''
and set $[t]_{z \bmod 2} = (z+1)$.)
For $s \in [-1,1]^m$, the linear height function $h^s \in M(\Z^m)$
is given by $h^s = h^{s \cdot x + 0}$, i.e.
\begin{equation}
	h^s(z) := [ s \cdot z ]_{z \bmod 2}
	\quad \text{for all $z \in \Z^m$} \,,
\end{equation}
\end{definition}

\begin{remark}
The symbol ``$x$'' in the superscript ``$s \cdot x + b$'' is a formal variable,
used so that the superscript resembles a meaningful expression
instead of, say, the less intuitive pair $(s,b)$.
It is not difficult to verify that
the functions defined above are graph homomorphisms.
We refer the reader to \cite[Lemma 7]{KrMeTa19} for the details.
\end{remark}

Until now, the setup has been the same as
in the companion article~\cite{KrMeTa19}.
Let us now turn to homogenization and to the new contributions of this article.
The main change in the model is that
instead of the uniform measure on $M(R_n, h_{\partial R_n})$
we consider a noisy perturbation~$\mu_{\omega}$ of the uniform measure,
where $\omega = (\omega_e)_{e \in E(\Z)}$ denotes a random field,
as described in Assumption~\ref{a_random_field}.

\begin{assumption}[Random field~$\omega$]\label{a_random_field}
 We consider a real-valued random potential
\begin{equation}
	\omega = (\omega_e)_{e \in E(\Z)} \in \R^{E(\Z)}
\end{equation}
defined on the set of edges~$E(\Z)$ of~$\Z$.
We assume that~$\omega$ satisfies the following assumptions:
\begin{itemize}
	\item The random field~$\omega$ is almost surely finite,
	and moreover the random variable $C_\omega$ defined by
	\begin{align}\label{e_boundedness_random_field}
		 C_\omega := 1 \vee \sup_{e \in E(\Z)} \abs{ \omega_{e} }
	\end{align}
	is in $L^1$,
	i.e.\ $\mathbb{E}[C_\omega] < \infty$.

\item The random field~$\omega$ is shift invariant.
	This means that for any finite number of edges
	$e_1, \ldots e_k \in E(\Z)$, any integer~$z \in \Z$,
	and any bounded and measurable function
	$\xi: \mathbb{R}^k \to \mathbb{R}$,
	\begin{equation}\label{e_shift_invariance}
		\E \bigl[ \xi(\omega_{e_1} , \ldots, \omega_{e_k} ) \bigr]
		= \E \bigl[ \xi(\omega_{\tau_z (e_1)} , \ldots,
			\omega_{\tau_z (e_k)} ) \bigr] \,,
	\end{equation}
	where~$\tau_z: E(\Z) \to E(\Z)$ is the shift by~$z$
	(as per the Notation and Conventions above).
\item Moreover, the random field~$\omega$ is ergodic
	with respect to the set of shifts
	$\{\tau_z \,|\, z \in \Z, \, z \equiv 0 \pmod 2 \}$.
	This means that if $E \subset \Omega$ is a shift invariant event,
	i.e.\ if $E = \tau_2^{-1}(E)$,
	then $\P(E) \in \{0, 1\}$.
\item We assume w.l.o.g. (as a matter of normalization) that
	\begin{equation} \label{e_wlog_mean_0}
		\E [ \omega_{e_{0,1}} ] = 0 \,,
	\end{equation}
	where $e_{0,1}$ is the edge from $0$ to $1$ in $\Z$.
\end{itemize}
\end{assumption}

\begin{example}
The simplest non-trivial example of a random field~$\omega$
that satisfies Assumption~\ref{a_random_field}
is the i.i.d.\ field.
Let~$X$ denote a bounded (real) random variable with mean~$0$,
and let $(\omega_e)_{e \in E(\Z)}$ denote a family of i.i.d.\ copies of $X$.
\end{example}

\begin{remark}
The assumptions of shift invariance and ergodicity are standard
in homogenization literature;
see for example the ``usual conditions'' for the random conductance model
from~\cite[Definition~3.1]{Bis11}.
However we point out one difference:
the random field $\omega$ is ergodic with respect to the even
shifts $\{ \tau_z \,|\, z \equiv 0 \pmod 2 \}$.
This is a stronger condition than being ergodic
with respect to the full set of shifts $\{ \tau_z \,|\, z \in \Z \}$.
This requirement is due to the earlier assumption
made in Definition~\ref{d_ht_func}
that height functions preserve parity.
As such, we cannot simply shift a height function
up or down by $1$ in the height space;
if $h_{S_n}(z) = k \in \Z$,
then there is no (parity-preserving) height function ``$\tau_1 h_{S_n}$''
such that $\tau_1 h_{S_n}(z) = k+1$.
More concretely, the family of measure-preserving translations
used in the proof of Lemma~\ref{l_que_local_surface_tension} below
includes all of the shifts $\{ \tau_z \,|\, z \equiv 0 \pmod 2 \}$
and none of the shifts $\{ \tau_z \,|\,z \equiv 1 \pmod 2 \}$,
hence the stronger ergodicity assumption is technically required.
\end{remark}

In this article we study the random surfaces in the random potential
defined by $\omega$.
In homogenization one considers two different situations:
In the quenched case, one considers the measure~$\mu_\omega$ for fixed~$\omega$.
In the annealed case, one takes the expectation with respect to~$\omega$.
Our goal is to show that the variational principle holds
with high probability.
With that context in mind,
we define the quenched Hamiltonian~$H_{R_n}(\cdot) = H_{R_n}(\cdot, \omega)$
and the quenched measure~$\mu_{\omega}$ as follows:

\begin{definition}[The quenched Hamiltonian]
For finite subsets $R_n \subset \Z^m$,
We define the Hamiltonian $H_{R_n}$ as follows:
for a fixed boundary height function $h_{\partial R_n}: \partial R_n \to \Z$,
and for any height function $h_{R_n} \in M(R_n, h_{\partial R_n})$
and any realization $\omega$ of the random field,
\begin{equation} \label{e_def_hamiltonian}
	H_{R_n}(h_{R_n}, \omega) = \sum_{e \in E(R_n)} \omega_{h_{R_n}(e)},
\end{equation}
where~$E(R_n) = \{e_{x,y} \,|\, x, y \in R_n\}$ is the edge set
of the subgraph of $\Z^m$ induced by $R_n$.
\end{definition}

\begin{definition}[Quenched Gibbs measure]
Given a realization~$\omega$ of the random field
and a set $A \subset M(R_n)$ of height functions,
the partition function $Z_\omega(A)$ is given by
\begin{equation} \label{e_d_partition_function}
	Z_\omega(A)
	= \sum_{h_{R_n} \in A} \exp \bigl( H_{R_n}(h_{R_n}, \omega) \bigr) \,.
\end{equation}
For a fixed boundary data function $h_{\partial R_n} \in M(\partial R_n)$,
the quenched Gibbs measure $\mu_\omega$ on~$M(R_n, h_{\partial R_n})$
is defined by
\begin{equation}\label{e_d_gibbs_measure}
	\mu_{\omega}(h_{R_n})
	= \frac{1}{Z_\omega \bigl( M(R_n,h_{\partial R_n}) \bigr) }
	\exp \bigl( H_{R_n}(h_{R_n}, \omega) \bigr) \,.
\end{equation}
\end{definition}

\begin{remark}\label{r_gibbs_measure_uniform_measure}
If one chooses the constant field $\omega= \mathbf{0} = (0)_{e \in E(\Z)}$,
then the associated quenched Gibbs measure~$\mu_\mathbf{0}$
is the uniform measure on $M(R_n, h_{\partial R_n})$.
In this case one recovers the variational principle of~\cite{KrMeTa19}.
\end{remark}

Now let us introduce the microscopic entropy of our model.
Again there are two situations: first, the quenched case, defined for a fixed
realization~$\omega$ and the annealed case.

\begin{definition}[Quenched and annealed microscopic entropy]
\label{d_micro_entropy}
Given a domain $R_n \subset \Z^m$
and a finite non-empty subset $A \subset M(R_n)$,
the quenched microscopic entropy~$\Ent_{R_n}(A, \omega)$ is given by
\begin{align}\label{e_d_que_micro_entropy}
	\Ent_{R_n}(A, \omega)
	&:= - \frac{1}{\abs{R_n}} \log Z_\omega(A) \\
	& \biggl( = -\frac{1}{\abs{R_n}} \log \sum_{h_{R_n} \in A}
		\exp \bigl( H_{R_n}(h_{R_n}, \omega) \bigr) \biggr) \,.
\end{align}
The annealed microscopic entropy~$\Ent( R_n, h_{\partial R_n}) $ is given by
\begin{align}\label{e_d_ann_micro_entropy}
	\Ent_{R_n, \an}(A)
	:= \mathbb{E} \bigl[ \Ent_{R_n}(A, \omega) \bigr].
\end{align}
\end{definition}

\begin{remark}\label{r_uniform_micro_entropy}
As in Remark~\ref{r_gibbs_measure_uniform_measure},
if one chooses the constant field~$\omega = \mathbf{0}$,
then the quenched microscopic entropy
$\Ent_{R_n}(M(R_n, h_{\partial R_n}), \mathbf{0})$
is the same as the microscopic entropy of~\cite{KrMeTa19}.
\end{remark}

Next, we define the local surface tension.
As with the microscopic entropy, the local surface tension
admits both a quenched and an annealed version.

\begin{definition}[Quenched microscopic and local surface tension]
\label{d_que_local_surf_tens}
The quenched local surface tension is the a.s.-limit
\begin{equation} \label{e_que_local_surf_tens}
	\ent(s, \omega) := \lim_{n \to \infty} \ent_n(s, \omega) \,,
\end{equation}
where $\ent_n(s, \omega)$ is the quenched microscopic surface tension,
defined by
\begin{align}\label{e_que_micro_surf_tens}
	\ent_n (s, \omega)
	: = \Ent_{S_n} \bigl( M(S_n, h_{\partial S_n}^s), \omega \bigr)
	\,.
\end{align}

Recall from Notation and Conventions above that $S_n = \{-n, \dotsc, n\}^m$,
and note that the existence of the limit in~\eqref{e_que_local_surf_tens}
is the content of Lemma~\ref{l_que_local_surface_tension}.
\end{definition}

\begin{definition}[Annealed microscopic and local surface tension]
\label{d_annealed_local_surface_tension}
The annealed microscopic surface tension~$\ent_{n, \an}(s)$
is given by
\begin{align}
  \ent_{n, \an}(s) := \mathbb{E} \left[ \ent_n(s, \omega) \right],
\end{align}
and the annealed local surface tension~$\ent(s)$ is given by
\begin{align}
  \ent_\an(s) :=  \mathbb{E} \left[ \ent (s, \omega) \right].
\end{align}
\end{definition}

\begin{remark} \label{r_local_surface_tension}
Similarly to Remark~\ref{r_gibbs_measure_uniform_measure}
and Remark~\ref{r_uniform_micro_entropy},
we obtain back the local surface tension for the uniform measure
if we consider a constant random field~$\omega = \mathbf{0}$.
In the case of random potential,
it follows from Assumption~\ref{a_random_field}
and Lemma~\ref{p_entropy_bounds} that $\ent_n(s,\omega)$ is uniformly integrable
and therefore that $\ent_{n,\an}$ and $\ent_\an$ are well-defined.
\end{remark}

\begin{remark}
It is not hard to see that the annealed local surface tension
is also the limit of the annealed microscopic surface tension.
Indeed, from Assumption~\ref{a_random_field}
the quenched microscopic surface tension $\ent_n(s, \omega)$
is dominated by an $L^1$ function (see Lemma~\ref{p_entropy_bounds}).
Therefore, the dominated convergence theorem implies that
\begin{align}
	\lim_{n \to \infty} \ent_{n, \an}(s)
	= \lim_{n \to \infty} \mathbb{E} \left[ \ent_n(s, \omega) \right]
	= \mathbb{E} \left[ \lim_{n \to \infty} \ent_n(s, \omega) \right]
	= \ent_\an(s) \,.
\end{align}
\end{remark}

\medskip

The annealed macroscopic entropy is defined by:

\begin{definition}[Annealed macroscopic entropy]
\label{d_ann_macro_ent}
Given an asymptotic height function~$h_R \in M(R, h_{\partial R})$,
the annealed macroscopic entropy~$\Ent_{R,\an}(h_R)$
is defined by
\begin{equation}
	\Ent_{R,\an}(h_R) :=
	\int_R \ent_{\an} (\nabla h(x)) \, dx \,.
\end{equation}
\end{definition}

The first main result of this article is the profile theorem:

\begin{theorem}[Profile theorem] \label{th_profile}
Recall that $C_\omega : =1 \vee \sup_{e \in E(\Z)} \abs{\omega_e}$
is by Assumption~\ref{a_random_field} an $L^1$ random variable.
Then for any $h_R \in M(R, h_{\partial R})$ and any $\eta > 0$,
there exist functions $\theta_{h_R}(\delta)$
and $\theta_{h_R,\delta}(\tfrac{1}{n})$
with $\theta_{h_R}(\delta) \to 0$ as $\delta \to 0$
and $\theta_{h_R,\delta}(\tfrac{1}{n}) \to 0$ as $n \to \infty$
such that
\begin{equation} \label{e_profile_thm_ineq} \begin{aligned}
	\lim_{n \to \infty} \mathbb{P} \biggl(
		&\abs[\Big]{
			\Ent_{R_n} \bigl( B(R_n, h_R, \delta), \omega \bigr)
			- \Ent_\an(R, h_R) } \\
		&\qquad\ge \eta
		+ C_\omega \theta_{h_R}(\delta)
		+ C_\omega \theta_{h_R,\delta} \bigl( \tfrac{1}{n} \bigr)
	\biggr) = 0 \,.
\end{aligned} \end{equation}
\end{theorem}

The second main result is the variational principle:

\begin{theorem}[Variational principle] \label{th_varprin}
The random variables
\begin{equation}
	\Ent_{R_n}(M(R_n, h_{\partial R_n}, \delta), \omega)
\end{equation}
converge in probability to the infimum of $\Ent_\an(R, h_R)$
over asymptotic height functions $h_R \in M(R, h_{\partial R})$,
i.e.\@ for every~$\eta > 0$,

\begin{equation} \label{e_main_result_minimizer_homogenization} \begin{aligned}
	\limsup_{\delta \to 0} \, \limsup_{n \to \infty} \,
	\mathbb{P} \biggl( \, \abs[\Big] {
	&
		\Ent_{R_n} \bigl( M(R_n, h_{\partial R_n}), \omega \bigr) \\
	& \quad
		- \inf_{h_R \in M(R, h_{\partial R})} \Ent_{R,\an}(h_R) }
		\ge \eta
	\biggr)
	= 0 \,.
\end{aligned} \end{equation}
\end{theorem}

The third main result, which we state but do not prove,
is the large deviations principle.
The notation introduced below is standard for large deviations theory.
\begin{theorem}[Large deviations principle]
Consider the space $M(R)$ of asymptotic height functions on $R$,
endowed with the topology of uniform convergence.
For $\delta > 0$ and $n \in \N$,
define a random probability measure $\mu_{\delta,n}(\cdot, \omega)$ on $M(R)$ by
\begin{equation} \label{th_ldp}
	\mu_{\delta,n}(A, \omega)
	:= \frac{
		Z_\omega\bigl( \bigl\{
			h_{R_n} \in M(R_n, h_{\partial R_n}, \delta)
			\,\big|\, \tilde h_{R_n} \in A \bigr\} \bigr)
	}{
		Z_\omega\bigl( M(R_n, h_{\partial R_n}, \delta) \bigr)
	} \,,
\end{equation}
where $\tilde h_{R_n} \in M(R)$ denotes the asymptotic height function
given by rescaling and interpolating $h_{R_n} \in M(R_n)$,
i.e.\ $\tilde h_{R_n}(\tfrac{1}{n} z) = \frac{1}{n} h_{R_n}(z)$
for $z \in R_n$.

Then the measures $\mu_{\delta,n}$ satisfy a large deviations principle
in probability with rate functional $I$ given by
\[
	I(h_R) := \begin{cases}
		\Ent_{R,\an}(h_R) - E
		&\text{if $h_R \in M(R,h_{\partial R})$} \,, \\
		+\infty &\text{otherwise} \,,
	\end{cases}
\]
where $E := \inf_{h_R \in M(R, h_{\partial R})} \Ent_{R,\an}(h_R)$.
Specifically, this means that for any Borel set $A \subset M(R)$,
\begin{equation}
	\limsup_{\delta \to 0} \, \limsup_{n \to \infty} \, \P \biggl(
		\frac{1}{\abs{R_n}} \log \mu_{\delta,n}(A)
		\:\ge\: -\inf_{h_R \in A^\circ} I(h_R)
	\biggr) = 0
\end{equation}
and
\begin{equation}
	\limsup_{\delta \to 0} \, \limsup_{n \to \infty} \, \P \biggl(
		\frac{1}{\abs{R_n}} \log \mu_{\delta,n}(A)
		\:\le\: -\inf_{h_R \in \overline{A}} I(h_R)
	\biggr) = 0 \,,
\end{equation}
where $A^\circ$ denotes the interior of $A$
and $\overline{A}$ denotes the closure.
\end{theorem}

\section{The quenched and annealed local surface tension}
\label{s_local_surface_tension}

The purpose of this section is
to establish several fundamental properties of the
quenched entropy and local surface tension of our model.
We proceed as follows:

\begin{itemize}
\item In Section~\ref{ss_kirszbraun} we state the Kirszbraun theorem, used
heavily in the rest of this section and beyond.
\item In Section~\ref{ss_robustness} we derive robustness of the entropy and
local surface tension under boundary value changes.
\item In Section~\ref{ss_lst_exist} we prove the existence of the quenched local
surface tension and the equivalence between the quenched and annealed local
surface tension.
\item In Section~\ref{ss_lst_other} we study the local surface tension
as a function $s \mapsto \ent_\an(s)$,
and we show that this function is convex and continuous.
\end{itemize}

\subsection{Kirszbraun theorem}
\label{ss_kirszbraun}

The Kirszbraun theorem for $\Z$-homomorphisms is a discrete analogue
of the classical Kirszbraun theorem of~\cite{Kir34}.
The classical theorem gives a condition under which
a Lipschitz continuous function can be extended
from a subset of a domain to the entirety of that domain.
Likewise, the Kirszbraun theorem for graph homomorphisms
gives a condition under which a $\Z$-valued graph homomorphism may be extended
from a subset of a domain to the entire domain.
Note that the property of being a $\Z$-valued graph homomorphism is stronger
than the Lipschitz property with constant $1$,
since if $z \sim \tilde z$ are two adjacent points in the domain
of a graph homomorphism $h: S \to \Z$, then $h(z) \ne h(\tilde z)$.

\begin{theorem} \label{p_kirszbraun}
Let~$\Lambda$ be a connected region of $\Z^m$,
let~$S$ be a subset of~$\Lambda$,
and let~$\bar{h}: S \to \Z$ be a graph homomorphism that preserves parity.
There exists a graph homomorphism~$h: \Lambda \to \Z$
such that~$h = \bar{h}$ on~$S$ if and only if for all~$x, y \in S$,
\begin{equation} \label{e_kirszbraun}
	d_\Z ( \bar{h}(x), \bar{h}(y) ) \le d_\Lambda(x, y),
\end{equation}
where~$d_\Z$ and~$d_\Lambda$ denote respectively the graph distance
on~$\Z$ and on~$\Lambda \subset \Z^m$.
\end{theorem}

This is a well-known result (see e.g.\ \cite[Lemma 4.3.1]{She05}),
and we omit the proof from this article.
As an illustration of the usefulness of the Kirszbraun theorem,
we prove the following lemma,
which justifies the choice of the normalizing factor $\frac{1}{\abs{R_n}}$
in Definition~\ref{d_micro_entropy}:

\begin{lemma} \label{p_entropy_bounds}
Almost surely
(in terms of the distribution~$\mathbb{P}$ of the random field~$ \omega$),
\begin{equation}
	-\log(2) - 2m C_\omega
	\le \Ent_{R_n} \bigl( M(R_n, h_{\partial R_n}), \omega \bigr)
	\le 2m C_\omega \,.
\end{equation}
\end{lemma}

\begin{proof}
As a corollary of the Kirszbraun theorem
(Theorem~\ref{p_kirszbraun}),
there is always at least one height function~$h_0 \in M(R_n, h_{\partial R_n})$.
So,
\begin{align}
	\Ent_{R_n} \bigl( M(R_n, h_{\partial R_n}), \omega \bigr)
	&\le -\frac{1}{\abs{R_n}} \log \sum_{h \in \{h_0\}} \exp \left(
		\sum_{e \in E(R_n)} \omega_{e_{h(x),h(y)}} \right) \\
	&\le \frac{\abs{E(R_n)}}{\abs{R_n}} C_\omega \\
	&\le 2 m C_\omega \,.
\end{align}

On the other hand, we overestimate the cardinality of~$M(R_n, h_{\partial R_n})$
as follows: enumerate the points of the interior of~$R_n$,
in such a way that each point~$x_i$ is adjacent to the previous point~$x_{i-1}$
(and the first point~$x_1$ is adjacent to~$x_0 \in \partial R_n$).
For each point~$x_i$ in the enumeration,
we require that~$h(x_i) = h(x_{i-1}) \pm 1$,
so there are at most~$2$ choices for~$h(x_i)$.
All together, $\abs{ M(R_n, h_{\partial R_n}) } \le 2^{\abs{R_n}}$.
It follows that
\begin{align}
	\Ent(R_n, h_{\partial R_n}, \omega)
	&\ge - \frac{1}{\abs{R_n}} \log \Bigl(
		\abs{ M(R_n, h_{\partial R_n}) }
		\exp \bigl( C_\omega \abs{ E(R_n) } \bigr)
	\Bigr) \\
	&\ge -\frac{1}{\abs{ R_n }} \log 2^{\abs{ R_n }}
		- \frac{\abs{ E(R_n) }}{\abs{ R_n }} C_\omega \\
	&\ge -\log(2) - m C_\omega.
\end{align}
\end{proof}

In the sequel, we will usually use the Kirszbraun theorem
in the following setting.
Given two domains $R_{n_1} \subset R_{n_2} \subset \Z^m$,
a height function $h_{R_{n_1}} \in M(R_{n_1})$,
and a boundary height function $h_{\partial R_{n_2}} \in M(\partial R_{n_2})$,
there exists an extension $\tilde h_{R_{n_2}} \in M(R_{n_2})$
with $\tilde h_{R_{n_2}}|_{R_{n_1}} = h_{R_{n_1}}$
and $\tilde h_{R_{n_2}}|_{\partial R_{n_2}} = h_{\partial R_{n_2}}$
if and only if
\[
	\abs{ h_{R_{n_1}}(z_1) - h_{\partial R_{n_2}}(z_2) }
	\le \abs{ z_1 - z_2 }_1
	\quad
	\text{for all $z_1 \in \partial R_{n_1}, z_2 \in \partial R_{n_2}$} \,.
\]

\subsection{Robustness of the quenched entropy}
\label{ss_robustness}

The quenched microscopic entropy and local surface tensions are robust,
in the sense that small changes in boundary values
cause small changes in the numeric value of the entropy.
There are two steps in proving these robustness results:
First, just as for the unperturbed model of~\cite{KrMeTa19},
compare the two sets of height functions
associated with the two boundary value functions,
perhaps by exhibiting an injection from one set into the second
or by estimating cardinalities directly.
Second, show that individual height functions from each of the two sets
contribute comparable amounts to the entropy
after applying the random potential,
e.g.\ by showing that every height function in one set
admits a ``similar'' height function in the second set,
whose Hamiltonian value is not much different;
this step is sometimes straightforward and other times quite subtle.

\begin{lemma} \label{p_equivalence_of_quenched_ent_free_fixed}
Let $\alpha > 0$,
let $s \in \R^m$ with $\abs{s}_\infty \le 1-\alpha$,
let~$\varepsilon \in (0, \tfrac{\alpha}{2})$,
let~$n \in \N$ with $n \ge (1 - \tfrac{2\varepsilon}{\alpha})^{-1}$),
and let $h_{\partial S_n} \in M(\partial S_n, s, \varepsilon)$.
Write
\[
	n^+ := \ceil[\big]{(1 + \tfrac{2\varepsilon}{\alpha}) n}
	\qquad \text{and} \qquad
	n^- := \floor[\big]{(1 - \tfrac{2\varepsilon}{\alpha}) n} \,.
\]
(We remark that $1 \le n^- < n < n^+$.) Then,
\begin{equation} \begin{aligned} \label{e_robustbulk_concl}
	\ent_{n^+}(s, \omega)
	- C_\omega \, \theta \bigl( \tfrac{\varepsilon}{\alpha} \bigr)
	&\le
	\Ent_{S_n} \bigl( M(S_n, h_{\partial S_n}), \omega \bigr) \\
	&\qquad \le
	\ent_{n^-}(s, \omega)
	+ C_\omega \, \theta_m \bigl( \tfrac{\varepsilon}{\alpha} \bigr) \,.
\end{aligned} \end{equation}
\end{lemma}

\begin{figure}
	\centering
	\begin{tikzpicture}[
		x=0.01in,y=0.01in,
	]
		\newcommand*\smallrad{65}
		\newcommand*\medrad{95}
		\newcommand*\bigrad{125}
		\path[draw] (-\smallrad,-\smallrad)
			node [anchor=south west]{$S_{n^-}$}
			-- (\smallrad,-\smallrad) -- (\smallrad,\smallrad)
			-- (-\smallrad,\smallrad) -- cycle;
		\path[draw] (-\medrad,-\medrad)
			node [anchor=south west]{$S_n$}
			-- (\medrad,-\medrad) -- (\medrad,\medrad)
			-- (-\medrad,\medrad) -- cycle;
		\path[draw] (-\bigrad,-\bigrad)
			node [anchor=south west] {$S_{n^+}$}
			-- (\bigrad,-\bigrad) -- (\bigrad,\bigrad)
			-- (-\bigrad,\bigrad) -- cycle;
		\path[draw,decorate,decoration=brace]
			(\bigrad,-\medrad) + (0.05in,0) --
			node [anchor=west] {$\delta n$}
			+ (0.05in,-30);
	\end{tikzpicture}
	\caption{Nested domains from
	Lemma~\ref{p_equivalence_of_quenched_ent_free_fixed}.
	}
	\label{f_ent_kirsz_domains}
\end{figure}

\begin{proof}[Proof of Lemma~\ref{p_equivalence_of_quenched_ent_free_fixed}]
We prove the inequality
\[
	\ent_{n^+}(s, \omega)
	- C_\omega \theta \bigl( \tfrac{\varepsilon}{\alpha} \bigr)
	\le \Ent_{S_n}\bigl( M(S_n, h_{\partial S_n}), \omega \bigr).
\]
The proof of the reverse inequality is similar.
\\

Note that the smaller square~$S_n = \{-n, -(n-1), \dotsc, n-1, n\}^m$
is contained inside the larger square~$S_{n^+}$,
and that
\begin{equation} \label{e_robustbulk_bdy}
	\abs{ x-y }_1 \ge \tfrac{2\varepsilon}{\alpha} n
	\qquad \text{whenever $x \in \partial S_n$
		and $y \in \partial S_{n^+}$} \,.
\end{equation}
\\

We construct an injection from~$M(S_n, h_{\partial S_n})$
into~$M(S_{n^+}, h_{\partial S_{n^+}}^s)$
using the Kirszbraun theorem, Theorem~\ref{p_kirszbraun}.
Let $h_{S_n} \in M(S_n, h_{\partial S_n})$,
let $x \in \partial S_n$, and let $y \in \partial S_{n^+}$.
By the definitions of $M(S_n, h_{\partial S_n})$
and of $h_{\partial S_n}^s$,
\begin{align} \label{e_equivalence_quenched_ent_free_fixed_kirszbraun}
       \hskip3em & \hskip-3em
       \abs[\big]{ h_{S_n}(x) - h_{S_{n^+}}^s(y) } \\
       &\le \abs[\big]{ h_{S_n}(x) - s \cdot x }
       \,+\, \abs[\big]{ s \cdot (x - y) }
       \,+\, \abs[\big]{ h_{S_{n^+}}^s(y) - s \cdot y } \\
       &\le \varepsilon n + \abs{s}_\infty \abs{ x-y }_1 + 1.
\end{align}
By hypothesis $\abs{s}_\infty \le 1 - \alpha$
and by~\eqref{e_robustbulk_bdy},
$\varepsilon n \le \tfrac{\alpha}{2} \abs{ x-y }_1$.
Therefore for $n \ge \tfrac{2}{\alpha}$,
\[
	\abs{ h_{S_n}(x) - h_{S_{n^+}}^s(y) } \le \abs{ x-y }_1 \,,
\]
so $h_{S_n}$ admits an extension
$h_{S_{n^+}} \in M(S_{n^+}, h_{\partial S_{n^+}}^s)$.
The map $h_{S_n} \mapsto h_{S_{n^+}}$ is an injection
from~$M(S_n, h_{\partial S_n})$ into~$M(S_{n^+}, h_{\partial S_{n^+}}^s)$.
The existence of such an injection implies immediately that
\begin{equation} \begin{aligned}
	\Ent_{S_n} \bigl( M(S_n, h_{\partial S_n}), \omega \bigr)
	&\ge \frac{\abs{S_n}}{\abs{S_{n^+}}} \Ent_{S_{n^+}} \bigl(
		M(S_{n^+}, h_{\partial S_{n^+}}^s), \omega \bigr) \\
	&\qquad- \frac{2mC_\omega (\abs{S_{n^+}} - \abs{S_n})}{\abs{S_n}} \\
	&= \Ent_{S_{n^+}} \bigl( M(S_{n^+}, h_{\partial S_{n^+}}^s),
		\omega \bigr)
	- C_\omega \theta_m \bigl( \tfrac{\varepsilon}{\alpha} \bigr) \,.
\end{aligned} \end{equation}

This proves the first inequality of~\eqref{e_robustbulk_concl}.
As mentioned at the beginning of the proof, the other inequality is similar.
Since $n^- < n$, one extends height functions
from~$M(S_{n^-}, h_{\partial S_{n^-}^s})$ to~$M(S_n, h_{\partial S_n})$.
We omit the details.
\end{proof}

Lemma~\ref{p_equivalence_of_quenched_ent_free_fixed}
does not extend to the case where~$\abs{s}_\infty = 1$.
As~$\abs{s}_\infty \to 1$ the ratio of the box sizes
$\frac{\abs{S_{n^+}}}{\abs{S_n}} \approx 1 + \frac{\varepsilon}{\alpha}$
and the error bound~$\theta(\tfrac{\varepsilon}{\alpha})$ both diverge.
Fundamentally these difficulties come from the Kirszbraun theorem.
When~$\abs{s}_\infty$ is close to~$1$,
the ``margin'' $S_{n^+} \setminus S_n$ must be large in order to connect
$h_{\partial S_n}$ to $h_{\partial S_{n^+}}$
and when $\abs{s}_\infty = 1$, such an extension is not generally possible.
Therefore we take a different approach for $\abs{s}_\infty \approx 1$,
using elementary combinatorics to count the number of height functions.
The two following calculations are intermediate results
used to prove the robustness lemma, Lemma~\ref{p_continuity_of_que_ent_slope_1}.

\begin{lemma}[Counting height functions near $\abs{s}_\infty = 1$]
\label{calc_comboedge}
Let $\varepsilon > 0$.
Let~$s \in \R^m$ with $1 - \varepsilon < \abs{s}_\infty \le 1$,
and let~$h_{\partial S_n} \in M(\partial S_n, h_{\partial S_n}^s, \varepsilon)$.
Then,
\begin{equation} \label{e_comboedge_concl}
	\frac{1}{\abs{S_n}} \log \abs[\big]{ M(S_n, h_{\partial S_n}) }
	= \theta(\varepsilon) \,.
\end{equation}
\end{lemma}

\begin{proof}[Proof of Lemma~\ref{calc_comboedge}]
Fix a coordinate index $1 \le i \le m$
such that $\abs{s_i} > 1 - \varepsilon$,
and assume without loss of generality that $s_i > 1 - \varepsilon$.
Decompose $S_n$ into $(2n+1)^{m-1}$ lines
in the $i$\textsuperscript{th} coordinate direction.
Along each such line $h_{S_n}$ must increase
by at least $2(1-2\varepsilon)n$.
Therefore, the $2n$ edges in the line split into two subsets:
at least $2(1-2\varepsilon)n$ ``increasing'' edges,
and at most $4 \varepsilon n$ ``decreasing'' edges.
Counting each line independently, we conclude that
\[
	\abs[\big]{ M(S_n, h_{\partial S_n}) }
	\le \binom{2n}{\ceil{4 \varepsilon n}}^{(2n+1)^{m-1}} \,.
\]
The conclusion~\eqref{e_comboedge_concl} follows immediately.
For a more verbose version of this proof,
see \cite[Lemma 21]{KrMeTa19}.
\end{proof}

\begin{lemma}[Height functions at slope $\abs{s}_\infty = 1$]
\label{calc_combo1}
Let $s' \in \R^m$ with $\abs{s'}_\infty = 1$.
Then $\abs{ M(S_n, h_{\partial S_n}^{s'}) } = 1$,
and the sole element of $M(S_n, h_{\partial S_n}^{s'})$
is the canonical height function $h_{S_n}^{s'}$.
\end{lemma}

\begin{proof}[Proof of Lemma~\ref{calc_combo1}]
As in the proof of Lemma~\ref{calc_comboedge},
fix a coordinate index $1 \le i \le m$ such that $\abs{s_i} = 1$.
Decompose $S_n$ into lines in the $i$\textsuperscript{th} coordinate direction.
Along each line, any height function $h_{S_n} \in M(S_n, h_{\partial S_n}^{s'})$
must increase by exactly $2n$.
Since $h_{S_n}$ is a graph homomorphism,
that is only possible if $h_{S_n}$ increases along every edge,
i.e.\@ $h_{S_n}(x+1,y) - h_{S_n}(x,y) = 1$ for $x = -n, \dotsc, n-1$.
It follows that $\abs{ M(S_n h_{\partial S_n}^{s'}) } \le 1$.
To complete the proof, observe that
$h_{S_n}^{s'} \in M(S_n, h_{\partial S_n}^{s'})$.
\end{proof}

Having recorded Lemma~\ref{calc_comboedge} and~\ref{calc_combo1},
we return to establishing robustness results.
As in Lemma~\ref{p_equivalence_of_quenched_ent_free_fixed},
our goal is to compare the microscopic surface tension
$\ent_n(s, \omega) := \Ent_{S_n}(M(S_n, h_{\partial S_n}^s), \omega)$
and the entropy $\Ent_{S_n}(M(S_n, h_{\partial S_n}), \omega)$
associated to an ``approximately affine'' boundary height function
$h_{\partial S_n} \in M(\partial S_n, h_{\partial S_n}^s, \varepsilon)$.
The difference is that Lemma~\ref{p_equivalence_of_quenched_ent_free_fixed}
took $\abs{s}_\infty \le 1 - \alpha$
and the lemma below takes $\abs{s}_\infty > 1 - \alpha$.
\\

\begin{lemma} \label{p_continuity_of_que_ent_slope_1}
Let $\varepsilon > 0$.
Let~$s, s' \in \R^m$ with~$\abs{s}_\infty \le 1$, $\abs{s'}_\infty = 1$,
and $\abs{s - s'}_\infty < \varepsilon$.
Let $n \in \N$ be sufficiently large
(specifically, $n \ge \tfrac{1}{\varepsilon}$)
and let $h_{\partial S_n} \in M(\partial S_n, h_{\partial S_n}^s, \varepsilon)$.
Then:
\begin{equation} \label{e_robustedge_concl}
	\abs[\big]{ \Ent(S_n, h_{\partial S_n}, \omega) - \ent_n(s', \omega) }
	\,\le\, C_\omega \, \theta(\varepsilon) \,.
\end{equation}
\end{lemma}

Because of the~$\theta(\varepsilon)$ error term,
Lemma~\ref{p_continuity_of_que_ent_slope_1}
will not be useful for slopes~$s$ with~$\abs{s}_\infty$ far from~$1$.

\begin{remark}[Comment about the proof]
There are two ingredients to the proof.
The first is counting results
of Lemma~\ref{calc_comboedge} and Lemma~\ref{calc_combo1},
and the second is a comparison between
the Hamiltonian $H_{S_n}(h_{S_n},\omega)$
of a generic height function $h_{S_n} \in M(S_n, h_{\partial S_n})$
and the Hamiltonian $H_{S_n}(h_{S_n}^{s'}, \omega)$
of the unique element $h_{S_n}^{s'} \in M(S_n, h_{\partial S_n}^{s'})$.
Since proofs were already given for the two lemmas,
most of the argument below is spent on the comparison of Hamiltonians.
\\

The comparison of Hamiltonians is also fundamentally a combinatorial argument
that relies on the rigidity caused by the slopes $s$ and $s'$ being close to
(or on) the boundary of the slope space $[-1,1]^m$.
It is surprising that such a subtle argument is (apparently) needed
in the case of homogenization,
since the two counting lemmas are sufficient in the uniform case,
and these lemmas are not very complicated to prove.
\\

The subtlety is similar to that of the proof of Lemma~\ref{lem_robustbdy} below.
In both cases, the subtlety arises when comparing Hamiltonians
for two height functions defined on the same domain $S_n$.
In comparison, the proof Lemma~\ref{p_equivalence_of_quenched_ent_free_fixed}
(which has a similar statement to the current
Lemma~\ref{p_continuity_of_que_ent_slope_1}))
is based on extending height functions from one domain to another larger domain
via the Kirszbraun theorem.
Comparing the Hamiltonian of a height function on a large domain
to the Hamiltonian of the same function on a restricted domain is simple,
since the difference is exactly relatable to the difference in domains.
\end{remark}

\begin{proof}[Proof of Lemma~\ref{p_continuity_of_que_ent_slope_1}]
As mentioned above, we will compare the Hamiltonians
$H_{S_n}(h_{S_n}, \omega)$ and $H_{S_n}(h_{S_n}^{s'}, \omega)$,
where $h_{S_n} \in M(S_n, h_{\partial S_n})$
and $h_{S_n}^{s'} \in M(S_n, h_{\partial S_n}^{s'})$.
More precisely, we will later deduce the inequality
\begin{equation} \label{e_robustedge_hamil}
	\abs[\big]{ H_{S_n}(h_{S_n}, \omega) - H_{S_n}(h_{S_n}^{s'}, \omega) }
	\le 210 m^2 (2n+1)^m C_\omega \varepsilon \,.
\end{equation}

Given that~\eqref{e_robustedge_hamil} holds, the proof is straight-forward:
For one inequality, we calculate
\begin{align}
	\hskip3em&\hskip-3em
	\Ent_{S_n} \bigl( M(S_n, h_{\partial S_n}), \omega \bigr) \\
	&= -\frac{1}{\abs{S_n}} \log \sum_{h_{S_n} \in M(S_n, h_{\partial S_n})}
		\exp \bigl( H_{S_n}(h_{S_n}, \omega) \bigr) \\
	&\overset{\mathclap{\eqref{e_robustedge_hamil}}}{\le}
	-\frac{1}{\abs{S_n}} \log \sum_{h_{S_n} \in M(S_n, h_{\partial S_n})}
		\exp \biggl( H_{S_n}(h_{S_n}^{s'}, \omega) \\
	&\hspace{16em}
			- 210 m^2 (2n+1)^m C_\omega \varepsilon \biggr) \\
	&\overset{\mathclap{Lemma~\ref{calc_comboedge}}}{\le} \qquad
	-\frac{1}{\abs{S_n}} H_{S_n} \bigl( h_{S_n}^{s'}, \omega \bigr)
	+ \theta(\varepsilon)
	+ 210 m^2 C_\omega \varepsilon \\
	\label{e_robustedge_ent}
	&= \Ent_{S_n} \bigl( M(S_n, h_{S_n}^{s'}), \omega \bigr)
	+ \theta(\varepsilon) \,.
\end{align}

The opposite inequality is derived in the same way,
which concludes the proof of Lemma~\ref{p_continuity_of_que_ent_slope_1}
up to the verification of~\eqref{e_robustedge_hamil}.\\

For convenience, let us use for the remaining argument the following convention:
When denoting the Hamiltonian of~$H(h_{S_n}, \omega)$
we just write~$H(h_{S_n})$,
omitting the dependency on the random field~$\omega$. \\

Verification of~\eqref{e_robustedge_hamil}:
Heuristically, the estimate~\eqref{e_robustedge_hamil} makes sense.
Because the slopes~$s$ and~$s'$ are ~$\varepsilon$-close to each other,
and~$s'$ has slope 1,
every height function~$h_{S_n} \in M(S_n, h_{\partial S_n})$
has to behave similar to the canonical height function
$h_{S_n}^{s'}$ of slope $s'$.
Therefore, the difference in the associated energies,
as measured by the Hamiltonian~$H_{S_n}(h_{S_n})$ and $H_{S_n}(h_{S_n}^{s'})$,
should vanish as~$\varepsilon \to 0$.\medskip

To make this argument rigorous one needs to precisely estimate
the number of heights that each height function $h_{S_n}$ visits,
i.e.\ the set $\{h_{S_n}(e) \,|\, e \in E(S_n)\}$ with multiplicities,
and compare to the corresponding set for $h_{S_n}^{s'}$.
This is relatively straight-forward on a one-dimensional lattice
but unfortunately becomes much more subtle on a higher-dimensional lattice.
To see why, consider the decomposition of the box $S_n$ into lines.
This leads a decomposition of the edges in $E(S_n)$ into
\emph{parallel edges} within a line,
and \emph{cross edges} connecting two lines.
Without cross edges the one-dimensional argument would easily carry over,
but controlling the cross edges is necessary as well.
This control is accomplished by the sets $G_y$ below.\\

To begin the rigorous verification of~\eqref{e_robustedge_hamil},
pick an arbitrary height function $h_{S_n} \in M(S_n, h_{\partial S_n})$.
As mentioned above, we decompose $S_n$ into lines
parallel to one of the coordinate axes.
Assume by symmetry that
$s = (s_1, s_2, \dotsc, s_m)$ and $s' = (s_1', \dotsc, s_m')$
satisfy $s'_1 = 1$ and (therefore) $s_1 > 1 - \varepsilon$.
For $y \in \{-n,\dotsc,n\}^{m-1}$
let $\ell_y$ denote the line in the first coordinate direction
through $(0,y)$ in $S_n$, i.e.
\[
	\ell_y
	:= \bigl\{ (-n,y),\, (-n+1,y),\, \dotsc,\, (n-1,y),\, (n,y) \bigr\} \,.
\]
Observe that $S_n$ is the disjoint union of the $(2n+1)^{m-1}$ lines $\ell_y$.
In particular, the Hamiltonian $H_{S_n}(h_{S_n})$
decomposes with respect to the lines $\ell_y$ as
\begin{equation} \begin{aligned} \llabel{e_tilde_hamil}
	H_{S_n}(h_{S_n})
	&:= \sum_{e \in E(S_n)} \omega_{h_{S_n}(e)} \\
	&= \sum_y \Biggl(
		\sum_{e \in E(\ell_y)} \omega_{h_{S_n}(e)}
		+ \frac{1}{2} \sum_{y' \sim y} \, \sum_{e \in \tilde E_{y,y'}} \omega_{h_{S_n}(e)}
	\Biggr) \\
	&= \sum_y \tilde H_y(h_{S_n}) \,,
\end{aligned} \end{equation}
where $\tilde E_{y,y'}$ is the set of edges in $E(S_n)$
with one endpoint in $\ell_y$ and the other in $\ell_{y'}$
(we call these \emph{cross edges}),
and where $\tilde H_y$ is defined to be the parenthesized quantity
from the line above.
Note that the factor $\tfrac{1}{2}$ is necessary
because each cross edge in $\tilde E_{y,y'}$
also contributes to $\tilde H_{y'}(h_{S_n})$,
so without the factor $\tfrac{1}{2}$
the contributions from the cross edges would be double-counted.
\\

We define two families of sets $A_y \subset E(\Z)$ and $G_y \subset A_y$,
indexed by points $y \in \{-n, \dotsc, n\}^{m-1}$.
In terms of the heuristic argument above,
these sets roughly correspond to the heights visited
by $h_{S_n}$ and $h_{S_n}^{s'}$,
although in fact both $A_y$ and $G_y$ are subsets of
$\{h_{S_n}(e) \,|\, e \in E(S_n)\}$.
\\

Let $A_y$ denote the edges $e \in E(\Z)$ that lie inside the interval
from $(s \cdot (-n, y) + 2 \varepsilon n)$
to $(s \cdot (+n, y) - 2 \varepsilon n)$.
Based on the boundary conditions and homomorphism property
of $h_{S_n}$ and $h_{S_n}^{s'}$,
every edge $e \in A_y$ occurs
both in the image $\{h_{S_n}(\tilde e) \,|\, \tilde e \in E(\ell_y)\}$
and in the image $\{h_{S_n}^{s'}(\tilde e) \,|\, \tilde e \in E(\ell_y)\}$.
(The factors of $2$ in the definition of $A_y$ are necessary
since the boundary height function $h_{\partial S_n}$
may differ from $h_{\partial S_n}^s$ by up to $\varepsilon n$,
in addition to $s_1$ differing from $1$ by up to $\varepsilon$.)
The situation in dimension $m=1$ is illustrated in
Figure~\ref{ff_que_ent_cts_edge_ay}.
\\

\begin{figure}
	\subcaptionbox{
		Here $h_{S_n}$ is a one-dimensional height function
		with slope $s \ge 1 - \varepsilon$.
		The set $A_y$ comprises the $1-4\varepsilon$ fraction
		of the $2n$ edges in $\ell_y$, centered around 0.
		(The central height in higher dimensions
		is instead $s \cdot (0,y)$.)
		Both $h_{S_n}$ and $h_{S_n}^{s'}$
		must contain all of these edges in their image.
		They might contain additional edges.
		\label{ff_que_ent_cts_edge_ay}
	}[0.45\textwidth]{ \begin{tikzpicture}[
		x=0.007in,y=0.007in,
	]
		\newcommand*\ptrad{2pt}
		\draw[<->] (-100,0) node[left]{$x$} -- (100,0);
		\draw[<->] (0,-100) -- (0,100) node[above]{$h_{S_n}(x)$};
		\draw	   (-100,-90)
			-- (- 32,-22)	
			-- (- 28,-26)	
			-- (  28, 30)	
			-- (  32, 26)	
			-- (  49, 43)	
			-- (  51, 41)	
			-- ( 100, 90);	
		\fill[black] (-100,-90) circle [radius=\ptrad]
			node [below] {$(-n, -sn)$};
		\fill[black] ( 100, 90) circle [radius=\ptrad]
			node [above] {$(n, sn)$};
		\draw[decorate,decoration=brace] (100,80) + (0.05in,0)
			-- node [anchor=west] {$A_y$}
			+(0.05in,-160);
	\end{tikzpicture} }
	\hfill
	\subcaptionbox{
		The three lines are $\ell_y$ in the center
		and two of its neighbors, $\ell_{y'}$ and $\ell_{y''}$.
		The highlighted edge is the edge $e_s \in E(\ell_y)$
		for $e \in G_y$,
		i.e.\ the unique edge in $\ell_y$ with $h_{S_n}(e_s) = e$.
		There is also an edge $e_{s'}$ (not shown),
		satisfying the corresponding uniqueness property
		for $h_{S_n}^{s'}$.
		Finally, all six highlighted vertices are good,
		i.e.\ each vertex has a unique height within its line.
		\label{ff_que_ent_cts_edge_gy}
	}[0.45\textwidth]{ \begin{tikzpicture}[
		x=0.07in,y=0.07in,
	]
		\newcommand*\ptrad{2pt}
		\draw[<->] (0,-10) node[below]{$\ell_y$} -- (0,10);
		\draw[ultra thick] (0,-3) -- (0,3);
		\path (0,0) node [right] {$e_s$};
		\draw[<->] (6,-10) node[below]{$\ell_{y''}$} -- (6,10);
		\draw[<->] (-6,-10) node[below]{$\ell_{y'}$} -- (-6,10);
		\fill[black]	( 0, -3) circle [radius=\ptrad]
				( 0,  3) circle [radius=\ptrad]
				( 6, -3) circle [radius=\ptrad]
				( 6,  3) circle [radius=\ptrad]
				(-6, -3) circle [radius=\ptrad]
				(-6,  3) circle [radius=\ptrad];
	\end{tikzpicture} }
	\caption{
		Figures relating to the proof
		of Lemma~\ref{p_continuity_of_que_ent_slope_1}.
	}
	\label{f_que_ent_cts_edge_defs}
\end{figure}

We define $G_y \subset A_y$ in the following way:
These are the edges $e \in A_y \subset E(\Z)$
satisfying these three constraints with respect to $h_{S_n}$
(illustrated in Figure~\ref{ff_que_ent_cts_edge_gy}):

\begin{itemize}
\item $e$ occurs with multiplicity $1$ in the multi-set
	$\{h_{S_n}(\tilde e) \,|\, \tilde e \in E(\ell_y)\}$.
	(By choice of $A_y$, $e$ occurs with multiplicity $\ge 1$.)
	Write $e_s$ for the unique edge $e_s \in E(\ell_y)$ such that
	$h_{S_n}(e_s) = e$.
\item Both endpoints of $e$ occur with multiplicity $1$
	in the multi-set $\{h_{S_n}(z) \,|\, z \in \ell_y\}$.
\item For each endpoint $z$ of $e_s$ and each neighboring vertex $z' \sim z$
	that lies in $S_n \setminus \ell_y$,
	$h_{S_n}(z')$ occurs with multiplicity $1$
	in the multi-set $\{h_{S_n}(\tilde z) \,|\, \tilde z \in \ell_{y'}\}$
	for the line $\ell_{y'}$ that contains it.
\end{itemize}
Further on in the argument, we will call elements of~$G_y$ ``good'' edges.
We will call a vertex $z \in \ell_y$ ``good''
if its height $h_{S_n}(z)$ occurs in with multiplicity $1$ in
$\{h_{S_n}(\tilde z) \,|\, \tilde z \in \ell_y\}$,
and likewise for $z' \in \ell_{y'}$.
\\

Later on, we will need that for an arbitrary ``good'' edge~$e \in G_y$ it holds:
\begin{align} \label{e_good_edge_counting}
    \sum_{y' \sim y} \,
    \sum_{\substack{\tilde e \in \tilde E_{y,y'}\\ h_{S_n}(\tilde e)=e}} \,
        \omega_{h_{S_n} (\tilde e)}
    = \abs[\big]{\{y' \sim y\}} \, \omega_e .
\end{align}
Note that $\abs{\{y' \sim y\}} \le 2m$ for all $y$,
with equality unless $y$ is a boundary point
(implicitly we assume that $y' \in \{-n,\dotsc,n\}^{m-1}$).
Argument for~\eqref{e_good_edge_counting}:
We observe that for each $y' \sim y$,
by using the second and third constraints and considering cases,
there is a unique cross edge $e_{s,y'}$ between $\ell_y$ and $\ell_{y'}$
such that $h_{S_n}(e_{s,y'}) = e$.
For a proof of this simple fact we refer to Figure~\ref{f_que_ent_cts_edge_cases}.
The identity~\eqref{e_good_edge_counting} follows then immediately.\\

\begin{figure}
	\subcaptionbox{
		Case 1 (both adjacent height values larger):
		Clearly there is one edge between $\ell_y$ and $\ell_{y'}$
		that is mapped to $e = e_{k,k+1}$.
		Suppose that another cross edge has heights $k$ and $k+1$.
		Then its left endpoint would have
		either height $k$ or height $k+1$,
		which contradicts the fact that
		the two labelled vertices in $\ell_y$ are ``good,''
		i.e. that their heights occur only once in $\ell_y$.
		\label{ff_que_ent_cts_case_1}
	}[0.3\textwidth]{ \begin{tikzpicture}[
		x=0.05in,y=0.07in,
	]
		\newcommand*\ptrad{2pt}
		\draw[<->] (0,-7) node[below]{$\ell_y$} -- (0,11);
		\draw[ultra thick] (0,-3) -- (0,3);
		\draw[ultra thick] (0,-3) -- (6,-3);
		\draw[<->] (6,-7) node[below]{$\ell_{y'}$} -- (6,11);
		\fill[black]
			( 0, -3) circle [radius=\ptrad] node[left] {$k$}
			( 0,  3) circle [radius=\ptrad] node[left] {$k+1$}
			( 6, -3) circle [radius=\ptrad] node[right]{$k+1$}
			( 6,  3) circle [radius=\ptrad] node[right]{$k+2$};
	\end{tikzpicture} }
	\hfill
	\subcaptionbox{
		Case 2 (both adjacent height values smaller):
		Again there is one edge between $\ell_y$ and $\ell_{y'}$
		that is mapped to $e = e_{k,k+1}$,
		and again no other vertices in $\ell_y$
		can have either height $k$ or height $k+1$.
		\label{ff_que_ent_cts_case_2}
	}[0.3\textwidth]{ \begin{tikzpicture}[
		x=0.05in,y=0.07in,
	]
		\newcommand*\ptrad{2pt}
		\draw[<->] (0,-7) node[below]{$\ell_y$} -- (0,11);
		\draw[ultra thick] (0,-3) -- (0,3);
		\draw[ultra thick] (0,3) -- (6,3);
		\draw[<->] (6,-7) node[below]{$\ell_{y'}$} -- (6,11);
		\fill[black]
			( 0, -3) circle [radius=\ptrad] node[left] {$k$}
			( 0,  3) circle [radius=\ptrad] node[left] {$k+1$}
			( 6, -3) circle [radius=\ptrad] node[right]{$k-1$}
			( 6,  3) circle [radius=\ptrad] node[right]{$k$};
	\end{tikzpicture} }
	\hfill
	\subcaptionbox{
		Case 3 (cannot occur because $e \in G_y$):
		Here there would be two edges between the lines
		that both map to $e_{k,k+1}$.
		But since the vertex at height $k$ in $\ell_y$ is ``good'',
		the vertex labelled $\alpha$ must have height $k+2$.
		Likewise since the vertex at height $k+1$ in $\ell_{y'}$
		is ``good'', vertex $\beta$ must have height $k-1$.
		Since $\alpha \sim \beta$, this violates the
		graph homomorphism property.
		\label{ff_que_ent_cts_case_3}
	}[0.3\textwidth]{ \begin{tikzpicture}[
		x=0.05in,y=0.07in,
	]
		\newcommand*\ptrad{2pt}
		\draw[<->] (0,-7) node[below]{$\ell_y$} -- (0,11);
		\draw[ultra thick] (0,-3) -- (0,3);
		\draw[ultra thick] (0,-3) -- (6,-3);
		\draw[ultra thick] (0,3) -- (6,3);
		\draw[<->] (6,-7) node[below]{$\ell_{y'}$} -- (6,11);
		\fill[black]
			( 0, -3) circle [radius=\ptrad] node[left] {$k$}
			( 0,  3) circle [radius=\ptrad] node[left] {$k+1$}
			( 0,  9) circle [radius=\ptrad] node[left] {$\alpha$}
			( 6, -3) circle [radius=\ptrad] node[right]{$k+1$}
			( 6,  3) circle [radius=\ptrad] node[right]{$k$}
			( 6,  9) circle [radius=\ptrad] node[right]{$\beta$};
	\end{tikzpicture} }
	\caption{
		Consideration of cases for part of the proof
		of Lemma~\ref{p_continuity_of_que_ent_slope_1}.
		The claim to be shown is: given $e \in G_y$
		(say $e = e_{k,k+1}$),
		there is a unique cross edge $e_{s,y'} \in \tilde E_{y,y'}$
		which is mapped to $e$ by the height function $h_{S_n}$.
		In the figure, the vertices are labelled by their heights,
		i.e.\ by the values of $h_{S_n}$.
		The bolded edge in $\ell_y$ is $e_s \in E(\ell_y)$,
		i.e.\ the unique edge in $\ell_y$ with $h_{S_n}(e_s) = e$.
		In Figure~\ref{ff_que_ent_cts_case_1}
		and Figure~\ref{ff_que_ent_cts_case_2},
		the bolded edge between the lines is the unique
		edge between the lines with height $e_{k,k+1}$.
		Figure~\ref{ff_que_ent_cts_case_3} shows two such edges,
		but in fact this case cannot occur.
		By the homomorphism property,
		these three cases exhaust the possibilities for heights
		on the two vertices in $\ell_{y'}$ that are adjacent to
		the endpoints of $e_s$.}
	\label{f_que_ent_cts_edge_cases}
\end{figure}

We will also need count $\abs{G_y}$.
Heuristically, since the slope $s$ is close to $1$,
$G_y$ must be a large subset of $E(\ell_y)$.
To be precise, recall that
$\abs{A_y} \ge 2n - 4\ceil{\varepsilon n}$ by construction,
and that $G_y$ is the subset of edges $e \in A_y$
that satisfy the three constraints above.
The second constraint actually implies the first,
so to count $G_y$ we simply count how many edges in $A_y$
satisfy the last two constraints.
Actually we count the complement,
i.e.\ how many edges do not satisfy these two constraints.
Indeed, each ``bad'' vertex in $\ell_y$
(in the sense described after the constraints)
causes at most two edges in $E(\ell_y)$ to violate the second constraint.
Likewise, each ``bad'' vertex in an adjacent line $\ell_y'$
causes at most two edges in $E(\ell_y)$ to violate the second constraint.
All other edges in $A_y$ are ``good,'' i.e. are included in $G_y$.

It remains to count the ``bad'' vertices in any line $\ell_y$.
Since $s_1 > 1 - \varepsilon$ and since $h_{S_n}$
approximates the slope-$s$ height function $h_{S_n}^s$ on $\partial S_n$,
the height values $h_{S_n}(-n,y)$ and $h_{S_n}(+n,y)$
on the endpoints of $\ell_y$ differ by at least $2 n - 4 \varepsilon n$.
Since $h_{S_n}$ is a graph homomorphism, it maps
the $2n+1$ vertices in $\ell_y$
surjectively onto the set of $\ge 2 n - 4 \ceil{\varepsilon n} + 1$ integers
between the heights of the endpoints.
By the pigeonhole principle, at most $8 \ceil{\varepsilon n}$ of these integers
occur with multiplicity $\ge 2$,
i.e.\ at most $8 \ceil{\varepsilon n}$ vertices are ``bad.''
Thus
\begin{align}
	\abs{G_y}
	&\ge \abs{A_y} - 2 \, \abs[\big]{ \bigl\{
		\text{``bad'' vertices in $\ell_y$
		or $\ell_{y'}$ (for $y' \sim y$)} \bigr\} } \\
	\:&\ge\:
	\underbrace{2n - 4 \ceil{\varepsilon n}}_{\abs{A_y}}
	\,-\, 2
	\: \cdot \:
	\underbrace{(2m+1)}_{\text{\# lines}}
	\quad \cdot \quad
	\underbrace{8 \ceil{\varepsilon n}}_{\mathclap{\text{
			``bad'' vertices per line}
	}} \\
	&= 2n - (32m + 20) \ceil{\varepsilon n} \\
	\llabel{e_size_gy}
	&\ge 2n - 52m\ceil{\varepsilon n} \,.
\end{align}

Now we work towards the Hamiltonian estimate~\eqref{e_robustedge_hamil}.
Let $e \in G_y$,
and recall that $e_s$ is the unique edge in $E(\ell_y)$
such that $h_{S_n}(e_s) = e$,
and that $e_{s,y'}$ is the unique cross edge
between $\ell_y$ and $\ell_{y'}$ such that $h_{S_n}(e_{s,y'}) = e$.
As a result (recall the definitions of $\tilde H_{\ell_y}$ and $\tilde E_{y,y'}$
from~\leqref{e_tilde_hamil} above):
\begin{equation} \begin{aligned}
	\tilde H_{\ell_y} ( h_{S_n} )
	&= \Bigl( \sum_{\tilde e \in E(\ell_y)}
		\omega_{h_{S_n}(\tilde e)} \Bigr)
	+ \frac{1}{2} \Bigl( \sum_{y' \sim y} \, \sum_{\tilde e \in \tilde E_y}
		\omega_{h_{S_n}(\tilde e)} \Bigr)
	\\
	& \overset{\mathclap{\eqref{e_good_edge_counting}}}{=} \Bigl(
		\sum_{e \in G_y} \omega_e
		+ \sum_{\substack{
			\tilde e \in E(\ell_y) \\
			h_{S_n}(\tilde e) \not\in G_y
		}} \omega_{h_{S_n}(\tilde e)} \Bigr) \\
	&\qquad + \frac{1}{2} \Bigr(
		\abs[\big]{\{y' \sim y\}} \sum_{e \in G_y} \omega_e
		+ \sum_{y' \sim y} \, \sum_{\substack{
			\tilde e \in \tilde E_{y,y'} \\
			h_{S_n}(\tilde e) \not\in G_y
		}} \omega_{h_{S_n}(\tilde e)} \Bigr) \,,
\end{aligned} \end{equation}
so
\begin{align}
	\hskip3em&\hskip-3em
	\abs[\Big]{
		\tilde H_{\ell_y}(h_{S_n})
		- \bigl( \tfrac{1}{2} \abs[\big]{\{y' \sim y\}} + 1 \bigr)
			\sum_{e \in G_y} \omega_e
	} \\
	&\le C_\omega \bigl( \abs{E(\ell_y)} - \abs{G_y} \bigr)
	+ \frac{1}{2} \sum_{y' \sim y} C_\omega \bigl(
		\abs{\tilde E_{y,y'}} - \abs{G_y} \bigr) \\
	&\overset{\mathclap{\leqref{e_size_gy}}}{\le}
	52m C_\omega \ceil{\varepsilon n}
	+ \frac{1}{2} \sum_{y' \sim y} C_\omega \bigl(
		52m \ceil{\varepsilon n} + 1 \bigr) \\
	&\le 52m C_\omega \ceil{\varepsilon n} (1+m) + m C_\omega \\
	&\le 104m^2 C_\omega \ceil{\varepsilon n} + m C_\omega \\
	&\le 105 m^2 (2n+1) C_\omega \varepsilon \,.
\end{align}

(In the last line, we assume that $n \ge \tfrac{1}{\varepsilon}$,
so that $(2n+1) \varepsilon \ge \ceil{\varepsilon n} \ge 1$.)
\\

Because $s'_1 = 1$, $h_{S_n}^{s'}|_{\ell_y}$ is an injection,
the three bullet points above are also satisfied with $h_{S_n}^{s'}$
in place of $h_{S_n}$.
Therefore the calculation above also applies with $h_{S_n}^{'s}$
in place of $h_{S_n}$, so
\begin{equation}
	\abs[\Big]{
		\tilde H_{\ell_y}(h_{S_n}^{s'})
		- \bigl( \tfrac{1}{2} \abs[\big]{\{y' \sim y\}} + 1 \bigr)
			\sum_{e \in G_y} \omega_e
	}
	\le 105 m^2 (2n+1) C_\omega \varepsilon \,.
\end{equation}

By the triangle inequality,
\begin{equation}
	\abs[\big]{ \tilde H_{\ell_y}(h_{S_n})
		- \tilde H_{\ell_y}(h_{S_n}^{s'}) }
	\le 210 m^2 (2n+1) C_\omega \varepsilon \,.
\end{equation}
By summing over $y \in \{-n,\cdots,n\}^{m-1}$,
we get the desired inequality~\eqref{e_robustedge_hamil}, i.e.
\begin{equation}
	\abs[\big]{ H_{S_n}(h_{S_n}) - H_{S_n}(h_{S_n}^{s'}) }
	\le 210  m^2 (2n+1)^m C_\omega \varepsilon \,.
\end{equation}
\end{proof}

Both Lemma~\ref{p_equivalence_of_quenched_ent_free_fixed}
and Lemma~\ref{p_continuity_of_que_ent_slope_1}
imply that the microscopic entropy is robust to changes in boundary data,
but they apply in different regimes.
The former result applies when the boundary data has slope~$s$
with norm~$\abs{s}_\infty$ bounded away from~$1$,
and the latter when the slope~$s$ has norm close to~$1$.
For convenience later on, we combine the two results into a single theorem.

\begin{theorem} \label{p_robustness}
For any~$\varepsilon \in (0,\tfrac{1}{9})$ and any slope~$s \in [-1,1]^m$,
there exist~$A = A(s, \varepsilon) > 0$, $B = B(s, \varepsilon) > 0$,
and~$n_0 = \ceil{\tfrac{1}{\varepsilon}} \in \N$ such that,
for any~$n \ge n_0$ and any boundary height function~%
$h_{\partial S_n} \in M(\partial S_n, h_{\partial S_n}^s, \varepsilon)$,
\begin{equation} \begin{aligned}
	\ent_{An}(s, \omega) - C_\omega \, \theta(\varepsilon)
	&\le \Ent \bigl( M(S_n, h_{\partial S_n}), \omega \bigr) \\
	&\qquad \le \ent_{Bn}(s, \omega) + C_\omega \, \theta(\varepsilon).
\end{aligned} \end{equation}
Moreover, the functions~$A(s, \varepsilon)$ and~$B(s, \varepsilon)$
are bounded away from~$0$ and~$\infty$ uniformly in~$s$ and~$\varepsilon$.
More precisely,
\begin{equation}
	1
	\le A(s, \varepsilon)
	\le \bigl( 1 + 2\varepsilon^{1/2} + \tfrac{1}{n} \bigr)
	< \infty
\end{equation}
and
\begin{equation}
	0
	< \bigl( 1 - 2\varepsilon^{1/2} - \tfrac{1}{n} \bigr)
	< B(s, \varepsilon)
	\le 1.
\end{equation}
\end{theorem}

\begin{proof}[Proof of Theorem~\ref{p_robustness}]
Take~$\alpha = \varepsilon^{1/2}$, and proceed according to two cases.
For slopes~$s$ with~$\abs{s}_\infty \le 1-\alpha$,
use Lemma~\ref{p_equivalence_of_quenched_ent_free_fixed}
to choose~$A = n^+ / n \approx (1 + 2\varepsilon^{1/2})$
and~$B = n^- / n \approx (1 - 2\varepsilon^{1/2})$.
Note that $\varepsilon < \tfrac{1}{9}$
implies that $\varepsilon < \tfrac{\alpha}{2}$
and $n \ge \tfrac{1}{\varepsilon} \ge (1 - 2\varepsilon^{1/2})^{-1}$,
as required by the lemma.
Moreover $1 - 2\varepsilon^{1/2} - \tfrac{1}{n} > \tfrac{2}{9}$,
so $B$ is indeed bounded away from $0$.
The error terms~$\theta(\tfrac{\varepsilon}{\alpha})$ from the lemma
are equivalent to~$\theta(\varepsilon^{1/2}) = \theta(\varepsilon)$.
\\

For slopes with~$\abs{s}_\infty > 1-\alpha$, take~$A = B = 1$
and apply Lemma~\ref{p_continuity_of_que_ent_slope_1} twice,
using $\alpha = \varepsilon^{1/2}$ in place of $\varepsilon$:
once for the boundary height function~$h_{\partial S_n}$ given in the statement
of the theorem,
and once for the canonical boundary height function~$h_{\partial S_n}^s$.
The estimate on
\[
	\abs[\big]{
		\Ent(M(S_n, h_{\partial S_n}), \omega) - \ent_n(s, \omega)
	}
\]
follows from the triangle inequality.
\end{proof}

The robustness results above focused on boundary height functions
that differed at macroscopic scale,
i.e.\ $\abs{h_{\partial S_n} - \tilde h_{\partial S_n}}_u \le \varepsilon n$.
For boundary height functions with sub-linear differences,
we will derive stronger robustness results.
Lemma~\ref{lem_robustbdy} addresses the case where the two
boundary height functions differ at only a single point on $\partial S_n$,
and Corollary~\ref{cor_robustbdy} extends to the sub-linear case
(actually, only to
$\abs{h_{\partial S_n} - \tilde h_{\partial S_n}}_u = o(\tfrac{\log n}{n})$,
but that is sufficient for our purposes.)

\begin{lemma}[Robustness for minimally different boundary height functions]
\label{lem_robustbdy}
Fix $n \in \N$,
and let $h_{\partial S_n}^+$ and $\tilde h_{\partial S_n}^-$
be two boundary height functions on the hypercube $S_n$
which differ at exactly one point $z_0 \in \partial S_n$,
i.e.\ $
	h_{\partial S_n}^+|_{S_n \setminus \{z_0\}}
	= h_{\partial S_n}^-|_{S_n \setminus \{z_0\}}
$ and $h_{\partial S_n}^+(z_0) = h_{\partial S_n}^-(z_0) + 2$.

Then,
\begin{equation} \label{e_robustbdy_concl}
	\abs[\Big]{
		\Ent_{S_n} \bigl(
			M(S_n, h_{\partial S_n}^+), \omega \bigr)
		- \Ent_{S_n} \bigl(
			M(S_n, h_{\partial S_n}^-), \omega \bigr)
	}
	\le \frac{4mC_\omega + \log(2n)}{\abs{S_n}} \,.
\end{equation}
\end{lemma}

\begin{remark}
The $\log(2n)$ term is necessary at least in some extreme cases.
For example, suppose that $\omega \equiv 0$, $m=1$, $z_0 = -n$,
$h_{\partial S_n}^+(-n) = 2$, $h_{\partial S_n}^-(-n) = 0$,
and $h_{\partial S_n}^\pm(n) = 2n$. Then
$\Ent_{S_n}(M(S_n, h_{\partial S_n}^+),\mathbf{0}) = -\frac{1}{n} \log (2n)$
and $\Ent_{S_n}(M(S_n, h_{\partial S_n}^-),\mathbf{0}) = 0$;
cf.\ Lemma~\ref{calc_comboedge} and Lemma~\ref{calc_combo1} for calculations.
\end{remark}

\begin{proof}[Proof of Lemma~\ref{lem_robustbdy}]
For concreteness and w.l.o.g., assume that the boundary values at $z_0$ are
$h_{\partial S_n}^-(z_0) = 0$ and $h_{\partial S_n}^+(z_0) = 2$.
(Technically this assumption is only valid if $z_0$ has even parity
because we require that height functions preserve parity,
and one should instead assume e.g. that $h_{\partial S_n}^\pm(z_0) \in \{1,3\}$
in the other case.
For simplicity we ignore this detail in the rest of the proof.)
\\

Consider the line $z_0, z_1, \dotsc, z_{2n}$ of points in $S_n$
starting from $z_0$ and going into $S_n$, perpendicular to the boundary.
Classify each height function $h_{S_n}^+ \in M(S_n, h_{\partial S_n}^+)$
based on the number of initial ``up'' steps, i.e.
\[
	k_\up(h_{S_n}^+)
	:= \max \bigl \{ \tilde k \ge 0 \,|\,
		h_{S_n}^+(z_k) = h_{S_n}^+(z_{k-1}) + 1 \,
		\text{for $1 \le k \le \tilde k$} \bigr\} \,.
\]

Note that from our initial assumption, $h_{S_n}^+(z_k) = k+2$
for $0 \le k \le k_\up$.
Necessarily $k_\up(h_{S_n}^+) < 2n$,
since if $h_{S_n}^+$ went up along all $2n$ edges,
then the values $h_{S_n}^-(z_{2n}) = h_{S_n}^+(z_{2n}) = 2n+2$
and $h_{S_n}^-(z_0) = 0$ would violate the Kirszbraun theorem.
\\

On the line segment $\{z_0, \dotsc, z_{n_\texttt{up}}\} \subset S_n$,
$h_{S_n}^+$ is ``too high,''
in the sense that no height function
in $M(S_n, h_{\partial S_n}^-)$ can match it.
But by the Kirszbraun theorem,
there exists $h_{S_n}^- \in M(S_n, h_{\partial S_n}^-)$
such that $h_{S_n}^-(z_{n_\texttt{up}+1}) = h_{S_n}^+(z_{n_\texttt{up}+1})$.
In fact, we may define $h_{S_n}^-$ by
\[
	h_{S_n}^-(z) = \begin{cases}
		k = h_{S_n}^+(z)-2, &
		\text{if $z = z_k$ for $0 \le k \le k_\up$, and} \\
		h_{S_n}^+(z), &
		\text{otherwise} \,.
	\end{cases}
\]

It follows that $h_{S_n}^+$ and $h_{S_n}^-$ have the same Hamiltonian,
except for the contribution from the edges
incident to a vertex $z_k$ ($0 \le k \le k_\up$).
There are $(2m-1) (k_\up+1)$ such edges,
which leads to the naive estimate
$
	\abs{H_{S_n}(h_{S_n}^+,\omega) - H_{S_n}(h_{S_n}^-,\omega)}
	\le (2m-1) (k_\up+1) C_\omega
$.
This estimate is not useful because $k_\up$ on the right-hand side
leads to an error of order $n$ in the worst case.
However, as shown in Figure~\ref{f_1point},
a more careful estimate is possible.
Indeed, both $h_{S_n}^+$ and $h_{S_n}^-$ map the edges $e$ in question
to the same collection of edges
$\{e_{k,k+1} \,|\, 0 \le k \le k_\up\} \subset E(\Z)$,
with each $e_{k,k+1}$ repeated about $2m-1$ times.
We omit the details, but a careful count of the edge heights
yields the inequality
\begin{equation} \llabel{e_hamil}
	\abs[\big]{H_{S_n}(h_{S_n}^+,\omega) - H_{S_n}(h_{S_n}^-,\omega)}
	\le 4mC_\omega \,.
\end{equation}

\begin{figure}
	\subcaptionbox{
		The values of the height functions $h_{S_n}^+$ and $h_{S_n}^-$
		from the proof of Lemma~\ref{lem_robustbdy}.
		On the vertices $z_0, \dotsc, z_{k_\up}$
		where the two height functions differ,
		the larger value is the height that $h_{S_n}^+$ takes
		and the smaller value is $h_{S_n}^-$.
		Here $k_\up = 3$, since $h_{S_n}^+$ increases across
		the first three edges in the center line.
		The $(2m-1)(k_\up+1)$ shaded edges are exactly the set up edges
		incident to any of $z_0,\dotsc,z_{k_\up}$,
		and these are the only edges on which $h_{S_n}^\pm$ differ.
		\label{f_1point_pic}
	}[0.45\textwidth]{

\begin{tikzpicture}[
	x=0.5in,
	y=0.5in,
	every node/.style={minimum size=0.35in,font=\tiny,inner sep=3pt},
]
	\draw (-1,0) node [circle      ,draw] (l0) {$1$};
	\draw ( 0,0) node [circle split,draw] (m0) {$2$ \nodepart{lower} $0$};
	\node [below right=0.05in and -0.03in of m0,minimum size=0, inner sep=0]
		{$z_0$};
	\draw ( 1,0) node [circle      ,draw] (r0) {$1$};
	\draw  (l0.west) -- +(-0.1in,0) (r0.east) -- +(+0.1in,0);
	\draw[ultra thick] (l0.east) -- (m0.west) (m0.east) -- (r0.west);
	\foreach \pre/\lo/\adj/\hi/\lbl in {%
		0/1/2/3/$z_1$,%
		1/2/3/4/$z_2$,%
		2/3/4/5/$z_3 = z_{k_\up}$%
	} {
		\draw (-1,\lo) node [circle      ,draw] (l\lo) {$\adj$};
		\draw ( 0,\lo) node [circle split,draw] (m\lo)
			{$\hi$ \nodepart{lower} $\lo$};
		\draw ( 1,\lo) node [circle      ,draw] (r\lo) {$\adj$};
		\draw[ultra thick] (l\lo) -- (m\lo) (m\lo) -- (r\lo)
			(m\pre) -- (m\lo);
		\draw (l\pre) -- (l\lo) (r\pre) -- (r\lo);
		\draw (l\lo.west) -- +(-0.1in,0) (r\lo.east) -- +(0.1in,0);
		\node [below right=0.05in and -0.03in of m\lo,minimum size=0,
			inner sep=0,fill=white] {\lbl};
	}
	\draw (-1,4) node [circle      ,draw] (l4) {$?$};
	\draw ( 0,4) node [circle      ,draw] (m4) {$4$};
	\node [below right=0.05in and -0.03in of m4,minimum size=0, inner sep=0]
		{$z_4$};
	\draw ( 1,4) node [circle      ,draw] (r4) {$?$};
	\draw[ultra thick] (m3) -- (m4);
	\draw (l4) -- (m4) (m4) -- (r4) (l3) -- (l4) (r3) -- (r4);
	\draw (l4.west) -- +(-0.1in,0) (r4.east) -- +(0.1in,0)
		(l4.north) -- +(0,0.1in) (m4.north) -- +(0,0.1in)
		(r4.north) -- +(0,0.1in);
\end{tikzpicture}
	}
	\hfill
	\subcaptionbox{
		Number of shaded edges on which $h_{S_n}^+$, $h_{S_n}^-$
		attain certain heights.
		For example, from the last row of the table:
		$h_{S_n}^+(e) = e_{k_\up+1,k_\up+2}$
		for all $2m$ edges incident on $z_{k_\up}$.
		In the difference $H_{S_n}(h_{S_n}^+) - H_{S_n}(h_{S_n}^-)$,
		the bulk of the height values in the table cancel,
		leaving only boundary terms.
		That is why the bound in~\leqref{e_hamil}
		does not depend on $k_\up$.
		\label{f_1point_table}
	}[0.45\textwidth]{
		\begin{tabular}{|c|c|c|}
			\hline
			$e \in E(\Z)$ & $h_{S_n}^+$ & $h_{S_n}^-$ \\
			\hline
			$e_{0,1}$ & $0$    & $2m-1$ \\
			$e_{1,2}$ & $2m-2$ & $2m-1$ \\
			$e_{2,3}$ & $2m-1$ & $2m-1$ \\
			$e_{3,4}$ & $2m-1$ & $2m-1$ \\
			\vdots & \vdots & \vdots \\
			$e_{k_\up,k_\up+1}$ & $2m-1$ & $2m-1$ \\
			$e_{k_\up+1,k_\up+2}$ & $2m$ & $0$ \\
			\hline
		\end{tabular}
	}
	\caption{Explanation of inequality~\leqref{e_hamil}
	from the proof of Lemma~\ref{lem_robustbdy}.}
	\label{f_1point}
\end{figure}

Now we turn to the entropy inequality. For $0 \le k < 2n$, let
\[
	M_k := \bigl\{
		h_{S_n}^+ \in M(S_n, h_{\partial S_n}^+)
		\,\big|\, k_\up(h_{S_n}^+) = k
	\bigr\}.
\]
Then the sets $M_k$ ($0 \le k < 2n$) partition $M(S_n, h_{\partial S_n}^+)$,
so
\begin{align}
	\hskip3em&\hskip-3em
	\Ent_{S_n} \bigl( M(S_n, h_{\partial S_n}^+), \omega \bigr) \\
	&= - \frac{1}{\abs{S_n}} \log \sum_{k=0}^{2n-1} \,
		\sum_{h_{S_n}^+ \in M_k} \exp \bigl(
			H_{S_n}(h_{S_n}^+, \omega) \bigr) \\
	&\overset{\mathclap{\leqref{e_hamil}}}{\ge}
	- \frac{1}{\abs{S_n}} \log \sum_{k=0}^{2n-1} \,
		\sum_{h_{S_n}^+ \in M_k} \exp \bigl(
			H_{S_n}(h_{S_n}^-, \omega) + 4mC_\omega \bigr) \\
	&\ge - \frac{1}{\abs{S_n}} \log \sum_{k=0}^{2n-1} \,
		\sum_{h_{S_n}^- \in M(S_n, h_{\partial S_n}^-)} \exp \bigl(
			H_{S_n}(h_{S_n}^-, \omega) + 4mC_\omega \bigr) \\
	&= \Ent_{S_n} \bigl( M(S_n, h_{\partial S_n}^-), \omega \bigr)
		- \frac{4mC_\omega + \log(2n)}{\abs{S_n}} \,.
\end{align}

The reverse inequality is derived
by exchanging the roles of $h_{\partial S_n}^\pm$,
considering the number $k_\down$ of initial downward steps
of $h_{\partial S_n}^-$ on the line $\{z_0,\dotsc,z_{2n}\}$,
and proceeding as before with the necessary changes.
\end{proof}

Lemma~\ref{lem_robustbdy} applies only when the two boundary height function
$h_{S_n}^+$ and $h_{S_n}^-$ differ minimally.
However by applying Lemma~\ref{lem_robustbdy} repeatedly,
we can compare two more different height functions.
That idea is captured in the following corollary.

\begin{corollary}[Robustness with respect to sub-linear height differences]
\label{cor_robustbdy}
Let $h_{\partial S_n}$ and $\tilde h_{\partial S_n}$
be boundary height functions on $S_n$,
and let $M = \norm{h_{\partial S_n} - \tilde h_{\partial S_n}}_\infty$.
Then
\begin{equation} \begin{aligned}
	\hskip3em&\hskip-3em
	\abs[\Big]{
		\Ent_{S_n} \bigl(
			M(S_n, h_{\partial S_n}), \omega \bigr)
		- \Ent_{S_n} \bigl(
			M(S_n, \tilde h_{\partial S_n}), \omega \bigr)
	} \\
	&\le \frac{M}{2} \bigl( 4 m C_\omega +\log(2n) \bigr)
		\frac{\abs{\partial S_n}}{\abs{S_n}} \,.
\end{aligned} \end{equation}
\end{corollary}

\begin{remark}
The main idea of the proof is to interpolate the boundary height function
from $h_{\partial S_n}$ to $\tilde h_{\partial S_n}$,
where each step in the interpolation changes the value
of the boundary height function at exactly one boundary point.
Note that each interpolation step changes the height by $2$
at that distinguished boundary point,
which is the reason for the factor $\tfrac{M}{2}$ rather than simply $M$.
Given such an interpolation,
all that remains is to apply Lemma~\ref{lem_robustbdy}
and the triangle inequality.
\end{remark}

\begin{proof}[Proof of Corollary~\ref{cor_robustbdy}]
We claim that there exists a finite sequence
$h_{\partial S_n}^{(1)}, \dotsc, h_{\partial S_n}^{(k)}$
such that each pair $h_{\partial S_n}^{(j)}$ and $h_{\partial S_n}^{(j+1)}$
differ at exactly one point,
such that $h_{\partial S_n}^{(1)} = h_{\partial S_n}$
and $h_{\partial S_n}^{(k)} = \tilde h_{\partial S_n}$,
and such that $k \le \frac{M}{2} \abs{\partial S_n}$.
Each element of the sequence is constructed from the previous element
by a ``flip'' operation:
Given a (boundary) height function $h_{\partial S_n}^{(j)}$
and a vertex $z_j \in \partial S_n$
where all the neighboring vertices $z' \in \partial S_n, \, z' \sim z_j$
have the same height $h_{\partial S_n}(z') = a \in \Z$,
the height function $h_{\partial S_n}^{(j+1)}$
is identical to $h_{\partial S_n}^{(j)}$ on $\partial S_n \setminus \{z_j\}$
and takes the other valid value on $z_j$.
Specifically, if $h_{\partial S_n}^{(j)}(z_j) = a + 1$,
then $h_{\partial S_n}^{(j+1)}(z_j) = a - 1$;
otherwise $h_{\partial S_n}^{(j+1)}(z_j) = a + 1$.

It remains to show that the vertices $z_1, \dotsc, z_{k-1}$ can be chosen
so that $h_{\partial S_n}^{(k)} = \tilde h_{\partial S_n}$
and so that $k \le \frac{M}{2} \abs{\partial S_n}$.
To prove both these points, consider the metric
$d: M(\partial S_n) \times M(\partial S_n) \to \Z$
defined by
\[
	d(h_{\partial S_n}', h_{\partial S_n}'')
	:= \sum_{z \in \partial S_n} \abs[\big]{
		h_{\partial S_n}'(z) - h_{\partial S_n}''(z) } \,.
\]

As long as $d(h_{\partial S_n}^{(j)}, \tilde h_{\partial S_n}) > 0$,
we will find a vertex $z_j$ for which the flip operation both is valid
and decreases the distance $d$.
Towards this end, let $E_j := \{z \in \partial S_n
	\,|\, h_{\partial S_n}^{(j)}(z) > \tilde h_{\partial S_n}(z) \}$.
If $E_j \ne \varnothing$,
choose $z_j := \argmax_{z \in E_j} h_{\partial S_n}^{(j)}$.

We claim that flipping at $z_j$ is valid,
and more specifically that for all neighbors $z' \sim z_j$ in $\partial S_n$,
$h_{\partial S_n}^{(j)}(z') = h_{\partial S_n}^{(j)}(z_j)-1$.
Indeed, there are two cases.
If $h_{\partial S_n}^{(j)}(z') = \tilde h_{\partial S_n}(z')$
for any $z' \sim z_j$,
then necessarily $\tilde h_{\partial S_n}(z_j) = h_{\partial S_n}^{(j)}(z_j)-2$
and $\tilde h_{\partial S_n}(z') = h_{\partial S_n}^{(j)}(z')
	= h_{\partial S_n}^{(j)}(z_j)-1$ for all $z' \sim z$.
Otherwise all $z' \sim z$ are also in $E_j$, so the claim follows
since $z_j$ maximizes $h_{\partial S_n}^{(j)}$ over $E_j$.
So as claimed, it is valid to flip the height function $h_{\partial S_n}^{(j)}$
at $z_j$, and this flip decreases the difference
$\abs{h_{\partial S_n}^{(j+1)}(z_j) - \tilde h_{\partial S_n}(z_j)}$ by two,
and therefore decreases the distance
$d(h_{\partial S_n}^{(j+1)}, \tilde h_{\partial S_n})$ by two.

If $E_j$ is empty, use instead the set $F_j := \{z \in \partial S_n
	\,|\, h_{\partial S_n}^{(j)}(z) < \tilde h_{\partial S_n}(z) \}$,
pick $z_j := \argmin_{z \in F_j} h_{\partial S_n}^{(j)}$,
and repeat the argument, changing inequalities and signs accordingly.
If $F_j$ is also empty, then $h_{\partial S_n}^{(j)} = \tilde h_{\partial S_n}$
and the process is complete.

At most $\tfrac{1}{2} d(h_{\partial S_n}, \tilde h_{\partial S_n})
	\le \tfrac{M}{2} \abs{\partial S_n}$
steps are needed in total,
since each step decreases the distance by $2$.
\\

To complete the proof of the corollary, apply Lemma~\ref{lem_robustbdy}
to each pair $\{h_{\partial S_n}^{(j)}, h_{\partial S_n}^{(j+1)}\}$
and use the triangle inequality.
\end{proof}

\subsection{Existence and equivalence of quenched and annealed
local surface tension}
\label{ss_lst_exist}

Recall from Definition~\ref{d_que_local_surf_tens}
that the quenched local surface tension is defined as
the limit of the quenched microscopic surface tension.
Because of the random potential $\omega$,
the existence of this limit is not obvious.
We prove the existence of the limit
using an ergodic theorem for almost superadditive random families.\\

First, we introduce the notation needed for stating the ergodic theorem.
Let~$\mathcal{B}$ denote the set of all (non-empty) boxes in~$\Z^m$, i.e.
$$
	\mathcal{B} = \Bigl\{
		\bigl( \, [a_1, b_1) \times \dotsb \times [a_m, b_m) \bigr)
		\cap \Z^m
	\Bigm|
		a_1 < b_1, \, \dotsc, \, a_m < b_m \in \Z^m
	\Bigr\} \,.
$$
Note that the sets $S_n := [-n,n]^m \cap \Z^m$ are included in~$\mathcal{B}$.
We say that a family of~$L^1$ random variables~$F = (F_B)_{B \in \mathcal{B}}$
is \defn{almost superadditive} if,
for any finitely many disjoint boxes~$B_1, \dotsc, B_n \in \mathcal{B}$
whose union~$B = B_1 \cup \dotsb \cup B_n$ also lies in~$\mathcal{B}$,
\begin{equation} \label{e_almost_superadd}
	F_B \ge \sum_{i=1}^n F_{B_i} - A \sum_{i=1}^n \abs{\partial B_i}
    \quad \as \,,
\end{equation}
where~$A = A(\omega): \Omega \to [0, \infty)$ is an~$L^1$ random variable,
and where
$
	\partial B_i
	= \{ x \in B_i \,|\,  \exists y \in \Z^m \setminus B_i, \, x \sim y \}
$
is the inner boundary of~$B_i$.

\begin{theorem}[Ergodic theorem for almost superadditive random families]
\label{p_strong_superadd_erg_thm}
Let~$(\Omega, \mathcal{F}, \mathbb{P})$ be a probability space,
let~$\tau = (\tau_u)_{u \in \Z^m}$ be a family of measure-preserving
transformations on~$\Omega$,
and let~$F = (F_B)_{B \in \mathcal{B}}$ be a family of~$L^1$ random variables
satisfying the following three conditions:
\begin{itemize}
	\item
	$F$ is almost superadditive,
	i.e.~$F$ satisfies~\eqref{e_almost_superadd},
	\item
	For all~$u \in \Z^m$,
	\begin{equation} \label{e_almost_transinvar}
		\lim_{n \to \infty} \sup_{u \in \Z^m}
			\frac{1}{\abs{S_n}} \abs[\Big]{
				F_{u+S_n} - F_{S_n} \circ \tau_u}
		= 0 \,,
	\end{equation}
	where~$u+B = \{ u+x \,|\, x \in B \}$ is the translation of~$B$ by~$u$.
	\item
	The quantity~$
		\tilde\gamma(F)
		= \limsup_{n \to \infty} \frac{1}{\abs{S_n}} \,
			\mathbb{E}[F_{S_n}]
	$ is finite.
\end{itemize}
Then the limit~$\lim_{n \to \infty} \frac{1}{\abs{S_n}} \, F_{S_n}$
exists almost surely and in~$L^1$.
If moreover~$\{\tau_u\}_{u \in \Z^m}$ is ergodic, then the limit is
\begin{equation} \label{e_strong_superadd_erg_limit}
	\lim_{n \to \infty} \frac{1}{\abs{S_n}} \, F_{S_n}
	=
	\tilde\gamma(F).
\end{equation}
\end{theorem}

This theorem is based on \cite[Theorem~2.4]{AK81},
which is a multidimensional extension of the subadditive ergodic theorem
proven in~\cite{Kin68,Lig85} among many other sources.
The version stated here is adapted to notion of almost superadditivity
that the quenched microscopic entropy satisfies.
For completeness, we give a proof of this version of the ergodic theorem
in Appendix~\ref{a_almost_superadd_ergodic_thm}.
Now let us turn to the application of this ergodic theorem:
\\

\begin{lemma}[Existence of the quenched local surface tension]
  \label{l_que_local_surface_tension}
For almost every realization~$\omega$ of the random field,
the limit~\eqref{e_que_local_surf_tens} exists.
\end{lemma}

The proof is a straightforward application of the ergodic theorem.

\begin{proof}[Proof of Lemma~\ref{l_que_local_surface_tension}]
Fix $s \in [-1,1]^m$.
Let the family of measure-preserving transformations
$\tau = (\tau_u)_{u \in \Z^m}$
be given by
\begin{equation} \label{e_qlst_tau}
	\bigl( \tau_u \omega \bigr)_e
	:=
	\omega_{e - [s \cdot u]_{u \bmod 2}}
	\quad
	\text{for~$e \in E(\Z)$ and~$u \in \Z^m$}.
\end{equation}

Define the random process $F = (F_B)_{B \in \mathcal{B}}$ by
\begin{equation}
	F_B
	:= - \abs{B} \, \Ent \bigl( M(B, h_{\partial B}^s), \omega \bigr)
	= \log Z_\omega \bigl( M( B, h_{\partial B}^s) \bigr) \,.
\end{equation}

Now we verify the hypotheses of the ergodic theorem
(Theorem~\ref{p_strong_superadd_erg_thm}).
First, the fact that~$\abs{\omega_e} \le C_\omega$ for all edges~$e \in E(\Z)$
implies that each variable $F_B$ ($B \in \mathcal{B}$) is in $L^1$.
\\

Next, the almost superadditivity property~\eqref{e_almost_superadd}
follows from distributivity:
\begin{equation} \begin{aligned}
	\sum_{i=1}^n F_{B_i}
	&= \log \, \prod_{i=1}^n \,
		\sum_{h_{B_i} \in M(B_i, h_{\partial B_i}^s)}
			\exp \bigl( H_{B_i}(h_{B_i}, \omega) \bigr) \\
	&= \log \sum_{\substack{
		h_{B_1} \in M(B_1, h_{\partial B_1}^s) \\
		\dotsb \\
		h_{B_n} \in M(B_n, h_{\partial B_n}^s)
	}} \exp \Biggl( \, \sum_{i=1}^n H_{B_i} ( h_{B_i}, \omega ) \Biggr) \,.
\end{aligned} \end{equation}

The final sum is indexed by $n$-tuples of height functions,
i.e.\@ it is the sum over the Cartesian product
of the sets $M(B_i, h_{\partial B_i}^s)$.
This Cartesian product is a subset of $M(B, h_B)$, so
\begin{equation} \label{e_almost_superadd_ineq}
	\sum_{i=1}^n F_{B_i}
	\le \log \sum_{h_B \in M(B, h_{\partial B}^s)} \exp \Biggl( \,
		\sum_{i=1}^n H_{B_i} ( h_B|_{B_i}, \omega) \Biggr) \,.
\end{equation}

The quantity on the right-hand side of~\eqref{e_almost_superadd_ineq}
differs from $F_B$ by at most $m C_\omega \sum_{i=1}^n \abs{\partial B_i}$,
since the Hamiltonian terms in~\eqref{e_almost_superadd_ineq}
do not include edges that cross from one box $B_i$ to another box $B_j$.
This error term satisfies~\eqref{e_almost_superadd}.
\\

Now let us show that~$F$~satisfies
the translation invariance estimate~\eqref{e_almost_transinvar}.
For $h_{\partial (u+B)} \in M(\partial (u+B))$,
consider the shifted boundary height function
$\Psi_u h_{\partial (u+B)} \in M(\partial B)$ defined by
\begin{equation}
	(\Psi_u h_{\partial (u+B)})(z)
	:= h_{\partial (u+B)}(u+z) - \floor{s \cdot u}
	\qquad \text{for $z \in \partial B$} \,.
\end{equation}

Since both $h_{\partial B}^s$ and $h_{\partial (u+B)}^s$
are rounded to the nearest integer (of appropriate parity),
the shifted boundary height function $\Psi_u h_{\partial (u+B)}^s$
may not agree exactly with $h_{\partial B}^s$.
However it holds that
\[
	\abs[\big]{ \Psi_u h_{\partial (u+B)}(z) - h_{\partial B}(z) }
	\le 4
	\qquad \text{for all $z \in \partial B$} \,.
\]
Therefore by Corollary~\ref{cor_robustbdy},
\[
	\abs[\big]{ F_{u+B} - F_B \circ \tau_u }
	\le \abs{B} \, \theta \bigl( \tfrac{1}{n} \bigr) \,.
\]

The last condition to check is~
$
	\tilde\gamma(F)
	=
	\limsup_{n \to \infty} \frac{1}{\abs{S_n}} \mathbb{E}[F_{S_n}]
	< \infty
$,
which follows from boundedness of the quenched entropy.
Indeed by Lemma~\ref{p_entropy_bounds},
the inequality $F_B \le m \abs{B} C_\omega$ holds almost surely,
so $\tilde \gamma(F) \le \E(C_\omega) < \infty$.
\\

At this point we have checked all the hypotheses of the ergodic theorem
(Theorem~\ref{p_strong_superadd_erg_thm}).
From the ergodic theorem we conclude that the pointwise limit
\begin{equation}
	\ent(s, \omega)
	= \lim_{n \to \infty} \ent_n(s, \omega)
	= \lim_{n \to \infty} \frac{1}{\abs{S_n}} F_{S_n} (\omega)
\end{equation}
exists almost surely.
In addition, when $s \ne 0$, the family of measure-preserving transformations
$(\tau_u)_{u \in \Z^m}$ is ergodic with respect to $\P$,
since the family includes every shift
$\omega \mapsto (\omega_{k+e})_{e \in E(\Z)}$ for $k \in \Z$.
Therefore whenever $s \ne 0$, the limit $\ent(s,\omega)$
is almost surely equal to its expectation, $\E[\ent(s,\omega)] = \ent_\an(s)$.
\end{proof}

The failure of ergodicity in the case $s=0$ is evident
from the definition of $(\tau_u)_{u \in \Z^m}$ in~\eqref{e_qlst_tau}:
there we have $(\tau_u \omega)_e := \omega_{e - [s \cdot u]_{u \bmod 2}}$
for each $e \in E(\Z)$.
When $s = 0$ the quantity $s \cdot u$ is zero even as $u \to \infty$,
so the entire family of transformations $(\tau_u)_{u \in \Z^m}$
is actually finite rather than ergodic.
As such, a different argument is needed for $s = 0$.
The authors would like to thank Marek Biskup
for suggesting the following argument.

\begin{lemma}[Equivalence of quenched and annealed local surface tension]
  \label{l_equivalence_quenched_annealed_local_surface_tension}
For almost every $\omega$, it holds that
\begin{equation} \label{e_equivalence_quenched_annealed_local_surface_tension}
	\ent(s, \omega) = \ent_\an(s).
\end{equation}
Moreover, the quenched microscopic surface tension~$\ent_n(s, \omega)$
converges in~$L^1$ to~$\ent_\an(s)$.
\end{lemma}

\begin{proof}[Proof of
Lemma~\ref{l_equivalence_quenched_annealed_local_surface_tension}]
For $s \ne 0$, the desired
identity~\eqref{e_equivalence_quenched_annealed_local_surface_tension}
follows from the ergodic theorem,
as mentioned at the end of the proof of
Lemma~\ref{l_que_local_surface_tension}.
\\
For $s = 0$,
we will establish translation invariance of $\ent(s, \omega)$ directly.
First we replace the environmental shift $\tau_2$ by a shift in heights,
i.e.
\begin{equation} \llabel{e_shift}
	\Ent_{S_n} \bigl( M(S_n, h_{\partial S_n}^0) \bigr) \circ \tau_2
	= \Ent_{S_n} \bigl( M(S_n, h_{\partial S_n}^{0 \cdot x + 2}) \bigr) \,.
\end{equation}
This identity is justified simply by expanding definitions;
both sides are equal to
$-\tfrac{1}{\abs{S_n}} \log \sum_{h_{S_n}} \exp( \sum_e \omega_{h_{S_n}(e)+2})$,
where the first sum runs over $h_{S_n} \in M(S_n, h_{\partial S_n}^0)$
and the second runs over $e \in E(S_n)$.

Now, the square $S_n$ sits inside of $S_{n+2}$.
The boundary values $h_{\partial S_n}^{0 \cdot x + 2}$
and $h_{\partial S_{n+2}}^0$ satisfy
the Kirszbraun criterion~\eqref{e_kirszbraun};
in fact, each $h \in M(S_n, h_{\partial S_n}^{0 \cdot x + 2})$
admits a unique extension $\tilde h$ in $M(S_{n+2}, h_{\partial S_{n+2}}^0)$.
Since $\tilde h$ is an extension of $h$ to a domain with
$O(n^{m-1})$ more points and $O(n^{m-1})$ more edges,
the Hamiltonians satisfy
\[
	\abs[\big]{ H_{S_n}(h, \omega) - H_{S_{n+2}}(\tilde h, \omega) }
	\le c n^{m-1} C_\omega
\]
for some $c > 0$. Therefore
\begin{equation} \llabel{e_ext} \begin{aligned}
	\hskip3em&\hskip-3em
	\Ent_{S_n} \bigl(
		M(S_n, h_{\partial S_n}^{0 \cdot x + 2}), \omega \bigr) \\
	&\ge -\frac{1}{\abs{S_n}} \log
		\sum_{h \in M(S_n, h_{\partial S_n}^{0 \cdot x + 2})}
			\exp \bigl( H_{S_{n+2}}(\tilde h, \omega) \bigr)
	- \frac{c C_\omega}{n} \\
	&\ge \Ent_{S_{n+2}} \bigl(
		M(S_{n+2}, h_{\partial S_{n+2}}^0), \omega \bigr)
		- \frac{c C_\omega}{n} \,.
\end{aligned} \end{equation}

%
%

Now we combine~\leqref{e_shift} and~\leqref{e_ext} and send $n \to \infty$,
which yields
\[
	\ent(0,\omega) \circ \tau_2 \ge \ent(0,\omega) \,.
\]

By a similar argument with $\tau_2$ replaced by $\tau_{-2}$,
we conclude that $\ent(0,\tau_2\omega) = \ent(0,\omega)$,
i.e.\ $\ent(0,\omega)$ is invariant under $\tau_2$.
Since the distribution $\P$ of $\omega$ is ergodic with respect to $\tau_2$
(cf.\ Assumption~\ref{a_random_field}),
this implies that
$\ent(0,\omega) = \E[\ent(0,\omega)] = \ent_\an(0)$ almost surely.
\end{proof}

\subsection{Convexity and continuity}
\label{ss_lst_other}

The last results that we need about
the annealed local surface tension~$\ent_\an(s)$
are that is is convex and continuous as a function of the slope~$s$.
\\

Convexity allows us to apply standard analytic techniques to conclude that
the macroscopic entropy functional~$\Ent_{R,\an}(\cdot)$
is lower semi-continuous (see, for example, \cite[Section 2]{CKP01}).
By semi-continuity, there exists a (perhaps non-unique) minimizer of the
entropy functional, so the minimum in the variational principle
(Theorem~\ref{th_varprin}) is achieved.
\\

\begin{lemma} \label{p_convexity}
The function~$s \mapsto \ent_\an(s)$ is convex for~$s \in (-1,1)^m$.
\end{lemma}

\begin{remark}
The proof follows a standard argument based on buckled height functions;
see e.g.~\cite{KrMeTa19,She05} for the uniform case.
The energetic effect of the random potential contributes only on the boundary scale,
and so is negligible in the limit.
The proof could be considered an exercise for the reader;
we work out the details below.
\end{remark}

\begin{proof}[Proof of Lemma~\ref{p_convexity}]
We shall prove that
for any choice of fixed
coordinates~$s_1, \dotsc, s_{i-1}, s_{i+1}, \dotsc s_m \in [-1,1]^m$,
the single-variate functions~%
$s_i \mapsto \ent_\an((s_1, \dotsc, s_{i-1}, s_i, s_{i+1}, \dotsc, s_m)$
are convex.
It follows from elementary analysis that~$s \mapsto \ent_\an(s)$
is a convex function on the~$m$-dimensional domain~$[-1,1]^m$.
To simplify notation, we state the proof in the case~$m = 2$,
The proof generalizes to higher dimensions.
\\

So, choose~$u_0, u_1, u_2, v \in [-1, 1]$ such that
such that
\begin{equation}
	u_1 = \frac{1}{2} u_0 + \frac{1}{2} u_2.
\end{equation}

Our goal is to prove that
\begin{equation} \label{e_convexity_goal}
	\ent_\an((u_1,v))
	\le \frac{1}{2} \ent_\an((u_0, v))
	+ \frac{1}{2} \ent_\an((u_2, v)).
\end{equation}
We proceed as follows, in four steps.
\begin{itemize}
	\item First, consider a discrete hypercube~$S_{2n+1}$, which we recall
	is the hypercube~$\{-(2n+1), \dotsc, (2n+1)\}^2$ of side
	length~$2(2n+1)+1$ centered at the origin.
	We subdivide it into~$2m = 4$ smaller boxes.
	We choose height functions with slope~$(u_0, v)$
	or~$(u_2, v)$ on the smaller boxes,
	and we construct a bijection which maps
	from a choice of height functions on the four smaller boxes
	to a height function on the larger box.
	\item Second, we use the bijection to derive an inequality between the
	microscopic entropy on the four smaller boxes and an entropy-like
	quantity on the larger box.
	\item Third, we relate this ``entropy-like quantity'' to the annealed
	surface tension~$\ent_\an((u_1, v))$.
	\item Fourth, we relate the entropy on the smaller boxes to the
	right-hand side of~\eqref{e_convexity_goal},
	which concludes our proof.
\end{itemize}

\begin{figure}
	\subcaptionbox{
		The large box~$S_{2n+1}$ is divided into five subsets.
		The smaller boxes~$S_n^k$, for~$k=1,2,3,4$,
		are translated copies of the box~$S_n$ centered at the origin.
		The set~$S'$, indicated by dashed lines,
		is the intersection of~$S_{2n+1}$ and the~$x$- and~$y$-axes.
		\label{ff_convex_boxes}
	}[0.45\textwidth]{
		{ 

\newcommand*\ptrad{1.25pt}

\begin{tikzpicture}[
	point/.style = { shape = circle, },
	x=0.4in, y=0.4in,
	labelsty/.style = { font = \tiny },
	every label/.style = { labelsty },
]
	\useasboundingbox (-3,-3) rectangle (3,3);
	\coordinate (oo) at ( 0, 0);
	\coordinate (po) at ( 2, 0);
	\coordinate (mo) at (-2, 0);
	\coordinate (op) at ( 0, 2);
	\coordinate (om) at ( 0,-2);
	\foreach \v in {oo,po,mo,op,om}
		\fill [black] (\v) circle (\ptrad);

	\foreach \k/\xs/\xl/\ys/\yl in {
		1/-/m/+/p,
		2/-/m/-/m,
		3/+/p/-/m,
		4/+/p/+/p%
	} {
		\coordinate [label={center:$S_n^\k$}]
			(ctr\xl\yl) at (\xs1,\ys1);
		\coordinate (crn\xl\yl) at (\xs2,\ys2);
		\fill [black] (crn\xl\yl) circle (\ptrad);
		\coordinate (inn\xl\yl) at (\xs0.1,\ys0.1);
		\fill [black] (inn\xl\yl) circle (\ptrad);
	}

	\node [labelsty, below] at (om) {$S'$};

	\foreach \v [remember = \v as \u (initially crnpp)]
	in { po, crnpm, om, crnmm, mo, crnmp, op, crnpp}
		\draw [thin] (\u) -- (\v);
	\draw [thin,dashed] (po) -- (oo) -- (mo);
	\draw [thin,dashed] (op) -- (oo) -- (om);

	\foreach \xl in {p,m}
	\foreach \yl in {p,m}
		\draw [thick] (crn\xl\yl) -- (crn\xl\yl -| inn\xl\yl)
			-- (inn\xl\yl) -- (inn\xl\yl -| crn\xl\yl)
			-- (crn\xl\yl);
\end{tikzpicture} }
	}
	\hfill
	\subcaptionbox{
		To each of the four smaller boxes~$S_n^k$,
		we associate a slope~$s_k$.
		The two boxes on the left have~$s_k = (u_0,v)$
		and the two on the right have~$s_k = (u_2,v)$.
		If~$h_b$ is a \defn{buckled height function}
		(defined after~\eqref{e_convexity_tuple}),
		then~$h_b$ satisfies the indicated boundary conditions
		on the four smaller boxes.
		\label{ff_convex_slopes}
	}[0.45\textwidth]{
		{ 

\newcommand*\ptrad{1.25pt}

\begin{tikzpicture}[
	point/.style = { shape = circle, },
	x=0.4in, y=0.4in,
	labelsty/.style = { font = \tiny },
	every label/.style = { labelsty },
]
	\useasboundingbox (-3,-3) rectangle (3,3);
	\coordinate (oo) at ( 0, 0);
	\coordinate (po) at ( 2, 0);
	\coordinate (mo) at (-2, 0);
	\coordinate (op) at ( 0, 2);
	\coordinate (om) at ( 0,-2);
	\foreach \v in {oo,po,mo,op,om}
		\fill [black] (\v) circle (\ptrad);
	\foreach \xs/\xl in {-/m,+/p}
	\foreach \ys/\yl in {-/m,+/p} {
		\coordinate (ctr\xl\yl) at (\xs1,\ys1);
		\coordinate (crn\xl\yl) at (\xs2,\ys2);
		\fill [black] (crn\xl\yl) circle (\ptrad);
		\coordinate (inn\xl\yl) at (\xs0.1,\ys0.1);
		\fill [black] (inn\xl\yl) circle (\ptrad);
	}

	\foreach \v [remember = \v as \u (initially crnpp)]
	in { po, crnpm, om, crnmm, mo, crnmp, op, crnpp}
		\draw [thin] (\u) -- (\v);
	\draw [thin,dashed] (po) -- (oo) -- (mo);
	\draw [thin,dashed] (op) -- (oo) -- (om);

	\foreach \xl in {p,m}
	\foreach \yl in {p,m}
		\draw [thick] (crn\xl\yl) -- (crn\xl\yl -| inn\xl\yl)
			-- (inn\xl\yl) -- (inn\xl\yl -| crn\xl\yl)
			-- (crn\xl\yl);

	\node [labelsty] at (-1,1)
		{$\begin{gathered} s_1 \\ = (u_0,v) \end{gathered}$};

	\node [labelsty] at (-1,-1)
		{$\begin{gathered} s_2 \\ = (u_0,v) \end{gathered}$};

	\node [labelsty] at (1,-1)
		{$\begin{gathered} s_3 \\ = (u_2,v) \end{gathered}$};

	\node [labelsty] at (1,1)
		{$\begin{gathered} s_4 \\ = (u_2,v) \end{gathered}$};

	\draw [decoration={brace, amplitude=5}, decorate]
		($(innmp -| crnmp) - (6pt,0)$)
		-- node [labelsty, left=0.25] {$\begin{gathered}
			h_b|_{\partial S_n^1} \\
			\ = h_{\partial S_n^1}^{s_1}
		\end{gathered} $}
		($(crnmp) - (6pt,0)$);

	\draw [decoration={brace, amplitude=5}, decorate]
		($(crnmm) - (6pt,0)$)
		-- node [labelsty, left=0.25] {$\begin{gathered}
			h_b|_{\partial S_n^2} \\
			\ = h_{\partial S_n^2}^{s_2}
		\end{gathered} $}
		($(innmm -| crnmm) - (6pt,0)$);

	\draw [decoration={brace, amplitude=5}, decorate]
		($(innpm -| crnpm) + (6pt,0)$)
		-- node [labelsty, right=0.25] {$\begin{gathered}
			h_b|_{\partial S_n^3} \\
			\ = h_{\partial S_n^3}^{s_3}
		\end{gathered} $}
		($(crnpm) + (6pt,0)$);

	\draw [decoration={brace, amplitude=5}, decorate]
		($(crnpp) + (6pt,0)$)
		-- node [labelsty, right=0.25] {$\begin{gathered}
			h_b|_{\partial S_n^4} \\
			\ = h_{\partial S_n^4}^{s_4}
		\end{gathered} $}
		($(innpp -| crnpp) + (6pt,0)$);

\end{tikzpicture} }
	}
	\caption{Decomposition of~$S_{2n+1}$ into subsets,
	as used in the proof of Lemma~\ref{p_convexity}}.
	\label{f_convex_decomp}
\end{figure}

So, let us make precise how we decompose~$S_{2n+1}$.
We write
\begin{equation} \label{e_convex_nonzero_boxes}
	S_{2n+1} =  S_n^1 \cup S_n^2
	\cup S_n^3 \cup S_n^4
	\cup S'
\end{equation}
where
\begin{align}
	S_n^1 &:= \tau_{(-n,+n)} S_n, \\
	S_n^2 &:= \tau_{(-n,-n)} S_n, \\
	S_n^3 &:= \tau_{(+n,-n)} S_n, \\
	S_n^4 &:= \tau_{(+n,+n)} S_n, \text{ and} \\
	S' &:= \{ (x_1, x_2) \in S_{2n+1} \,|\, x_1=0 \text{ or } x_2=0 \}.
\end{align}
This decomposition is illustrated in Figure~\ref{ff_convex_boxes}.
\\

As an aside, it would be simpler if we could decompose~$S_{2n+1}$
into just the four boxes~$S_n^k$ without needing the extra set~$S'$.
But, both~$S_{2n+1}$ and~$S_n^k$ are centered boxes,
with an odd number of points along their edges
($4n+3$ and~$2n+1$ points, respectively),
so such a decomposition is arithmetically impossible.
Centered boxes are a requirement of the ergodic theorem
that we used to prove Lemma~\ref{l_que_local_surface_tension}
(the existence of the quenched local surface tension)
and Lemma~\ref{l_equivalence_quenched_annealed_local_surface_tension}
(the equivalence of the quenched and annealed local surface tension).
One could state the results without requiring odd-sized boxes
centered exactly at the origin,
but the statements become more complicated.
We choose instead to keep the odd-sized boxes,
and to keep the extra set~$S'$.
Because~$\abs{S'} = o(\abs{S_n})$,~$S'$ will be asymptotically negligible.
\\

To continue with the current proof,
we consider boundary height functions
of slope~$s_k$ on the small boxes~$S_n^k$, where
\begin{align}
	s_1 = s_2 &= (u_0, v), \text{ and} \\
	s_3 = s_4 &= (u_2,v).
\end{align}
This assignment of slopes to the small boxes is illustrated in
Figure~\ref{ff_convex_slopes}.
\\

Fix a~$4$-tuple of height functions
\begin{equation} \label{e_convexity_tuple}
	(h_n^k)_{k=1,2,3,4}
	\ \in \ %
	\prod_{k=1}^4 M(S_n^k, h_{\partial S_n^k}^{s_k}).
\end{equation}
We claim that there exists a height function~$h_b: S_{2n+1} \to \Z$
such that for each of the four boxes~$S_n^k$,
$h_b|_{S_n^k} = h_n^k$.
We call~$h_b$ a \defn{buckled height function},
since if~$h_b$ stays close to the linear height functions~$h_{S_n^k}^{s_k}$
over the entirety of the small boxes~$S_n^k$,
and if we view the graph of~$h_b$ in profile from along the~$y$-axis,
we see a buckled shape:
slope~$(u_0,v)$ along the left half,
which changes abruptly to slope~$(u_2,v)$ along the right half.
Figure~\ref{ff_convex_slopes} illustrates
the boundary conditions that are imposed on a buckled height function
on the boundaries~$\partial S_n^k$ of the small boxes.
\\

One can prove the existence of the height function~$h_b$
that extends the~$4$-tuple~$(h_n^k)$ to all of~$S_{2n+1}$
by using the Kirszbraun theorem.
However it is also easy to construct a concrete extension
using the canonical height functions.
Briefly, on either side of a point on the~$x$-axis,
the slopes~$s_k$ are equal.
For a point on the~$y$-axis,
the adjacent slopes differ only in the first coordinate.
That is not a problem because for a point to be on the~$y$-axis
means that the value in its first coordinate is~$0$.
\\

We write~$M_b$ for the set of all height functions~$h_b: S_{2n+1} \to \Z$
that can be realized by the above extension process.
Clearly, the set~$M_b$ is in bijection with the Cartesian product
of the four sets~$M(S_n^k, h_n^k)$.
This bijection completes the first step of our proof.
\\

In the second step of the proof, we derive the following approximation:
\begin{equation} \label{e_convexity_step2}
	-\frac{1}{\abs{S_{2n+1}}} \log Z_\omega(M_b)
	\le \frac{1}{4} \sum_{k=1}^4
		\Ent_{S_n} \bigl( M(S_n^k, h_{\partial S_n^k}^{s_k}),
			\omega \bigr)
	+ C_\omega \, \theta \bigl( \tfrac{1}{n} \bigr).
\end{equation}

The essential idea is that for a height function~$h: S_{2n+1} \to \Z$, the
Hamiltonian~$H_{S_{2n+1}}(h, \omega) = \sum_{e \in E(S_{2n+1})} \omega_{h(e)}$
splits as
\begin{align}
	H_{S_{2n+1}}(h, \omega)
	&= \sum_{k=1}^4 \sum_{e \in E(S_n^k)} \omega_{h(e)}
	+ \sum_{e \in E(S')} \omega_{h(e)}
	+ \sum_{\tilde e \in \tilde E} \omega_{h(e)} \\
	&\geq \sum_{k=1}^4 \sum_{e \in E(S_n^k)} \omega_{h(e)}
	- C_\omega \, O(n^{m-1}), \\
\end{align}
where~$\tilde E$ is the set of edges from~$E(S_{2n+1})$
that cross between two distinct parts of
the decomposition~$S_{2n+1} = S_n^1 \cup S_n^2 \cup S_n^3 \cup S_n^4 \cup S'$.
It follows that
\begin{align}
	-\frac{1}{\abs{S_{2n+1}}} \log Z_\omega(M_b)
	&= -\frac{1}{\abs{S_{2n+1}}} \log \sum_{h_b \in M_b}
		\exp H_{S_{2n+1}}(h_b, \omega) \\
	&\leq -\frac{1}{\abs{S_{2n+1}}} \log \vbig{4}[
		\prod_{k=1}^4 \vbig{4}(
			\sum_{h_n \in M \bigl(
			S_n^k, h_{\partial S_n^k}^{s_k} \bigr)}
				\exp H_{S_n}(h_n, \omega)
		\vbig{4}) \\
	&\hskip8em \exp \bigl(- C_\omega O(n^{m-1}) \bigr) \vbig{4}] \\
	&= -\frac{1}{4} \sum_{k=1}^4 \frac{1}{\abs{S_n}} \log
			\sum_{h_n \in M \bigl(
			S_n^k, h_{\partial S_n^k}^{s_k} \bigr)}
				\exp H_{S_n}(h_n, \omega) \\
	&\hskip8em + \frac{1}{\abs{S_{2n+1}}} C_\omega \, O(n^{m-1}) \\
	&= \frac{1}{4} \sum_{k=1}^4 \Ent_{S_n} \bigl( M(S_n^k, h_{S_n^k}^{s_k}),
		\omega \bigr)
	+ C_\omega \, \theta \bigl( \tfrac{1}{n} \bigr).
\end{align}
This proves~\eqref{e_convexity_step2}
and completes the second step of the proof.
\\

Two steps remain.
The third step is to relate the expression~%
$-\frac{1}{\abs{S_{2n+1}}} \log Z(M_b)$
(which we described as ``entropy-like'' earlier
when describing the steps of this proof)
to the annealed surface tension~$\ent_\an((u_1, v))$.
The fourth and final step is to verify that the microscopic entropy
$\Ent_{S_n} ( M(S_n^k, h_{\partial S_n^k}^{s_k}), \omega )$
converges to the
annealed surface tension~$\ent_\an(s_k)$ for~$k = 1, 2, 3, 4$.
This will suffice to prove the convexity inequality~\ref{e_convexity_goal}.
\\

To relate~$-\frac{1}{\abs{S_{2n+1}}} \log Z(M_b)$ and~$\ent_\an((u_1, v))$,
we first pass to the microscopic entropy
$\Ent_{S_{2n+1}}(M(S_{2n+1}, h_{\partial S_{2n+1}}^b), \omega)$.
The boundary height function~$h_{\partial S_{2n+1}}^b$
is given by~$h_{\partial S_{2n+1}}^b = h_b|_{\partial S_{2n+1}}$
for any~$h_b \in M_b$.
All the buckled height function~$h_b$ have the same boundary data
because of the boundary conditions on~$\partial S_n^k$,
plus the consistent (albeit arbitrary)  choice of extension to~$S'$.
Obviously~$M_b \subseteq M(S_{2n+1}, h_{\partial S_{2n+1}}^b)$,
and therefore by monotonicity,
\begin{equation}
	-\frac{1}{\abs{S_{2n+1}}} \log Z_\omega[M_b]
	\ge \Ent_{S_{2n+1}} \bigl( M(S_{2n+1}, h_{\partial S_{2n+1}}^b),
		\omega \bigr).
\end{equation}

To estimate~$\Ent_{S_{2n+1}}(M(S_{2n+1}, h_{\partial S_{2n+1}}^b), \omega)$,
let us consider any boundary point~$x = (x_1, x_2) \in \partial S_{2n+1}$.
If~$x_1 \le 0$, then the boundary height function~$h_{\partial S_{2n+1}}^b(x)$
is equal to~$u_0 \cdot x$ up to a rounding error of at most~$1$, so
\begin{equation} \label{e_convex_buckled_approx_lin}
	\abs{ h_{\partial S_{2n+1}}^b(x) - (u_1, v) \cdot x }
	\le \abs{ u_0 - u_1 } \abs{x_1}
	\le \varepsilon n,
\end{equation}
where~$\varepsilon = \abs{ u_0 - u_1 } = \abs{ u_0 - u_2 }$.
If instead~$x_1 \ge 0$,
then~$h_{\partial S_{2n+1}}^b(x) = u_ \cdot x$ up to rounding error,
so still~\eqref{e_convex_buckled_approx_lin} holds.
Therefore, the boundary data~$h_b$ is approximately linear,
i.e.\ $h_b \in M(\partial S_{2n+1}, h_{\partial S_{2n+1}}^u, \varepsilon)$.
By Theorem~\ref{p_robustness}, there exists~$A = A((u_1, v), \varepsilon) > 0$
such that
\begin{equation} \label{e_convex_approx_lin_surf_tens}
	\Ent_{S_{2n+1}} \bigl( M(S_{2n+1}, h_{\partial S_{2n+1}}^b),
		\omega \bigr)
	\ge \ent_{A (2n+1)}((u_1, v), \omega).
\end{equation}

We have proved the following inequality, which concludes the third step:
\begin{equation} \label{e_convex_big_micro}
	\ent_{A (2n+1)}((u_1, v), \omega)
	\le -\frac{1}{\abs{S_{2n+1}}} \log Z_\omega(M_b).
\end{equation}

In the last step, we consider the
quenched microscopic entropy on the four sub-boxes, i.e.\
$\Ent_{S_n^k}( M(S_n^k, h_n^k), \omega)$.
Recall that each box~$S_n^k$ is a translation~$\tau_{(\pm n, \pm n)} S_n$
of the box~$S_n$ centered at the origin.
We transfer the translation over to the height function and environment.
Let~$\tau_1: E(\Z) \to E(\Z)$ denote the shift by~$s_k \cdot (-n, +n)$,
so that~$\tau_1 \circ h_n^1 = h_{S_n}^{s_1}$
and~$
	\Ent_{S_n^1} ( M(S_n^1, h_{\partial S_n^1}^{s_1}), \omega )
	= \ent_n(s_1, \tau_1 \omega)
$.
Likewise, define~$\tau_2, \tau_3, \tau_4: E(\Z) \to E(\Z)$
so that for each~$k = 1,2,3,4$, it holds that
\begin{equation}
	\Ent_{S_n^k} \bigl( M(S_n^k, h_{\partial S_n^k}^{s_k}), \omega \bigr)
	= \ent_n(s_k, \tau_1 \omega).
\end{equation}
Combining this identity with~\eqref{e_convexity_step2}
and~\eqref{e_convex_big_micro},
we deduce a quenched microscopic inequality
\begin{equation}
	\ent_{A (2n+1)}((u_1, v), \omega)
	\le \frac{1}{4} \sum_{k=1}^4 \ent_n(s_k, \tau_k \omega)
	+ C_\omega \, \theta \bigl( \tfrac{1}{n} \bigr).
\end{equation}

Taking expectations, the annealed microscopic inequality is
\begin{equation}
	\ent_{A (2n+1), {\tt an}}((u_1, v))
	\le \frac{1}{4} \sum_{k=1}^4 \ent_{n, {\tt an}}(s_k)
	+ \mathbb{E}(C_\omega) \, \theta \bigl( \tfrac{1}{n} \bigr).
\end{equation}

Inequality~\eqref{e_convexity_goal},
which states that the annealed local surface tension is convex,
follows immediately by sending~$n \to \infty$.
\end{proof}

\section{Profile theorem} \label{s_profile}

Before proving the profile theorem, Theorem~\ref{th_profile},
in its full generality,
it is useful to prove a special case of the theorem
with the extra assumptions that the asymptotic height function
is piecewise affine, on a domain which is of a collection of simplices.
In this special case it is not difficult to relate the microscopic entropy
$\Ent_{R_n}( B(R_n, h_R, \delta),\omega)$
to the quenched microscopic surface tension $\ent_n(s, \omega)$,
and then to derive the desired conclusion~\eqref{e_profile_thm_ineq}.
The special case is stated in Lemma~\ref{lem_simpl_profthm} below,
after some necessary notation is introduced in Definitions~\ref{d_simpl_dom}
and~\ref{d_simpl_hf}.

\begin{definition}[Simplices of scale $\ell$;
cf.\ {\cite[Definition 27]{KrMeTa19}} and {\cite[Section 5.2.1]{She05}}]
\label{d_simpl_dom}
Let $\Sym(m)$ denote the group of permutations on $\{1, \dotsc, m\}$,
and for $w = (w_1, \dotsc, w_m) \in \R^m$,
let $\floor{w}$ denote the integer point
$\floor{w} := (\floor{w_1}, \dotsc, \floor{w_m})$.
Let $v \in \Z^m$, let $\sigma \in \Sym(m)$, and let $\ell > 0$.
Define $C(v,\sigma)$ to be the closure of the set
\[
	\bigl\{ w \in \R^m
		\,\big|\, \floor{w} = v \textnormal{ and }
		w_{\sigma(1)} - \floor{w_{\sigma(1)}}
		> \dotsb
		> w_{\sigma(m)} - \floor{w_{\sigma(m)}} \bigr\} \,,
\]
and define the simplex of scale $\ell$ to the scaled set
\[
	\ell C(v,\sigma) := \{ \ell w \,|\, w \in C(v,\sigma) \} \,.
\]
\end{definition}

\begin{definition}[Piecewise affine asymptotic height functions]
\label{d_simpl_hf}
Let $\Delta_1, \dotsc, \Delta_k$ be simplices of scale $\ell$
and let $K = \Delta_1 \cup \dotsb \cup \Delta_k$ be their union.
We say that an asymptotic height function $h_K \in M(K)$ is piecewise affine
if each restriction $h_K|_{\Delta_i}$ is an affine function,
i.e.\ if there exist $s_i \in [-1,1]^m$ and $b_i \in \R$ such that
$h_K|_{\Delta_i}(x) = s_i \cdot x + b_i$ for all $x \in \Delta_i$.
We write
\begin{equation} \label{e_ht_func_sets_aff} \begin{aligned}
	M_\aff(K)
	&= \bigl\{ h_K \in M(K) \,\bigm|\, \text{$h_K$ is piecewise affine}
	\bigr\} \\
	M_\aff(K, h_{\partial K})
	&= M_\aff(K) \cap M(K, h_{\partial K}) \,.
\end{aligned} \end{equation}
\end{definition}

\begin{lemma}[Profile theorem, simplicial case]
\label{lem_simpl_profthm}
Let $\Delta_1, \dotsc, \Delta_k$ be simplices of scale $\ell$
and let $K = \Delta_1 \cup \dotsb \cup \Delta_k$ be their union.
\\

For any $h_K \in M_\aff(K, h_{\partial K})$ and any $\eta > 0$,
there exists $\varepsilon = \varepsilon_0(h_K, \eta)$ such that,
for any $\varepsilon \in (0, \varepsilon_0]$ and any $p_\tmax \in (0, 1)$,
there exists $n_0 = n_0(h_K, \eta, \varepsilon, p_\tmax)$ such that
for all $n \geq n_0$,
\begin{equation} \label{e_simpl_profile} \begin{aligned}
	\hskip3em&\hskip-3em
	\P \Bigl( \, \abs[\Big]{
		\Ent_{K_n} \bigl( B(K_n, h_K, \varepsilon \ell), \omega \bigr)
		- \Ent_{K,\an}( h_K ) } \\
	&> \eta
	+ C_\omega \theta_{h_K} (\varepsilon)
	+ C_\omega \theta_{h_K,\varepsilon} \bigl( \tfrac{1}{n} \bigr)
	\Bigr) < p_\tmax \,.
\end{aligned} \end{equation}
\end{lemma}

\begin{proof}
We will prove two bounds on the quenched microscopic entropy
$\Ent_{K_n}(B(K_n,h_{K_n},\varepsilon \ell),\omega)$:
an upper bound
\begin{equation} \llabel{e_upper} \begin{aligned}
	\hskip3em&\hskip-3em
	\P \Bigl( \Ent_{K_n} \bigl( B(K_n, h_K, \varepsilon \ell), \omega \bigr)
	> \Ent_{K,\an}(h_K) \\
	&+ \eta
	+ C_\omega \theta_{h_K}(\varepsilon)
	+ C_\omega \theta_{h_K,\varepsilon} \bigl( \tfrac{1}{n} \bigr) \Bigr)
	\le \theta_{h_K,\eta,\varepsilon} \bigl( \tfrac{1}{n} \bigr)
\end{aligned} \end{equation}
and a lower bound
\begin{equation} \llabel{e_lower} \begin{aligned}
	\hskip3em&\hskip-3em
	\P \Bigl( \Ent_{K_n} \bigl( B(K_n, h_K, \varepsilon \ell), \omega \bigr)
	< \Ent_{K,\an}(h_K) \\
	&- \eta
	- C_\omega \theta_{h_K}(\varepsilon)
	- C_\omega \theta_{h_K,\varepsilon} \bigl( \tfrac{1}{n} \bigr) \Bigr)
	\le \theta_{h_K,\eta,\varepsilon} \bigl( \tfrac{1}{n} \bigr) \,.
\end{aligned} \end{equation}

Assuming that both~\leqref{e_upper} and~\leqref{e_lower} hold,
the conclusion~\eqref{e_simpl_profile} follows immediately
by taking $n_0$ large enough based on
the two $\theta_{h_k,\eta,\varepsilon}(\tfrac{1}{n})$ terms
and applying the union bound on probabilities.
So first let us verify the upper bound~\leqref{e_upper},
and later we will verify the lower bound~\leqref{e_lower}.
For~\leqref{e_upper} we undercount
the set of height functions $B(K_n, h_{K_n}, \varepsilon \ell)$.
We choose a fine mesh of hypercubes $Q_{i,n}$ that approximate $K_n$
and consider only those height functions that agree
with the canonical boundary height functions
$h_{\partial Q_{i,n}}^{s_i \cdot x + b_i}$ on $\partial Q_{i,n}$,
where $s_i \in [-1,1]^m$ and $b_i \in \R$
are chosen such that $s_i \cdot x + b_i = h_K|_{Q_i}$.
The mesh size is small enough that every such height function
is in $B(K_n, h_{K_n}, \varepsilon \ell)$.
\\

To be precise, let $q = \tfrac{1}{4} \varepsilon \ell$ be the mesh size.
Let $Q_1, \dotsc, Q_k \subset \R^m$ enumerate
the set of hypercubes in $\R^m$ that have side length $q$,
have vertices in $q \Z^m$,
and lie entirely in one of the simplices $\Delta_j$.
That last property ensures that
there exist $s_i \in [-1,1]^m$ and $b_i \in \R$ such that
\[
	h_K(x) = s_i \cdot x + b_i
	\qquad \text{for all $x \in Q_i$}
	\,.
\]
For $n \in \N$, let $Q_{i,n} := \{z \in \Z^m \,|\, \tfrac{1}{n}z \in Q_i\}$.
Then as desired, for any choice of height functions
\[
	\bigl( h_{Q_{i,n}} \bigr)_{i=1}^k
	\in \prod_{i=1}^k M \bigl(
		Q_{i,n}, h_{\partial Q_{i,n}}^{s_i \cdot x + b_i} \bigr) \,,
\]
there exists at least one extension $h_{K_n} \in M(K_n)$ to the whole of $K_n$
(i.e.\ $h_{K_n}|_{Q_{i,n}} = h_{Q_{i,n}}$ for each $i=1,\dotsc,k$),
and any such extension lies in $B(K_n, h_K, \varepsilon \ell)$ by choice of $q$.
Therefore,
\begin{equation} \begin{aligned} \label{e_simpl_profile_upper_decomp}
	\Ent_{K_n} \bigl( B(K_n, h_K, \varepsilon \ell), \omega \bigr)
	&\le \frac{1}{k} \sum_{i=1}^k \Ent_{Q_{i,n}} \bigl(
		M(Q_{i,n}, h_{\partial Q_{i,n}}^{s_i \cdot x + b_i}),
	\omega \bigr) \\
	&\qquad + C_\omega \theta_m(\varepsilon)
	+ C_\omega \theta_{m,\varepsilon,\ell} \bigl( \tfrac{1}{n} \bigr) \,,
\end{aligned} \end{equation}
where the $\theta$ error terms come from the contribution of the
set $K_n \setminus \bigcup_{i=1}^k Q_{i,n}$.
For each $i=1, \dotsc, k$,
let us abuse notation and write ``$qn$'' to denote
the side length of the hypercube $Q_{i,n}$.
(In fact, the actual product $q \cdot n$ is generally not an integer,
but the quantity we call $qn$ satisfies $\abs{qn - q \cdot n} < 1$.)
Consider $Q_{i,n}$ as a translate
$Q_{i,n} = v_i + S_{qn}$ for $v_i \in \Z^m$.
Then the boundary values $h_{\partial Q_{i,n}}^{s_i \cdot x + b_i}$
are close to the translated values of $h_{\partial S_{qn}}^{s_i}$;
in particular, for $z \in \partial S_{qn}$,
\begin{equation} \label{e_simpl_profile_upper_bdy}
	\abs[\Big]{ h_{\partial Q_{i,n}}^{s_i \cdot x + b_i}(v_i + z) -
		\Bigl( h_{\partial S_{qn}}^{s_i}(z)
			+ \floor{s_i \cdot v_i + n b_i} \Bigr)
	} \le 4 \,.
\end{equation}
(A non-zero error occurs when $s_i$ is irrational,
or more generally when $qn s_i$ is not integral or has the wrong parity.)
By Corollary~\ref{cor_robustbdy} it follows that
\begin{equation} \label{e_simpl_profile_upper_trans} \begin{aligned}
	\Ent_{Q_{i,n}} \bigl(
		M(Q_{i,n}, h_{\partial Q_{i,n}}^{s_i \cdot x + b_i}), \omega
	\bigr)
	= \ent_{qn}(s_i, \tau_{\floor{s_i \cdot v_i + n b_i}} \omega)
	+ C_\omega \theta_m \bigl( \tfrac{1}{n} \bigr) \,.
\end{aligned} \end{equation}

Combining~\eqref{e_simpl_profile_upper_decomp}
and~\eqref{e_simpl_profile_upper_trans}
and abbreviating $\tau_{i,n} := \tau_{\floor{s_i \cdot v_i + n b_i}}$ yields
\begin{equation} \label{e_simpl_profile_upper_decomptrans} \begin{aligned}
	\hskip3em&\hskip-3em
	\Ent_{K_n} \bigl( B(K_n, h_K, \varepsilon \ell), \omega \bigr) \\
	&\le \frac{1}{k} \sum_{i=1}^k
		\ent_{qn}(s_i, \tau_{i,n} \omega)
	+ C_\omega \theta_m(\varepsilon)
	+ C_\omega \theta_{m,\varepsilon,\ell} \bigl( \tfrac{1}{n} \bigr) \,.
\end{aligned} \end{equation}

We note that the sequences
$\{ \ent_{qn}(s_i, \tau_{i,n} \omega) \}_{n \in \N}$
may not necessarily converge to $\ent_\an(s)$ as $n \to \infty$,
despite the almost-sure convergence result of
Lemma~\ref{l_equivalence_quenched_annealed_local_surface_tension},
due to the potential shifts $\tau_{i,n}$.
However, since each $\ent_{qn}(s_i, \cdot) \to \ent_\an(s_i)$ in $L^1$,
we can apply the Markov bound:
\begin{equation} \llabel{e_markov} \begin{aligned}
	\hskip3em&\hskip-3em
	\P \biggl(
		\abs[\bigg]{ \,
			\frac{1}{k} \sum_{i=1}^k
				\ent_{qn}(s_i, \tau_{i,n} \omega)
			- \frac{1}{k} \sum_{i=1}^k \ent_\an(s_i)
		}
		> \eta \biggr) \\[12pt]
	&\le \frac{1}{k} \sum_{i=1}^k \frac{1}{\eta} \, \norm[\big]{
		\ent_{qn}(s_i, \cdot) - \ent_\an(s_i) }_{L^1} \\
	&= \theta_{h_K,\eta,\varepsilon,\ell} \bigl(\tfrac{1}{n} \bigr) \,.
\end{aligned} \end{equation}

The last step in verifying~\leqref{e_upper}
is to compare $\Ent_{K,\an}(h_K)$ to a sum involving $\ent_\an(s_i)$.
This is straightforward:
because $h_K$ is affine on each hypercube $Q_i$,
the integrand $x \mapsto \ent_\an(\nabla h_K(x))$ in the macroscopic entropy
is constant on each $Q_i$,
so
\begin{equation} \llabel{e_int} \begin{aligned}
	\Ent_{K,\an}(h_K)
	&\overset{\text{def.}}=
		\frac{1}{\abs{K}} \int_K \ent_\an(\nabla h_K(x)) \, dx \\
	&= \frac{1}{k} \sum_{i=1}^k
		\frac{1}{\abs{Q_i}} \int_{Q_i} \ent_\an \bigl(
			\nabla h_K|_{Q_i} \bigr)
	+ \theta_K(\varepsilon) \\
	&= \frac{1}{k} \sum_{i=1}^k \ent_\an(s_i) + \theta_K(\varepsilon) \,.
\end{aligned} \end{equation}
The only error is from the contribution of the region
$K \setminus \bigcup_{i=1}^k Q_i$.
Combining inequalities~\eqref{e_simpl_profile_upper_decomptrans},
\leqref{e_markov}, and~\leqref{e_int}
proves the desired upper bound~\leqref{e_upper},
i.e.
\begin{equation} \begin{aligned}
	\hskip3em&\hskip-3em
	\P \Bigl(
		\Ent_{K_n} \bigl( B(K_n, h_K, \varepsilon \ell), \omega \bigr)
		> \Ent_{K,\an}(h_K) \\
	&
		+ \eta
		+ \theta_{h_K}(\varepsilon)
		+ \theta_{h_K,\varepsilon} \bigl( \tfrac{1}{n} \bigr) \Bigr)
	\le \theta_{h_K,\eta,\varepsilon} \bigl( \tfrac{1}{n} \bigr) \,.
\end{aligned} \end{equation}

\medskip

Now we turn to the lower bound~\leqref{e_lower}.
Similar to before,
let $q = \varepsilon^{1/2} \ell$
and let $Q_1, \dotsc, Q_{k}$ enumerate
the hypercubes that have side length $q$,
have vertices in $q \Z^m$,
and lie entirely inside of one of the simplices $\Delta_j$.
Note that the side length $q$ is different now
compared to above when we were justifying the upper bound~\leqref{e_upper},
and hence $Q_1, \dotsc, Q_{k}$ denotes a different set of hypercubes.

To prove~\leqref{e_lower} we overcount height functions,
using the same idea as in the companion article~\cite{KrMeTa19}.
In summary, define a subset of ``exceptional'' points
$E_n \subset K_n$ as follows: let
\begin{equation} \begin{aligned}
	G_n = \bigcup_{i=1}^{k} \partial Q_{i,n} \,,
	\quad
	U_n = K_n \setminus \bigcup_{i=1}^{k} Q_{i,n} \,,
	\quad \text{and} \quad
	E_n = G_n \cup U_n
	\,.
\end{aligned} \end{equation}

Informally, $G_n$ is the ``grid'' formed by the boundaries of the hypercubes
and $U_n$ is the ``uncovered'' region, i.e.\ the part of $K_n$
that is not covered by the hypercubes.
We group height functions $h_{K_n} \in B(K_n, h_K, \varepsilon \ell)$
based on their values on the set $E_n$.
For each fixed assignment of heights $h_{K_n}|_{E_n} \in M(E_n)$,
the entropy of the set of extensions to the hypercubes
$\bigcup_1^k Q_n \approx K_n \setminus E_n$
is asymptotically equal to the macroscopic entropy $\Ent_{K,\an}(h_K)$.
The set $E_n$ is not too large, so even after counting
all admissible assignments $h_{K_n}|_{E_n}$,
the resulting asymptotics match~\leqref{e_lower}.
\\

To make the above argument rigorous,
let $\Adm(E_n)$ denote the set of admissible height functions on $E_n$,
i.e.\ those height functions $h_{E_n} \in M(E_n)$ that admit an extension
to a height function in $B(K_n, h_K, \varepsilon \ell)$.
There is an obvious injection from $B(K_n, h_K, \varepsilon \ell)$ into
\begin{equation} \label{e_simpl_profile_lower_set}
	\biguplus_{h_{E_n} \in \Adm(E_n)} \prod_{i=1}^{k} M \bigl(
		Q_{i,n}, h_{E_n}|_{\partial Q_{i,n}} \bigr) \,,
\end{equation}
where ``$\biguplus$'' denotes the disjoint union
(so for distinct height functions $h_{E_n}$ and $h_{E_n}$ in $\Adm(E_n)$,
the product sets $
	\prod_1^{k} M \bigl( Q_{i,n},
		h_{E_n}|_{\partial Q_{i,n}} \bigr)
$ and $
	\prod_1^{k} M \bigl( Q_{i,n},
		h_{E_n}|_{\partial Q_{i,n}} \bigr)
$
are considered disjoint inside the set
from~\eqref{e_simpl_profile_lower_set}).
It follows that
\begin{equation} \begin{aligned}
	\hskip3em&\hskip-3em
	Z_\omega \bigl( B(K_n, h_K, \varepsilon \ell), \omega \bigr) \\
	&\le \sum_{h_{E_n} \in \Adm(E_n)}
		Z_\omega \Biggl( \, \prod_{i=1}^{k} M \bigl(
			Q_{i,n}, h_{E_n}|_{Q_{i,n}} \bigr) \Biggr)
	\\
	&\le \abs[\big]{ \Adm(E_n) }
	\, \max_{h_{E_n} \in \Adm(E_n)} Z_\omega \Biggl( \,
		\prod_{i=1}^{k} M \bigl(
			Q_{i,n}, h_{E_n}|_{Q_{i,n}} \bigr)
	\Biggr) \,.
\end{aligned} \end{equation}

Therefore
\begin{equation} \llabel{e_lower_entadm} \begin{aligned}
	\hskip3em&\hskip-3em
	\Ent_{K_n} \bigl( B(K_n, h_K, \varepsilon \ell), \omega \bigr) \\
	&\ge \min_{h_{E_n} \in \Adm(E_n)} \sum_{i=1}^{k}
		\frac{\abs{Q_{i,n}}}{\abs{K_n}}
		\Ent_{Q_{i,n}} \bigl(
			M(Q_{i,n}, h_{E_n}|_{Q_{i,n}}), \omega
		\bigr) \\
	& \qquad - \frac{\log \, \abs{\Adm(E_n)}}{\abs{K_n}}
	\,.
\end{aligned} \end{equation}

Clearly
$
	\frac{\abs{Q_{i,n}}}{\abs{K_n}}
	= \frac{1}{k} + \theta_m(\varepsilon)
		+ \theta_{m,\varepsilon,\ell} \bigl( \tfrac{1}{n} \bigr)
$.
\\

To control $\abs{\Adm(E_n)}$, we argue as follows.
First, $\tfrac{\abs{G_n}}{\abs{K_n}} = \theta_m(q) = \theta_m(\varepsilon)$
and $\tfrac{\abs{U_n}}{\abs{K_n}} = \theta_m(\varepsilon)$.
Second, for an arbitrary base point $z_0 \in E_n$,
there are at most $2 \varepsilon \ell n + 1$ admissible values
for $h_{E_n}(z_0)$ if $h_{E_n} \in \Adm(E_n)$,
since $h_{E_n}$ must extend to a height function in the ball
$B(K_n, h_K, \varepsilon \ell)$.
Third, the set $E_n$ is connected,
so for each of the admissible values of $h_{E_n}(z_0)$,
there are at most $2^{\abs{E_n}}$ height functions in $\Adm(E_n)$
taking that value at $z_0$.
Putting these observations together, we conclude that
$
	\tfrac{1}{\abs{K_n}} \log \, \abs{\Adm(E_n)}
	= \theta_m(\varepsilon)
		+ \theta_{m,\varepsilon,\ell} \bigl( \tfrac{1}{n} \bigr)
$.
\\

Applying these asymptotic results in~\leqref{e_lower_entadm} yields
\begin{equation} \llabel{e_lower_micro_min} \begin{aligned}
	\hskip3em&\hskip-3em
	\Ent_{K_n} \bigl( B(K_n, h_K, \varepsilon \ell), \omega \bigr) \\
	&\ge \min_{h_{E_n} \in \Adm(E_n)} \frac{1}{k} \sum_{i=1}^{k}
		\Ent_{Q_{i,n}} \bigl(
			M(Q_{i,n}, h_{E_n}|_{Q_{i,n}}),
				\omega \bigr) \\
	& \qquad - \theta_m(\varepsilon)
	- \theta_{m,\varepsilon,\ell} \bigl( \tfrac{1}{n} \bigr)
	\,.
\end{aligned} \end{equation}

Whenever $h_{E_n} \in \Adm(E_n)$,
\begin{equation}
	\max_{z \in E_n} \, \abs[\big]{
		h_K(\tfrac{1}{n} z) - \tfrac{1}{n}h_{E_n}(z)}
	< \varepsilon \ell \,,
\end{equation}
so for each $i=1, \dotsc, k$, by analogy to~\eqref{e_simpl_profile_upper_bdy},
\begin{equation}
	\max_{z \in \partial S_{q n}} \, \abs[\Big]{
		\bigl( h_{E_n}(v_i + z)
			- \floor{s_i \cdot v_i + q n b_i} \bigr)
		- h_{\partial S_{q n}}^{s_i}(z)}
	\le \varepsilon \ell n \,.
\end{equation}

We apply Theorem~\ref{p_robustness} to the height function
\[
	\bigl( z \mapsto h_{E_n}(v_i + z)
		- \floor{s_i \cdot v_i + qn b_i} \bigr)
	\in M(S_{qn})
\]
to conclude that
\begin{equation} \llabel{e_lower_micro_cts} \begin{aligned}
	\hskip3em&\hskip-3em
	\Ent_{Q_{i,n}} \bigl(
		M(Q_{i,n}, h_{E_n}|_{\partial Q_{i,n}}),
			\omega \bigr) \\
	&\ge \ent_{A qn} \bigl(
		s_i, \tau_{\floor{s_i \cdot v_i + qn b_i}} \omega \bigr)
	- C_\omega \theta(\varepsilon) \,.
\end{aligned} \end{equation}

The two almost-sure inequalities~\leqref{e_lower_micro_min}
and~\leqref{e_lower_micro_cts},
the probability estimate~\leqref{e_markov},
and the macroscopic bound~\leqref{e_int}
together imply the desired lower bound~\leqref{e_lower},
which completes the proof of Lemma~\ref{lem_simpl_profthm}.
\end{proof}

The remainder of the proof of the profile theorem (Theorem~\ref{th_profile})
for general asymptotic height functions
follows closely the proof in Section~6 of the companion article~\cite{KrMeTa19}.
Below we state an approximation result (Theorem~\ref{th_rade}),
which concludes that any asymptotic height function $h_R$
admits a ``good'' approximation $h_K$
satisfying the hypotheses of Lemma~\ref{lem_simpl_profthm} above.
Following that result are three robustness lemmas
(Lemmas~\ref{lem_approx_macro_ent}, Lemma~\ref{lem_approx_micro_ent_ahf},
and Lemma~\ref{lem_approx_micro_ent_dom}).
With these tools it is straightforward
to reduce the general case of Theorem~\ref{th_profile}
to the special case of Lemma~\ref{lem_simpl_profthm}.

The approximation result (Theorem~\ref{th_rade}
is unchanged from the companion article,
which should be expected because the random potential in the current model
does not affect the class of limit objects that our model admits,
i.e.\ domains satisfying Assumption~\ref{a_domain}
and asymptotic height functions.
It is similar to \cite[Lemma 2.2]{CKP01} or \cite[Theorem 1]{Sch14}.

\begin{theorem}[Simplicial Rademacher theorem] \label{th_rade}
Let $R \subseteq \R^m$ be a region satisfying Assumption~\ref{a_domain},
and let $h_R \in M(R, h_{\partial R})$ be an asymptotic height function on $R$.
For any $\varepsilon > 0$
and any $\ell > 0$ sufficiently small (depending on $\varepsilon$),
we may choose a simplex domain
$K = \Delta_1 \cup \dotsb \cup \Delta_k  \subseteq R$ of scale $\ell$
(see Definition~\ref{d_simpl_dom})
and a piecewise affine asymptotic height function $h_K: K \to \R$
(that is, an asymptotic height function
such that each restriction $h_K|_{\Delta_i}: \Delta_i \to \R$ is affine)
that satisfy the following properties:
\begin{enumerate}
\item \label{lem_approx_tri_volume}
	$\abs{ R \setminus K } < \varepsilon$ 
	and $d_H(K, R) < \varepsilon)$,
	where we recall that for subsets of $\R^m$,
	$\abs{\cdot}$ denotes the Lebesgue measure
	and $d_H(\cdot, \cdot)$ denotes Hausdorff metric;
\item \label{lem_approx_tri_approx}
	$\max_{x \in K} \abs{ h_K(x) - h_R(x) } < \tfrac{1}{2}\varepsilon \ell$;
	and
\item \label{lem_approx_tri_grad}
	on at least a $(1 - \varepsilon)$ fraction of the points in $K$
	(by Lebesgue measure), the gradients $\nabla h_K(x)$ and $\nabla h_R(x)$
	agree to within $\varepsilon$, i.e.\ %
	$
		\frac{1}{\abs{K}} \abs[\big]{ \bigl\{
			x \in K
			\,\bigm\vert\, \abs{ \nabla h_K(x) - \nabla h_R(x)}_2
			\ge \varepsilon
		\bigr\} }
		< \varepsilon
	$.
\end{enumerate}
\end{theorem}

A proof of this lemma is given in the companion article~\cite{KrMeTa19}.

Now we turn to the robustness lemmas,
which will be used when applying Theorem~\ref{th_rade}
to approximate $h_R$ by another asymptotic height function.
The three lemmas below are almost direct analogues of
Lemmas~35, 36, and~37 from~\cite{KrMeTa19} respectively.

\begin{lemma}[Robustness of macroscopic entropy under approximations]
\label{lem_approx_macro_ent}
Let $\varepsilon > 0$,
and let $\tilde R \subseteq R \subset \R^m$ be sets
meeting the assumptions from Assumption~\ref{a_domain}
with $\abs{ R \setminus \tilde R } < \varepsilon$.
Let $h_{\tilde R} \in M({\tilde R})$ and $h_R \in M(R)$ be such that
\begin{equation} \llabel{e_agree}
	\abs[\Big]{ \Bigl\{ x \in \tilde R \Bigm\vert
		\abs[\big]{\nabla h_{\tilde R}(x) - \nabla h_R(x)}_2
		\ge \varepsilon
	\Bigr\} } < \varepsilon \,.
\end{equation}
Then,
\[
	\Ent_{R,\an}(h_R)
	= \Ent_{\tilde R, \an}(h_{\tilde R})
	+ \theta_m(\varepsilon) \,.
\]
\end{lemma}

\begin{proof}
Recall from Definition~\ref{d_ann_macro_ent} that
\[
	\Ent_{R,\an}(h_R)
	:= \frac{1}{\abs{R}} \int_R
		\ent_\an \bigl( \nabla h_R(x) \bigr) \, dx \,,
\]
and likewise for $\Ent_{\tilde R,\an}(h_{\tilde R})$.
The conclusion follows from three observations:
first that the domains of integration are bounded sets
with small symmetric difference,
second that the function~$s \mapsto \ent_\an(s)$ is continuous,
and third that the functions $\nabla h_R$ and $\nabla h_{\tilde R}$
almost agree (as per~\leqref{e_agree})
on most of the intersection of their domains (by measure).
\end{proof}

\medskip

\begin{lemma}[Robustness of microscopic entropy under change in profile]
\label{lem_approx_micro_ent_ahf}
Let $\varepsilon > 0$ and $n \in \N$.
Let $R \subset \R^m$ satisfy Assumption~\ref{a_domain},
and let $R_n \subset \Z^m$ satisfy $\tfrac{1}{n} R_n \subset R$.
Let $h_R, \tilde h_R \in M(R)$ be two asymptotic height functions
such that $\sup_{x \in R} \abs{ h_R(x) - \tilde h_R(x) } \le \varepsilon$.
Then,
\[
	\Ent_{R_n} \bigl( B(R_n, h_R, 2 \varepsilon), \omega \bigr)
	\le \Ent_{R_n} \bigl( B(R_n, \tilde h_R, \varepsilon), \omega \bigr) \,.
\]
\end{lemma}

\begin{proof}
For any fixed $\omega$,
the functional $\Ent_{R_n}(\cdot, \omega): M(R_n) \to \R$ is monotonic,
and it follows from Definition~\ref{d_ht_func_sets} that
\[
	B(R_n, \tilde h_R, \varepsilon) \subseteq B(R_n, h_R, 2\varepsilon) \,.
\]
\end{proof}

\medskip

\begin{lemma}[Robustness of microscopic entropy under domain approximations]
\label{lem_approx_micro_ent_dom}
Let $c \in (0,1]$, $\varepsilon \in (0, 1]$, and $n \in \N$.
Let $\tilde R \subset R \subset \R^m$ and $\tilde R_n \subset R_n \subset \Z^m$
satisfy these assumptions:
\begin{align}
	\tfrac{1}{n} R_n &\subset R \,, &
	\tfrac{1}{n} \tilde R_n &\subset \tilde R \,, \\
	d_H(\tfrac{1}{n} R_n, R) &= \theta_R(\varepsilon) \,, &
	d_H(\tfrac{1}{n} \tilde R_n, \tilde R) &= \theta_R(\varepsilon) \,, \\
	\frac{\abs{R_n}}{n^m \abs{R} }
	&= 1
	+ \theta_R(\varepsilon)
	+ \theta_{R,\varepsilon} \bigl( \tfrac{1}{n} \bigr) \,, &
	\frac{\abs{\tilde R_n}}{n^m \abs{\tilde R}}
	&= 1
	+ \theta_R(\varepsilon)
	+ \theta_{R,\varepsilon} \bigl( \tfrac{1}{n} \bigr) \,, \\
	\frac{\abs{R}}{\abs{\tilde R}} &= 1 + \theta_R(\varepsilon) \,. &&
\end{align}
Let $h_R \in M(R)$ be an asymptotic height function
with $\Lip(h_R) \le 1 - c \varepsilon$.
Then,
\begin{equation} \llabel{e_concl} \begin{aligned}
	\hskip3em&\hskip-3em
	\Ent_{\tilde R_n} \bigl( B(\tilde R_n, h_R, \varepsilon), \omega \bigr)
	- C_\omega \theta_R(\varepsilon)
	- C_\omega \theta_{R,\varepsilon} \bigl( \tfrac{1}{n} \bigr) \\
	&\le \Ent_{R_n} \bigl( B(R_n, h_R, \varepsilon), \omega \bigr) \\
	&\le \Ent_{\tilde R_n} \bigl(
		B(\tilde R_n, h_R, \tfrac{c}{3} \varepsilon^2), \omega \bigr)
	+ C_\omega \theta_R(\varepsilon)
	+ C_\omega \theta_{R,\varepsilon} \bigl( \tfrac{1}{n} \bigr) \,.
\end{aligned} \end{equation}
\end{lemma}

\begin{proof}
We prove the two inequalities in~\leqref{e_concl} separately.
For the first inequality, observe that the map
\begin{equation} \begin{aligned}
	B(R_n, h_R, \varepsilon) &\to B(\tilde R_n, h_R, \varepsilon) \\
	h_R &\mapsto h_R|_{\tilde R}
\end{aligned} \end{equation}
is not generally an injection, but it is at most
$(2^{\abs{R_n \setminus \tilde R_n}})$-to-$1$
(by the graph homomorphism property and connectedness of $R_n$).
For any $h_{R_n} \in B(R_n, h_R, \varepsilon)$,
\[
	H_{R_n,\omega}(h_{R_n})
	\le H_{\tilde R_n,\omega}(h_{R_n}|_{\tilde R_n})
		+ C_\omega \abs{R_n \setminus \tilde R_n} \,,
\]
so
\[
	Z_\omega(B(R_n, h_R, \varepsilon))
	\le 2^{\abs{R_n \setminus \tilde R_n}}
		Z_\omega(B(\tilde R_n, h_R, \varepsilon))
		\exp \bigl( C_\omega \abs{R_n \setminus \tilde R_n} \bigr)
\]
and
\begin{equation} \begin{aligned}
	\hskip3em&\hskip-3em
	\Ent_{R_n} \bigl( B(R_n, h_R, \varepsilon), \omega \bigr) \\
	&\ge \frac{\abs{\tilde R_n}}{\abs{R_n}}
	\Ent_{\tilde R_n} \bigl( B(\tilde R_n, h_R, \varepsilon),
		\omega \bigr) \\
	&\qquad - \log(2) \, \frac{\abs{R_n \setminus \tilde R_n}}{\abs{R_n}}
	- C_\omega \abs{R_n \setminus \tilde R_n} \\
	&= \Ent_{\tilde R_n} \bigl( B(\tilde R_n, h_R, \varepsilon),
		\omega \bigr)
	- C_\omega \theta_R(\varepsilon)
	- C_\omega \theta_{R,\varepsilon} \bigl( \tfrac{1}{n} \bigr) \,.
\end{aligned} \end{equation}

To prove the second inequality in~\leqref{e_concl},
we first note that there exists an injection from
$B(\tilde R_n, h_R, \tfrac{c}{3} \varepsilon^2)$
into $B(R_n, h_R, \varepsilon)$.
A height function
$h_{\tilde R_n} \in B(\tilde R_n, h_R, \tfrac{c}{3} \varepsilon^2)$
is extended to $h_{R_n} \in B(R_n, h_R, \varepsilon)$
in such a way that
$\abs[\big]{ h_{R_n}(z) - n h_R \bigl( \tfrac{1}{n} z \bigr) } \le 1$
when $z$ is in $R_n$ and sufficiently far away from $\tilde R_n$;
the parameter value $\tfrac{c}{3} \varepsilon^2$ is chosen
so that such an extension is admissible by the Kirszbraun theorem.
For details, see the proof of \cite[Lemma 37]{KrMeTa19}.
For this injection $h_{\tilde R_n} \mapsto h_{R_n}$,
\[
	H_{\tilde R_n,\omega}(h_{\tilde R_n})
	\le H_{R_n,\omega}(h_{R_n})
		+ C_\omega \abs{R_n \setminus \tilde R_n} \,,
\]
so
\[
	Z_\omega \bigl( B(\tilde R_n, h_R, \tfrac{c}{3} \varepsilon^2) \bigr)
	\le Z_\omega \bigl( B(R_n, h_R, \varepsilon) \bigr)
	\exp \bigl( C_\omega \abs{R_n \setminus \tilde R_n} \bigr)
\]
and
\begin{equation} \begin{aligned}
	\hskip3em&\hskip-3em
	\Ent_{\tilde R_n} \bigl( B(\tilde R_n, h_R, \tfrac{c}{3} \varepsilon),
		\omega \bigr) \\
	&\ge \frac{\abs{R_n}}{\abs{\tilde R_n}}
	\Ent_{R_n} \bigl( B(R_n, h_R, \varepsilon), \omega \bigr) \\
	&\qquad - \log(2) \, \frac{\abs{R_n \setminus \tilde R_n}}{\abs{R_n}}
	- C_\omega \abs{R_n \setminus \tilde R_n} \\
	&= \Ent_{R_n} \bigl( B(R_n, h_R, \varepsilon), \omega \bigr)
	- \theta_R(\varepsilon)
	- \theta_{R,\varepsilon} \bigl( \tfrac{1}{n} \bigr) \,.
\end{aligned} \end{equation}
\end{proof}

To prove the profile theorem, we reduce to the special case
of Lemma~\ref{lem_simpl_profthm},
where the domain is a collection of simplices
and the asymptotic height function is piecewise affine.
Before that, in order to apply Lemma~\ref{lem_approx_micro_ent_dom},
we reduce to the case where $h_R$ has Lipschitz constant strictly less than $1$.
Both reductions are simple applications of the robustness results above.

\stepcounter{locallabel}

\begin{proof}[Proof of the profile theorem (Theorem~\ref{th_profile})]
For the reader's convenience we recall the conclusion of the theorem
that we are about to prove, namely: 
\begin{equation} \llabel{e_concl} \begin{aligned}
	\limsup_{n \to \infty}
	\mathbb{P} \biggl(
		&\abs[\Big]{
			\Ent_{R_n} \bigl( B(R_n, h_R, \delta), \omega \bigr)
			- \Ent_\an(R, h_R) } \\
		&\qquad\ge \eta
		+ C_\omega \theta_{h_R}(\delta)
		+ C_\omega \theta_{h_R,\delta} \bigl( \tfrac{1}{n} \bigr)
	\biggr) = 0 \,.
\end{aligned} \end{equation}

For the first step of the proof,
we reduce from the case of an arbitrary asymptotic height function
$h_R \in M(R, h_{\partial R})$,
i.e.\ a continuous function $h_R: R \to \R$
with Lipschitz constant at most $1$ (with respect to the $\ell^1$ norm on $R$),
to an asymptotic height function with Lipschitz constant strictly less than $1$.
Indeed, let $c := (2 \diam_1 R)^{-1} \wedge 1$,
where $\diam_1 R$ denotes the diameter of $R$ under the $\ell^1$ norm.
By translation invariance of the random field $\omega$,
we assume that there exists $x_0 \in R$ with $h_R(x_0) = 0$.
Define
\[
	\tilde h_R := (1 - c \delta) h_R \,.
\]

We make the following observations. First,
\begin{equation} \llabel{e_lip}
	\Lip(\tilde h_R)
	\;=\; (1-c\delta) \Lip(h_R)
	\;\le\; 1 - c \delta \,.
\end{equation}
Second, for any $x \in R$,
\begin{equation} \llabel{e_h}
	\abs{h_R(x) - \tilde h_R(x)}
	\;\le\; c \delta \abs{h_R(x)}
	\;\le\; c \delta \abs{x - x_0}_1
	\;\le\; \tfrac{\delta}{2} \,.
\end{equation}
Third, for any $x \in R$,
\begin{equation} \llabel{e_dh}
	\abs{\nabla h_R(x) - \nabla \tilde h_R(x)} \le c \delta \,.
\end{equation}

Lemma~\ref{lem_approx_macro_ent},
together with~\leqref{e_dh} and the choice of constant $c = c(R)$,
yields
\begin{equation} \llabel{e_lip_macro}
	\Ent_{R,\an}(h_R)
	= \Ent_{R,\an}(\tilde h_R) + \theta_R(\delta) \,.
\end{equation}

Similarly, Lemma~\ref{lem_approx_micro_ent_ahf} and~\leqref{e_h} imply that
almost surely,
\begin{equation} \begin{aligned}
	\hskip3em&\hskip-3em
	\Ent_{R_n} \bigl( B(R_n, \tilde h_R, 2\delta), \omega \bigr) \\
	&\le\Ent_{R_n} \bigl( B(R_n, h_R, \delta), \omega \bigr) \\
	&\le \Ent_{R_n} \bigl( B(R_n, \tilde h_R, \tfrac{1}{2} \delta),
		\omega \bigr) \,.
\end{aligned} \end{equation}

Assume for the sake of the proof that~\leqref{e_concl} holds for $\tilde h_R$.
Then almost surely,
\begin{equation} \begin{aligned}
	\Ent_{R_n} \bigl( B(R_n, h_R, \delta), \omega \bigr)
	&\le \Ent_{R_n} \bigl(
		B(R_n, \tilde h_R, \tfrac{\delta}{2}), \omega \bigr) \\
	&\le \Ent_{R,\an} (\tilde h_R)
	\,+\, \eta
	\,+\, \theta_{\tilde h_R} \bigl( \tfrac{\delta}{2} \bigr)
	\,+\, \theta_{\tilde h_R,\delta/2} \bigl( \tfrac{1}{n} \bigr) \\
	&= \Ent_{R,\an} (h_R)
	\,+\, \eta
	\,+\, \theta_{h_R} ( \delta )
	\,+\, \theta_{h_R,\delta} \bigl( \tfrac{1}{n} \bigr) ,
\end{aligned} \end{equation}
where in the last line,
we combine the $\theta_R(\delta)$ term from~\leqref{e_lip_macro} together
with the $C_\omega \theta_{\tilde h_R}(\tfrac{\delta}{2})$ term above;
this is admissible since $C_\omega \ge 1$ by definition
(recall that $C_\omega := 1 \vee \sup_{e \in E(\Z)} \abs{\omega_e}$)
and since the various factors of $\tfrac{1}{2}$ do not affect the asymptotics.
The reverse inequality is similar,
and so we have reduced to the problem of proving~\leqref{e_concl}
with the added assumption that $\Lip(h_R) \le 1 - c\delta$
for $c = c(R) \in (0,1)$.
\medskip

We reduce further to the special case from Lemma~\ref{lem_simpl_profthm},
i.e.\ a piecewise affine asymptotic height function
defined on a collection of simplices.
First, we choose parameter values $\varepsilon = \varepsilon(\delta)$
and $\ell = \ell(\varepsilon,\delta)$ satisfying three criteria:
\begin{enumerate}
\item $\varepsilon \to 0$ as $\delta \to 0$,
\item $\delta = \varepsilon \ell$,
\item $\ell$ is sufficiently small so that
the simplicial Rademacher theorem (Theorem~\ref{th_rade}) applies.
\end{enumerate}

The choices of $\varepsilon$ and $\ell$
may be realized as follows, from~\cite{KrMeTa19}:
Choose a sequence $\varepsilon_k \searrow 0$ arbitrarily,
e.g.\ $\varepsilon_k = \tfrac{1}{k}$.
Let $\ell_k$ be the largest admissible $\ell$ value based on $\varepsilon_k$,
but not larger than $1$.
For any given $\delta$ choose the smallest $\varepsilon_k$
such that $\varepsilon_k \ell_k > \delta$; this ensures the first criterion.
Set $\varepsilon = \varepsilon_k$
and $\ell = \frac{\delta}{\varepsilon_k} \le \ell_k$;
this ensures the last two criteria.
\\

For the remainder of the argument, fix $\delta > 0$.
Let $\varepsilon$ and $\ell$ satisfy the above criteria,
and let $K \subseteq R \subset \R^m$ be a simplicial domain
and $h_K \in M(K)$ an asymptotic height function
satisfying the conclusions of
the simplicial Rademacher theorem (Theorem~\ref{th_rade}).
Since $\nabla h_K \approx \nabla h_R$
(cf.\ conclusion~\ref{lem_approx_tri_grad} of Theorem~\ref{th_rade})
and since the macroscopic entropy is robust (Lemma~\ref{lem_approx_macro_ent}),
\begin{equation} \llabel{e_simpl_first} \begin{aligned}
	\abs[\Big]{ \Ent_{R,\an}(h_R) - \Ent_{K,\an}(h_K) }
	\,\le\, \theta_R(\varepsilon)
	\,=\, \theta_R(\delta) \,,
\end{aligned} \end{equation}
where we use the fact that $\varepsilon \to 0$ as $\delta \to 0$
in order to replace $\varepsilon$ by $\delta$ in the $\theta$ error term.
\\

Similarly, by conclusions~\ref{lem_approx_tri_volume}
and~\ref{lem_approx_tri_approx} of Theorem~\ref{th_rade}
and the microscopic entropy robustness,
\begin{equation} \llabel{e_simpl_second} \begin{aligned}
	\hskip3em&\hskip-3em
	\Ent_{R_n} \bigl( B(R_n, h_R, \varepsilon \ell), \omega \bigr) \\
	&\overset{\mathclap{(Lemma~\ref{lem_approx_micro_ent_dom})}}{\le} \qquad
	\Ent_{K_n} \bigl(
		B(K_n, h_R|_K, \tfrac{c}{3} (\varepsilon \ell)^2), \omega \bigr)
		+ C_\omega \theta(\varepsilon)
		+ C_\omega \theta_\varepsilon \bigl( \tfrac{1}{n} \bigr) \\
	&\overset{\mathclap{(Lemma~\ref{lem_approx_micro_ent_ahf})}}{\le} \qquad
	\Ent_{K_n} \bigl(
		B(K_n, h_K, \tfrac{c}{6} (\varepsilon \ell)^2), \omega \bigr)
		+ C_\omega \theta(\varepsilon)
		+ C_\omega \theta_\varepsilon \bigl( \tfrac{1}{n} \bigr)
\end{aligned} \end{equation}
and
\begin{equation} \llabel{e_simpl_third} \begin{aligned}
	\hskip3em&\hskip-3em
	\Ent_{R_n} \bigl( B(R_n, h_R, \varepsilon \ell), \omega \bigr) \\
	&\overset{\mathclap{(Lemma~\ref{lem_approx_micro_ent_dom})}}{\ge} \qquad
	\Ent_{K_n} \bigl(
		B(K_n, h_R|_K, \varepsilon \ell), \omega \bigr)
		- C_\omega \theta(\varepsilon)
		- C_\omega \theta_\varepsilon \bigl( \tfrac{1}{n} \bigr) \\
	&\overset{\mathclap{(Lemma~\ref{lem_approx_micro_ent_ahf})}}{\ge} \qquad
	\Ent_{K_n} \bigl(
		B(K_n, h_K, \tfrac{1}{2} \varepsilon \ell), \omega \bigr)
		- C_\omega \theta(\varepsilon)
		- C_\omega \theta_\varepsilon \bigl( \tfrac{1}{n} \bigr) \,.
\end{aligned} \end{equation}

Combining \leqref{e_simpl_first}, \leqref{e_simpl_second},
\leqref{e_simpl_third}, and the special case of the profile theorem proved
in Lemma~\ref{lem_simpl_profthm} completes the proof.
\end{proof}

\section{Variational principle} \label{s_varprin}

In this section we prove the variational principle (Theorem~\ref{th_varprin}).
The proof follows the steps of the corresponding proof for the uniform case
in~\cite{KrMeTa19}.
The main difference and the step that needs attention is
that the deterministic convergence needs to be lifted
to a convergence in probability.
The two main inequalities in the proof follow from
first comparing the set of height functions
$M(R_n, h_{\partial R_n}, \delta)$
to the subset $B(R_n, h_R^\ast, \delta)$ for a well-chosen
asymptotic height function $h_R^\ast$,
and second from comparing to a superset
$\bigcup_{i=1}^k B(R_n, h_R^{(i)}, \delta_i)$
for a collection of asymptotic height functions $h_R^{(1)}, \dotsc, h_R^{(k)}$.
Especially in the second part of the argument,
some care is needed in regards to the asymptotic parameters.
In particular:
\begin{itemize}
\item The choice (and number) of height functions $h_R^{(i)}$
	depends on $\delta$,
\item the radii $\delta_i$ of the balls around these height functions
	depends on $\eta$,
\item the probability that the profile theorem fails (i.e.\ the probability that
	$\Ent_{R_n}(B(R_n,h_R^{(i)},\delta_i),\omega)$
	and $\Ent_{R,\an}(h_R^{(i)})$ differ by a large amount
	due to the exact configuration $\omega$ of the random potential)
	depends not just on the error tolerance $\eta$
	but also on the number of height functions $h_R^{(i)}$.
\end{itemize}

\begin{proof}[Proof of Theorem~\ref{th_varprin}]
Let~$\eta > 0$ and $p_\tmax > 0$.
First we will establish that
\begin{equation} \llabel{e_upper} \begin{aligned}
	\limsup_{\delta \to 0} \, \limsup_{n \to \infty} \, \mathbb{P} \Bigl(
	&\Ent_{R_n} \bigl( M(R_n, h_{\partial R_n}, \delta),
		\omega \bigr) \\
	&\qquad >
	\inf_{h_R \in M(R, h_{\partial R})} \Ent_{R,\an}(h_R) + \eta
	\Bigr) \le p_\tmax \,.
\end{aligned} \end{equation}

Choose $h^* \in M(R, h_{\partial R})$ such that
\begin{equation} \llabel{e_upper_inf}
	\Ent_{R,\an}( h_R^* )
	\le \inf_{h_R \in M(R, h_{\partial R})} \Ent_{R,\an}( h_R )
	+ \tfrac{\eta}{4} \,.
\end{equation}

For any $\delta > 0$ and $n \in \N$,
$B(R_n, h_R^*, \delta) \subseteq M(R_n, h_{\partial R_n}, \delta)$.
Hence almost surely,
\begin{equation} \llabel{e_upper_star}
	\Ent_{R_n} \bigl( M(R_n, h_{\partial R_n}, \delta), \omega \bigr)
	\le \Ent_{R_n} \bigl( B(R_n, h_R^*, \delta), \omega \bigr) \,.
\end{equation}

By the profile theorem (applied to $h_R^*$),
\begin{equation} \llabel{e_upper_prof} \begin{aligned}
	\hskip3em&\hskip-3em
	\P \Bigl( \abs[\big]{
		\Ent_{R_n} \big( B(R_n, h_R^*, \delta), \omega \bigr)
		- \Ent_{R,\an} ( h_R^* ) } \\
	&> \tfrac{\eta}{4} + C_\omega \theta_{h_R^*}(\delta)
		+ C_\omega \theta_{h_R^*,\delta} \bigl( \tfrac{1}{n} \bigr)
	\Bigr) \underset{\mathclap{n \to \infty}}{\:\to\:} 0 \,.
\end{aligned} \end{equation}

Let us spend a part of the available probability $p_\tmax$
to establish a bound on $C_\omega$.
Specifically, since $C_\omega \in L^1$, Markov's inequality implies that
\begin{equation} \llabel{e_c}
	\P \bigl( C_\omega > \tfrac{2 \norm{C_\omega}_1}{p_\tmax} \bigr)
	\le \tfrac{1}{2} p_\tmax \,.
\end{equation}

Therefore as long as $\delta$ is small enough
so that the $\theta_{h_R^*}(\delta)$ term is less than
$\tfrac{\eta}{4} \cdot \tfrac{p_\tmax}{2 \norm{C_\omega}_1}$,
and as long as $n$ is large enough that
the $\theta_{h_R^*,\delta}(\tfrac{1}{n})$ term is less than
$\tfrac{\eta}{4} \cdot \tfrac{p_\tmax}{2 \norm{C_\omega}_1}$
and the probability in~\leqref{e_upper_prof}
is less than $\tfrac{1}{2} p_\tmax$,
we have
\begin{equation} \llabel{e_upper_prob}
	\P \Bigl( \Ent_{R_n} \bigl( B(R_n, h_R^*, \delta), \omega)
		> \Ent_{R,\an} (h_R^*) + \tfrac{3\eta}{4} \Bigr)
	< p_\tmax \,.
\end{equation}

The first desired inequality~\leqref{e_upper}
follows immediately from~\leqref{e_upper_star},~\leqref{e_upper_prob},
and~\leqref{e_upper_inf}.
\medskip

Now we turn to the second half of the variational principle, namely:
\begin{equation} \llabel{e_lower} \begin{aligned}
	\limsup_{\delta \to 0} \, \limsup_{n \to \infty} \, \mathbb{P} \Bigl(
	&
		\Ent_{R_n} \bigl( M(R_n, h_{\partial R_n}, \delta),
		\omega \bigr) \\
	&\qquad <
		\inf_{h_R \in M(R, h_{\partial R})} \Ent_{R,\an}(h_R) - \eta
	\Bigr) \le p_\tmax \,.
\end{aligned} \end{equation}

In order to establish~\leqref{e_lower},
we overcount the set $M(R_n, h_{\partial R}, \delta)$
using compactness of the space of asymptotic height functions
$M(R, h_{\partial R}, \delta)$
(with respect to the topology of uniform convergence).
Indeed, choose asymptotic height functions $h_R^{(1)}, \dotsc, h_R^{(k)}$
such that
\begin{equation} \llabel{e_lower_incl}
	M(R, h_{\partial R}, \delta)
	\subset \bigcup_{i=1}^k B(R, h_R^{(i)}, \delta_i) \,,
\end{equation}
where the values $\delta_i > 0$ are such that the
$\theta_{h_R^{(i)}}(\delta_i)$ terms from the profile theorem
(Theorem~\ref{th_profile}) are each less than
$\tfrac{\eta}{4} \cdot \tfrac{p_\tmax}{2 \norm{C_\omega}_1}$.

As in the first part of the proof,
we restrict to the event
\[
	\Omega'
	:= \Bigl\{ C_\omega < \frac{2 \norm{C_\omega}_1}{p_\tmax} \Bigr\} \,,
\]
which has $\P(\Omega') \ge 1 - \tfrac{p_\tmax}{2}$.
Furthermore, we assume implicitly that $n$ is large enough that:
\begin{itemize}
\item each of the $\theta_{h_R^{(i)},\delta_i}(\tfrac{1}{n})$ terms
	from the profile theorem is less than
	$\tfrac{\eta}{4} \cdot \tfrac{p_\tmax}{2 \norm{C_\omega}_1}$, and
\item the exceptional events \[
		E_{i,n} := \Omega' \cap \Bigl\{ \abs[\Big]{
			\Ent_{R_n} \bigl( B(R_n, h_R^{(i)}, \delta_i),
				\omega \bigr)
			- \Ent_{R,\an} ( h_R^{(i)} ) }
		> \tfrac{3\eta}{4} \Bigr\}
	\]
	satisfy $\P(E_{i,n}) < \tfrac{p_\tmax}{2k}$ for $i=1, \dotsc, k$.
\end{itemize}

Then for sufficiently small $\delta$ and sufficiently large $n$,
the ``good'' event
\[
	\Omega_{\delta,n} := \Omega' \cap E_{1,n}^c \cap \dotsb \cap E_{k,n}^c
\]
satisfies $\P(\Omega_{\delta,n}) \ge 1 - p_\tmax$
and, for $\omega \in \Omega_{\delta,n}$,
\begin{equation} \llabel{e_lower_good}
	\abs[\Big]{ \Ent_{R_n} \bigl( B(R_n, h_R^{(i)}, \delta_i), \omega \bigr)
		- \Ent_{R,\an} ( h_R^{(i)} ) }
	\le \tfrac{3\eta}{4} \,.
\end{equation}

Assume in the sequel that $\omega \in \Omega_{\delta,n}$.
By the set inclusion~\leqref{e_lower_incl},
\begin{equation} \llabel{e_lower_logsum}
	\Ent_{R_n} \bigl( M(R_n, h_{\partial R_n}, \delta), \omega \bigr)
	\ge - \frac{1}{\abs{R_n}} \log \biggl(
		\sum_{i=1}^k Z_\omega \bigl( B(R_n, h_R^{(i)}, \delta_i) \bigr)
	\biggr) \,.
\end{equation}

To handle the sum inside the logarithm,
we compare each summand $Z_\omega( B(R_n, h_R^{(i)}, \delta_i))$
against $\inf_{h_R} \Ent_{R,\an} (h_R)$.
Indeed,
\begin{equation} \begin{aligned}
	\Ent_{R_n} \bigl( B(R_n, h_R^{(i)}, \delta), \omega \bigr)
	&\overset{\mathclap{\leqref{e_lower_good}}}{\ge}\,
	\Ent_{R,\an} \bigl( h_R^{(i)} \bigr) - \frac{3\eta}{4} \\
	&\ge \inf_{h_R \in M(R, h_{\partial R})} \Ent_{R,\an} (h_R)
		- \frac{3\eta}{4} \,,
\end{aligned} \end{equation}
and so
\begin{equation}
	Z_\omega \bigl( B(R_n, h_R^{(i)}, \delta_i) \bigr)
	\le \exp \Bigl[ \abs{R_n} \Bigl( - \inf_{h_R \in M(R, h_{\partial R})}
		\Ent_{R,\an} ( h_R ) + \tfrac{3\eta}{4} \Bigr) \Bigr]
\end{equation}
and
\begin{equation}
	\sum_{i=1}^k Z_\omega \bigl( B(R_n, h_R^{(i)}, \delta_i) \bigr)
	\le k \exp \Bigl[ \abs{R_n} \Bigl( - \inf_{h_R \in M(R, h_{\partial R})}
		\Ent_{R,\an} ( h_R ) + \tfrac{3\eta}{4} \Bigr) \Bigr] \,.
\end{equation}

Returning to~\leqref{e_lower_logsum}, this yields
\begin{equation} \begin{aligned}
	\hskip3em&\hskip-3em
	\Ent_{R_n} \bigl( M(R_n, h_{\partial R_n}, \delta), \omega \bigr) \\
	&\ge \inf_{h_R \in M(R, h_{\partial R})} \Ent_{R,\an} ( h_R )
	- \frac{\log k}{\abs{R_n}} - \frac{3\eta}{4} \,.
\end{aligned} \end{equation}

As long as $n$ is large enough
(depending on $k$, which in turn depends on $\delta$),
we have $\tfrac{\log k}{\abs{R_n}} < \tfrac{\eta}{4}$, and so
\[
	\Ent_{R_n} \bigl( M(R_n, h_{\partial R_n}, \delta), \omega \bigr)
	\ge \inf_{h_R \in M(R, h_{\partial R})} \Ent_{R,\an} ( h_R ) - \eta \,,
\]
for any $\omega \in \Omega_{\delta,n}$.
This establishes~\leqref{e_lower} and thereby proves the variational principle
(Theorem~\ref{th_varprin}).
\end{proof}

\section{Open problems} \label{s_open}

\begin{itemize}

\item A natural question is whether
in the profile theorem (Theorem~\ref{th_profile})
and the variational principle (Theorem~\ref{th_varprin}),
the mode of convergence can be improved
from convergence in probability to almost-sure convergence.
The obstacle to achieving almost-sure convergence via the method of proof above
is the shifted environments $\tau_{i,n} \omega$
in~\eqref{e_simpl_profile_upper_decomptrans}.
Without the shifts $\tau_{i,n}$,
almost sure convergence would follow from the ergodic theorem,
applied individually for each index $i$ with $n \to \infty$.
It is possible that the ergodic theorem can be modified
to account for such shifts,
or that another method of proof can be used to improve the convergence result.
\item The proofs in this article assume that the random potential $\omega$
is almost surely bounded, and simulations provide evidence that the model does
not homogenize for some distributions of $\omega$ that are unbounded.
We conjecture that the model fails to homogenize
when additionally to our Assumption~\ref{a_random_field},
$\sup_{e \in E(\Z)} \abs{\omega_e} = \infty$ almost surely.
Alternatively, find the correct conditions on $\omega$ that ensure homogenization.
\item As mentioned above, we prove that the local surface tension is convex
(cf.\ Lemma~\ref{p_convexity}).
This is sufficient to conclude that the infima in the variational principle
and large deviations principle are attained
(as long as the set of height functions $A$
in the large deviations principle~\eqref{e_intro_ldp_nopot} is closed).
It would be useful to prove that the local surface tension is, moreover,
strictly convex.
Indeed, if the local surface tension is strictly convex,
then it follows that the minimizing height function
in the variational principle~\eqref{e_intro_varprin_nopot}
is unique, and hence is a limit shape.
Many random surface models are known to have
a strictly convex local surface tension,
e.g.\ domino tilings \cite{CKP01} and SAP models \cite{She05}.
For other models it is known that the local surface tension
is not strictly convex, e.g.\ the asymmetric five vertex model
(a degenerate case of the six-vertex model) \cite{dGeKeWa18}.
\item Characterize the fluctuations
of the perturbed probability measure $\mu_\omega$.
This is likely a complex problem.
By analogy to the dimer model studied in \cite{KOS06},
we expect that fluctuations may exhibit different asymptotics
in different parts of the domain (even asymptotically away from the boundary),
and by analogy to the random bridge model of \cite{GP11},
we expect non-trivial influences from the random potential.
\item Simulations suggest that the arctic circle phenomenon is universal,
i.e.\ that the shape of the boundary between the frozen and non-frozen regions
does not depend on the realization of the random field
or on the statistics of the random field.
This universality may even extend to unbounded random fields;
cf.\ Figure~\ref{f_sim}.
A promising method for studying the arctic circle is the tangent method
described in~\cite{CoSp16}.
\item We conjecture that concentration of measure holds,
at least in an appropriate asymptotic sense,
e.g.\ with high probability in the realization of the random field $\omega$.
It might be possible prove concentration
by adapting the idea of the harmonic embedding and corrector
from the study of random walks in random environment,
as explained in e.g.~\cite{Bis11}.
\end{itemize}

\appendix

\section{Ergodic theorem}
\label{a_almost_superadd_ergodic_thm}

Ergodic theory is a rich field of modern mathematics
with an extensive literature.
This includes several variants of the
superadditive (or subadditive) multidimensional ergodic theorem,
such as~\cite{Smy76,Ngu79,AK81},
which all propose technically different definitions of superadditivity
in the multidimensional setting.
The definition of superadditivity in~\cite{AK81} is a close match
for our application
(i.e.\ establishing that the limit almost surely exists
in our definition the quenched local surface tension).
However we actually need a version of the ergodic theorem
with weaker hypotheses, to allow for asymptotically negligible errors
in the superadditivity inequality~\eqref{e_almost_superadd}
and in the translation property~\eqref{e_almost_transinvar}.
These differences are not major or novel,
but neither are they so trivial that we are comfortable with omitting the proof
of the ergodic theorem under these weaker hypotheses.
At the time of writing we have not been able to find
this version of the ergodic theorem (or a stronger version) in the literature,
so we include a proof here.
The proof follows~\cite{AK81} closely;
for each step in the argument below,
we cite the corresponding step in~\cite{AK81}.

\begin{definition}[Boxes in $\Z^m$]
For $n \in \N$, let $S_n$ denote the box
\[
	S_n := [0,n)^m \cap \Z^m \,.
\]

Let~$\mathcal{B}$ denote the set of boxes
\begin{align*}
	\mathcal{B}
	&=
	\bigl\{
	( [a_1, b_1) \times \dotsb \times [a_m, b_m) ) \cap \Z^m
	\,\big| \\
	& \qquad \text{$a_i < b_i$ for all $i$, where $a_i, b_i \in \Z$}
	\bigr\},
\end{align*}
and for $k \in \N$, let $\mathcal{B}_k$ denote the set of boxes
\begin{align*}
	\mathcal{B}_k
	&=
	\bigl\{
	( [a_1, b_1) \times \dotsb \times [a_m, b_m) ) \cap \Z^m \in \mathcal{B}
	\,\big| \\
	& \qquad \text{all $a_i$ and $b_i$ are divisible by $k$}
	\bigr\}.
\end{align*}
\end{definition}

\begin{lemma}[A covering lemma; cf.\ {\cite[Lemma 3.1]{AK81}}]
\label{p_ak81_covering_lemma}
Let~$Z$ be a finite subset of~$\Z^m$.
For each~$z \in Z$ let~$n(z) \ge 1$ be an integer.
Then there is a set~$Z' \subseteq Z$
such that~$\{ z + S_{n(z)} \,|\, z \in Z' \}$ is a family of disjoint sets
and such that
$$ 3^m \sum_{z \in Z'} \abs{ S_{n(z)} } \ge \abs{Z}. $$
\end{lemma}

This is a modification of a common covering lemma due to Wiener.
The proof is standard.

\begin{theorem}[A maximal inequality; cf.\ {\cite[Theorem 3.2]{AK81}}]
\label{p_ak81_maximal_ineq}
For $\alpha > 0$, let $E_\alpha$ denote the event
$$
	E_\alpha :=
	\left\{ \limsup_{n \ge 1} \frac{1}{\abs{S_n}} F_{S_n} > \alpha \right\}
	\,.
$$

Then
\begin{equation} \label{e_ergo_maxl_concl}
	\mathbb{P}(E_\alpha) \le \frac{3^m}{\alpha/2} \, \gamma(F) .
\end{equation}
\end{theorem}

\begin{proof}
For~$N < M \in \N$, set
$$
	E_{N,M,\alpha} :=
	\left\{
		\sup_{N \le n \le M} \frac{1}{\abs{S_n}} F_{S_n} > \alpha
	\right\}.
$$

Clearly $E_\alpha = \cap_{N > 0} \cup_{M > N} E_{N,M,\alpha}$,
so it suffices to prove that
\begin{equation} \llabel{e_suff}
	\mathbb{P}(E_{N,M,\alpha})
	\le \frac{3^m}{\alpha/2} \, \gamma(F) + o(N) .
\end{equation}

Fix for now a larger integer~$K > M$. We will soon take~$K \to \infty$.
But first, consider a single~$\omega \in \Omega$.
Define the set~$Z$ as follows:
\begin{equation}
	Z = Z(\omega)
	:= \{z \in S_{K-M} \,|\, \tau_z \omega \in E_{N,M,\alpha} \}.
\end{equation}
We make two claims about $Z$:
first, that
$\frac{1}{\abs{S_K}} \E \abs{Z}$
is less than or equal to the right hand side of~\leqref{e_suff}
in the limit (see~\leqref{e_first_z_ineq} for the precise inequality),
and second,
that $\frac{1}{\abs{S_K}} \E \abs{Z} \ge \mathbb{P}(E_{N,M,\alpha})$
in the limit (see~\leqref{e_second_z_ineq} for the precise inequality).
After establishing these two claims, the result will follow quickly.
\\

Towards the first claim, consider any $z \in Z$.
There is an integer~$n(z)$ (implicitly depending on $\omega$)
such that~$N \le n(z) \le M$ and
\begin{equation}
	\frac{1}{\abs{S_{n(z)}}} \, F_{S_{n(z)}}(\tau_z \omega) > \alpha .
\end{equation}
By~\eqref{e_appx_trans_invar}, there exists $N_0 \in \N$
(independent of z and $\omega$) such that, whenever $N \ge N_0$,
\begin{equation} \label{e_ak81_maximal_ineq_F_gt_S}
	\frac{1}{\abs{S_{n(z)}}} \, F_{z + S_{n(z)}}(\omega)
	> \frac{\alpha}{2} .
\end{equation}

Apply the covering lemma (Lemma~\ref{p_ak81_covering_lemma}),
to pick~$z_1, \dotsc, z_l \in Z$ (again, implicitly depending on $\omega$)
such that the boxes~$z_i + S_{n(z_i)}$ are disjoint
but~$3^m \sum_{i=1}^l \abs{ S_{n(z_i)} } \ge \abs{Z}$.
Combining this with~\eqref{e_ak81_maximal_ineq_F_gt_S}
we get
\begin{equation*}
	\abs{Z}
	\le 3^m \sum_{i=1}^l \abs{ S_{n(z_i)} }
	\le \frac{3^m}{\alpha/2} \sum_{i=1}^l F_{z_i + S_{n(z_i)}} ,
\end{equation*}
and since $F \ge 0$ is almost superadditive,
\begin{equation*}
	\abs{Z}
	\le \frac{3^m}{\alpha/2} \, F_{S_K}
	+ A(\omega) \sum_{i=1}^l \abs{ \partial S_{n(z_i)}} .
\end{equation*}

Let $\varepsilon(N) = \sup_{n \ge N} \frac{\abs{\partial S_n}}{\abs{S_n}}$.
Note that~$\varepsilon(N) \to 0$ as~$N \to \infty$.
Since the boxes $z_i + S_{n(z_i)}$ are disjoint and contained in $S_K$,
\begin{equation} \label{e_ergo_maxl_zupp}
	\abs{Z}
	\le \frac{3^m}{\alpha/2} \, F_{S_K}
	+ A(\omega) \, \varepsilon(N) \, \abs{S_K} .
\end{equation}

Taking expectations and dividing by $\abs{S_K}$ yields the first claim, namely
\begin{equation} \llabel{e_first_z_ineq}
	\frac{1}{\abs{S_K}} \E \abs{Z}
	\le \frac{3^m}{\alpha/2} \cdot \frac{\E[F_{S_K}]}{\abs{S_K}}
	+ \norm{A}_{L^1} \, \varepsilon(N) .
\end{equation}

\medskip

Towards the second claim,
observe that as random variables,
\begin{equation}
	\abs{Z}
	= \sum_{z \in S_{K-M}} \mathbf{1}_{E_{N,M,\alpha}} \circ \tau_z .
\end{equation}
By translation invariance of the measure $\P$ on the random potential
(cf.\ Assumption~\ref{a_random_field}),
\begin{equation} \begin{aligned}
	\mathbb{E} \abs{Z}
	&= \sum_{z \in S_{K-M}}
		\mathbb{P} ( \tau_z \omega \in E_{N,M,\alpha} ) \\
	&= \sum_{z \in S_{K-M}}
		\mathbb{P} ( \omega \in E_{N,M,\alpha} ) \\
	&= \abs{ S_{K-M} } \, \mathbb{P}(E_{N,M,\alpha}) .
\end{aligned} \end{equation}

In other words,
\begin{equation} \llabel{e_second_z_ineq}
	\frac{1}{\abs{S_K}} \E \abs{Z}
	\ge \frac{\abs{S_{K-M}}}{\abs{S_K}} \P(E_{N,M,\alpha}) .
\end{equation}
\medskip

Combining the two claims that were just established,
namely~\leqref{e_first_z_ineq} and~\leqref{e_second_z_ineq},
we have
\[
	\frac{\abs{S_{K-M}}}{\abs{S_K}} \, \P(E_{N,M,\alpha})
	\,\le\, \frac{3^m}{\alpha/2} \, \frac{\mathbb{E}[F_{S_K}]}{\abs{S_K}}
	\,+\, \norm{A}_{L^1} \, \varepsilon(N) .
\]

Send $K$ to infinity:
$$
	\mathbb{P}(E_{N,M,\alpha})
	\,\le\, \frac{3^m}{\alpha/2} \, \gamma(F)
	\,+\, \norm{A}_{L^1} \, \varepsilon(N).
$$
This proves the desired inequality~\leqref{e_suff}
and completes the proof of Theorem~\ref{p_ak81_maximal_ineq}.
\end{proof}

\begin{lemma}[Convergence of expectations; cf.\ {\cite[Lemma 3.4]{AK81}}]
\label{p_ak81_convergence_of_expectations}
\begin{equation} \label{e_ergo_lim_concl}
	\gamma(F)
	= \lim_{n \to \infty} \frac{1}{\abs{S_n}} \, \mathbb{E}[F_{S_n}].
\end{equation}

Moreover, if~$H = (H_B)_{B \in \mathcal{B}_k}$ is almost superadditive but
defined only on boxes in~$\mathcal{B}_k$,
the same equality holds (except that both in the definition of~$\gamma(H)$
and in the right-hand side above, we only consider values of~$n$ that are
divisible by~$k$ as we take~$n \to \infty$).
\end{lemma}

\begin{proof}
By definition
$\gamma = \limsup_{n \to \infty} \frac{1}{\abs{S_n}} \mathbb{E}[F_{S_n}]$,
so it suffices to show that
$\liminf_{n \to \infty} \frac{1}{\abs{S_n}} \mathbb{E}[F_{S_n}] \ge \gamma$.
Let $k \in \N$.
For~$n \ge k$, we can subdivide the large box~$S_n$
into $r \ge 1$ translates of~$S_k$ and $s \ge 0$ translates of $S_1$,
say $S_n = \bigcup_{i=1}^r (u_i + S_k) \cup \bigcup_{j=1}^s (v_j + S_1)$.
By the superadditivity property~\eqref{e_appx_almost_superadd},
$$
	F_{S_n}
	\ge
	\sum_{i=1}^r F_{u_i+S_k} + \sum_{j=1}^s F_{v_j+S_1}
	- A \bigl( r \abs{\partial S_k} + s \abs{\partial S_1} \bigr).
$$
Taking expectations and dividing by~$\abs{S_n}$, we have
\begin{align*}
	\frac{1}{\abs{S_n}}\mathbb{E}[F_{S_n}]
	&\ge
	\frac{r}{\abs{S_n}} \mathbb{E}[F_{S_k}]
	+ \frac{s}{\abs{S_n}} \mathbb{E}[F_{S_1}] \\
	& \qquad
	- \mathbb{E}[A] \biggl(
		\frac{r \abs{\partial S_k}}{\abs{S_n}}
		+ \frac{s \abs{\partial S_1}}{\abs{S_n}}
	\biggr) \\
	&\qquad - \sup_{z \in \Z^m} \frac{1}{\abs{S_n}} \bigl(
		r \mathbb{E} \abs{F_{S_k} - F_{z+S_k}}
		+ s \mathbb{E} \abs{F_{S_1} - F_{z+S_1}} \bigr) \\
	&= \frac{1}{\abs{S_k}} \mathbb{E}[F_{S_k}] - o(n) - o(k) .
\end{align*}

Thus for every $k \ge 1$,
\[
	\liminf_{n \to \infty} \frac{1}{\abs{S_n}} \E[F_{S_n}]
	\ge \frac{1}{\abs{S_k}} \E[F_{S_k}] - o(k) ,
\]
and \eqref{e_ergo_lim_concl} follows by taking $k \to \infty$.

\medskip

Let us deal quickly with the case where the almost superadditive
process~$H = (H_B)_{B \in \mathcal{B}_k}$ is defined only on boxes
in~$\mathcal{B}_k$, i.e.\ only on boxes whose vertices lie on points of~$\Z^m$
whose every coordinate is divisible by~$k$.
We may define a process~$F = (F_B)_{B \in \mathcal{B}}$ by scaling,
i.e.\ $F_B = \frac{1}{\abs{S_k}} H_{kB}$,
where $kB = \{ ku \,|\, u \in B \}$ is the $k$-fold rescaling of $B$.
Then
$$
	\frac{1}{\abs{S_n}} F_{S_n}
	\, = \,
	\frac{1}{\abs{S_n} \abs{S_k}} H_{S_{kn}}
	\, = \,
	\frac{1}{\abs{S_{kn}}} F_{S_{kn}}
$$
so that $\gamma(F) = \gamma(H)$,
and the result just proven for~$F$ also carries over
(via linearity of the limit) to~$H$.
\end{proof}

\begin{theorem}[Ergodic theorem for almost superadditive random families]
\label{th_app_superadd_erg_thm}
Let~$(\Omega, \mathcal{F}, \mathbb{P})$ be a probability space,
let~$\tau = (\tau_u)_{u \in \Z^m}$ be a family of measure-preserving
transformations on~$\Omega$,
and let~$F = (F_B)_{B \in \mathcal{B}}$ be a family of~$L^1$ random variables
satisfying the following three conditions:
\begin{itemize}
	\item
	$F$ is almost superadditive, i.e.
	\begin{equation} \label{e_appx_almost_superadd}
		F_B \ge \sum_{i=1}^n F_{B_i} - A \sum_{i=1}^n \abs{\partial B_i}
		\quad \as \,,
	\end{equation}
	where~$A = A(\omega): \Omega \to [0, \infty)$
	is an~$L^1$ random variable.
	\item
	For all~$u \in \Z^m$,
	\begin{equation} \label{e_appx_trans_invar}
		\lim_{n \to \infty} \sup_{u \in \Z^m}
			\frac{1}{\abs{S_n}} \norm[\Big]{
				F_{u+S_n} - F_{S_n} \circ \tau_u}_{
				L^\infty(\omega)}
		= 0 \,,
	\end{equation}
	where~$u+B = \{ u+x \,|\, x \in B \}$ is the translation of~$B$ by~$u$.
	\item
	The quantity~$
		\gamma(F)
		= \limsup_{n \to \infty} \frac{1}{\abs{S_n}} \,
			\mathbb{E}[F_{S_n}]
	$ is finite.
\end{itemize}
Then the limit~$\lim_{n \to \infty} \frac{1}{\abs{S_n}} \, F_{S_n}$
exists almost surely and in~$L^1$.
If moreover~$\{\tau_u\}_{u \in \Z^m}$ is ergodic, then the limit is
\begin{equation} \label{e_appx_strong_superadd_erg_limit}
	\lim_{n \to \infty} \frac{1}{\abs{S_n}} \, F_{S_n}
	=
	\gamma(F).
\end{equation}
\end{theorem}

\begin{proof}[Proof of Theorem~\ref{th_app_superadd_erg_thm}]
The proof is in four steps.

\medskip

\paragraph{Step 1 (Reduction to~$F \ge 0$)}
Consider the to the additive process
$$
	G_B(\omega)
	:= \sum_{u \in B} F_{u+S_1}(\omega) - A(\omega) \, \abs{B}  .
$$
By the superadditivity property~\eqref{e_appx_almost_superadd},
$F' = F - G \ge 0$.
The desired convergence result is known for additive processes,
so it suffices to prove that $\frac{1}{\abs{S_n}} F'_{S_n}$
converges almost surely.
So, from this point on we shall assume that the process~$F$ is non-negative.

\medskip

\paragraph{Step 2 (Alternate rates of convergence)}
Let~$\overline{f} = \overline{f}(\omega)$
and~$\underline{f} = \underline{f}(\omega)$ denote respectively the
pointwise~$\limsup$ and~$\liminf$ of~$\frac{1}{\abs{S_n}} F_{S_n}$.
We shall show that, for~$m$ fixed, these two functions are also the
pointwise~$\limsup$ and~$\liminf$ of~$\frac{1}{\abs{S_{km}}} F_{S_{km}}$
as~$k \to \infty$.

For convenience, we write~$\overline{f}^{(m)}$ for the pointwise~$\limsup$ of
the sequence~$\frac{1}{\abs{S_{km}}} F_{S_{km}}$ as~$k \to \infty$.
Clearly~$\overline{f}^{(m)} \le \overline{f}$.
We must prove the opposite inequality.
Consider first any two boxes~$B \subseteq B'$.
Since~$F$ is almost superadditive and non-negative,
we have~$F_{B'} \ge F_B - O(\abs{B'})$.
In particular, when $k = \ceil{\tfrac{n}{m}}$,
$$
	\frac{1}{\abs{S_n}} F_{S_{m \ceil{n/m}}}
	\ge \frac{1}{\abs{S_n}} F_{S_n}
		- \frac{O( \abs{ S_{m \ceil{n/m}} } )}{\abs{S_n}} \,.
$$
Since $\abs{ S_{m \ceil{n/m}} } / \abs{S_n} \to 1$,
the left-hand side converges to $\overline{f}^{(m)}$ as $n \to \infty$,
and the right-hand side converges to $\overline f$.
The corresponding result for~$\underline{f}$ is proved similarly.

\medskip

\paragraph{Step 3 (Approximating~$F$)}
Fix~$\alpha > 0$. Let~$
	E = \{ \omega : \overline{f}(\omega) - \underline{f}(\omega) > \alpha \}
$.
In order to show that~$\mathbb{P}(E) = 0$, let~$\varepsilon > 0$.
By Lemma~\ref{p_ak81_convergence_of_expectations}, there exist~$k$ arbitrarily
large such that~$
	\frac{1}{\abs{S_k}} \mathbb{E}[F_{S_k}] >
	\gamma - \frac{\varepsilon}{2}
$.
Define an additive family~$H$ on~$\mathcal{B}_k$
(which, we recall, is the set of boxes whose vertices all lie in the
sub-lattice~$k \Z^m \subset \Z^m$) by
$$
	H_B =
	\sum_{u \in B \cap k\Z^m} F_{u+S_k}
	- A \, \abs{B \cap k\Z^m} \, \abs{\partial S_k}.
$$
By almost superadditivity~\eqref{e_appx_almost_superadd},
$$
	F_B
	\ge
	\sum_{u \in B \cap k \Z^m} F_{u + S_k}
	- \abs{A} \sum_{u \in B \cap k \Z^m} \abs{ \partial (u + S_k) }
	= H_B .
$$

Now let~$F' = F - H$, so that~$F'$ is a non-negative random family defined
on~$\mathcal{B}_k$.
It holds that
\begin{align}
	\overline{f}(\omega) - \underline{f}(\omega)
	&:= \limsup_{n \to \infty} \frac{1}{\abs{S_{kn}}} F_{S_{kn}}
	- \liminf_{n \to \infty} \frac{1}{\abs{S_{kn}}} F_{S_{kn}} \\
	&\overset{(*)}=
	\limsup_{n \to \infty} \frac{1}{\abs{S_{kn}}} F'_{S_{kn}}
	- \liminf_{n \to \infty} \frac{1}{\abs{S_{kn}}} F'_{S_{kn}} \\
	\label{e_ergo_pf_ffp}
	&\overset{(**)}\le
	\sup_{n \to \infty} \frac{1}{\abs{S_{kn}}} F'_{S_{kn}} .
\end{align}
In particular, $(*)$ holds because $H$ is additive,
so it converges pointwise almost surely,
and $(**)$ holds because $F' \ge 0$.

\medskip

Next, we compute~$\gamma(H)$ and $\gamma(F')$.
Applying Lemma~\ref{p_ak81_convergence_of_expectations}:
\begin{align*}
	\gamma(H)
	&=
	\lim_{n \to \infty} \left(
		\frac{1}{\abs{S_{kn}}} \mathbb{E}[H_{S_{kn}}]
	\right)
	\\
	&=
	\lim_{n \to \infty} \left(
		\frac{\sum_{u \in S_{kn} \cap k \Z^m}\mathbb{E}[F_{u+S_k}]}
			{k^m n^m}
		-
		\frac{
			\mathbb{E}[A] \,
			\abs{ S_{kn} \cap k \Z^m } \,
			\abs{ \partial S_k }
		}{
			k^m n^m
		}
	\right)
	\\
	&=
	\lim_{n \to \infty} \left(
		\frac{1}{\abs{S_k}} \mathbb{E}[F_{S_k}]
		-
		\mathbb{E}[A]  \frac{\abs{\partial S_k}}{\abs{S_k}}
		- o(k)
	\right).
\end{align*}
Note that~$n$ no longer appears in the final expression.
Taking $k \to \infty$,
we conclude that~$\gamma(H) > \gamma(F) - \varepsilon$.
Additionally from Lemma~\ref{p_ak81_convergence_of_expectations},
we can write~$\gamma(F')$ as a limit.
Importantly,~$\gamma$ is linear,
so~$\gamma(F') = \gamma(F) - \gamma(H) < \varepsilon$.

\medskip

By~\eqref{e_ergo_pf_ffp}
the event~$E := \{ \overline{f} - \underline{f} > \alpha \}$
is contained in~%
$\{ \sup_{n \ge 1} \frac{1}{\abs{S_{kn}}} F'_{S_{kn}} > \alpha \}$.
By Lemma~\ref{p_ak81_maximal_ineq},
$$
	\mathbb{P}(E)
	\le \frac{3^m \gamma(F')}{\alpha/2}
	\le \frac{3^m \varepsilon}{\alpha/2}.
$$
Taking~$\varepsilon \to 0$, we see that~$\mathbb{P}(E) = 0$.
Since~$\alpha > 0$ was arbitrary,
we conclude that~$\overline{f} = \underline{f}$ almost surely,
and thus that~$\frac{1}{\abs{S_n}} F_{S_n}$ converges pointwise almost surely.
\end{proof}

\section*{Acknowledgment}
The authors want to thank Tim Austin, Nathana\"el Berestycki, Marek Biskup, Antoine Gloria, Michelle Ledoux, Thomas Liggett, Igor Pak, Larent Saloff-Coste, and Tianyi Zheng for the many discussions on this topic.
\bibliographystyle{alpha}
\bibliography{bib}

\end{document}